\documentclass[12pt]{article}
\usepackage{amsmath}
\usepackage{amsthm}
\usepackage{amssymb}
\usepackage{amsbsy}
\usepackage{amsfonts}
\usepackage{amstext}
\usepackage{amscd}
 \usepackage{mathrsfs}
 \usepackage{stmaryrd}
\usepackage[dvipdfmx]{graphicx,psfrag}
\usepackage{enumitem,color,comment}
\usepackage{mathtools}
\usepackage{tikz}
\usepackage{geometry}

\numberwithin{equation}{section}
\theoremstyle{plain}
\newtheorem{thm}{Theorem}[section]
\newtheorem{prop}[thm]{Proposition}
\newtheorem{cor}[thm]{Corollary}
\newtheorem{lem}[thm]{Lemma}
\theoremstyle{definition}
\newtheorem{exa}[thm]{Example}

\newtheorem{rem}[thm]{Remark}
\newtheorem{defi}[thm]{Definition}
\newtheorem{Notation}[thm]{Notation}

\newcommand{\real}{\mathbb{R}}
\newcommand{\comp}{\mathbb{C}}
\newcommand{\N}{\mathbb{N}}
\newcommand\cA{\mathcal{A}}
\newcommand\cB{\mathcal{B}}
\newcommand\cS{\mathcal{S}}   % spreadability system 
\newcommand\cU{\mathcal{U}}   

\newcommand\alg[1]{{\rm alg}(#1)}
\newcommand\wh[1]{\widehat{#1}}   
\newcommand\wt{\widetilde}   
\newcommand\vp{\varphi}     

\newcommand\id{\mathrm{id}}  
\newcommand\di{{\rm d}}   % for integration, e.g. dx, dt. 
\newcommand{\Der}{\left. \frac{\di}{\di t} \right|_{t=0}}      % derivative 
\newcommand{\pDer}{\left. \frac{\partial}{\partial t} \right|_{t=0}}    % partial derivative

\newcommand\pp{p}  % linear functional 
\newcommand{\Restr}[2]{{#1}{\restriction}_{#2}}   % restriction 
\newcommand\sstar{\mathop{\text{\Large$\star$} \hspace{0.1em} }   }      % Nonunital free product
\newcommand\sstat[1]{\mathop{\text{\Large$\star$}}_{#1} \hspace{0.1em}    }      % Nonunital free product
\newcommand\freeprod{ \mathop{\text{\Large$\ast$} \hspace{0.1em}  } }    % unital free product
\newcommand\freeproe[1]{ \mathop{\text{\Large$\ast$}   }_{#1} \hspace{0.1em} }   % unital free product 

\newcommand\vara{u} % variance for \varphi
\newcommand\varb{v}   % variance for \psi
\newcommand\varc{w}  % variance  for \theta
\newcommand\fin[1]{{#1}^\ast}    %  the set of sequences of elements of $I$ of finite length
\newcommand\bfi{{\mathbf i}}    % A sequence of elements of $I$ of finite length  
\newcommand\bfj{{\mathbf{j}}}     % A sequence of elements of $I$ of finite length 
\newcommand\bfk{{\mathbf{k}}}     % A sequence of elements of $I$ of finite length 

\renewcommand\AA{\mathscr{B}}      % A set of inner blocks
\newcommand\BB{\mathscr{C}}     % A set of inner blocks
\newcommand\nA{B}      %  A direct summand of free product pre-Hilbert space
\newcommand\nB{C}      % A direct summand of free product pre-Hilbert space
\newcommand\Bottom{{\rm Bot}}   % The set of bottoms of a sequence (without endpoints)
\newcommand\Peak{{\rm Peak}}  % The set of peaks of a sequence  (without endpoints)

\newcommand\AP{\alpha}    %Alpha product / indep / cumulants / convolution / ... 
\newcommand\BGP{\beta\gamma}    % Beta-gamma product / ...
\newcommand\BP{\beta}   %  Beta product /...
\newcommand\GP{\gamma}  %   Gamma product / ...

\newcommand\CF{{\rm CF}}   %C-free
\newcommand\CM{{\rm CM}}   %c-monotone
\newcommand\AM{{\rm AM}}  %antimonotone
\newcommand\CAM{{\rm CAM}}   %c-antimonotone
\newcommand\Free{{\rm F}}   %Free
\newcommand\Mon{{\rm M}}  %Monotone
\newcommand\Boole{{\rm B}}  % Boole

\newcommand\Des{{\rm Des}}   % Descents
\newcommand\Asc{{\rm Asc}}    % Ascents  
\newcommand\Lin{\mathcal{L}^\ast}    % Adjointable linear operators
\newcommand{\SP}{\mathcal{P}} 
\newcommand{\NC}{\mathcal{NC}}  
\newcommand{\IP}{\mathcal{I}}      % Interval partitions 
\newcommand{\Inn}{\mathscr{I}}     % Inner blocks 
\newcommand{\Out}{\mathscr{O}}    % Outer blocks
\newcommand\Cov{C}        %  The nearest cover 
\newcommand{\OSP}{\mathcal{OP}}     
\newcommand{\ONCP}{\mathcal{ONC}}

\newcommand{\MP}{\mathcal{M}}  
\newcommand{\AMP}{\mathcal{AM}}  
\newcommand{\CSP}{\mathcal{CP}} 
\newcommand{\fmp}{\mathop{\leftthreetimes\hspace{0.1em}}}    % \alpha-product 
\newcommand{\cf}[4]{\mathop{{#1}\hspace{0.1em}{}_{#2}\hspace{-0.3em}\ast_{#3}\hspace{-0.1em}{#4} }}     % c-free product  
\newcommand{\cm}[3]{\mathop{{#1}\trr_{#2}\hspace{-0.1em}{#3} }}      % c-monotone 
\newcommand{\cam}[3]{\mathop{{#1} \hspace{0.1em} {}_{#2}\hspace{-0.3em}  \trl{#3} }}    % c-antimonotone 
\newcommand{\trr}{\triangleright}     % monotone  product of states (also for additive convolutions)
\newcommand{\trl}{\triangleleft}       % anti-monotone (also for additive convolutions)
\newcommand{\cfconv}[4]{\mathop{{#1}\hspace{0.1em}{}_{#2}\hspace{-0.3em}\boxplus_{#3}\hspace{-0.1em}{#4} }}    %additive c-free convolution 
\newcommand{\subcfconv}[4]{\mathop{{#1}{}_{#2}\hspace{-0.1em}\boxplus_{#3}{#4} }}   % additive c-free convolution (for subscripts)
\newcommand{\mcfconv}[4]{\mathop{{#1}\hspace{0.1em}{}_{#2}\hspace{-0.3em}\boxtimes_{#3}\hspace{-0.1em}{#4} }}     % multiplicative c-free convolution 
\newcommand{\submcfconv}[4]{\mathop{{#1}{}_{#2}\hspace{-0.1em}\boxtimes_{#3}{#4} }}     % multiplicative c-free convolution (for subscripts)
\renewcommand{\subsetplus}{ \mathop{\hspace{0.07em}\text{\large$\sqsubset$}\hspace{-0.72em} \raisebox{0.15em}[0ex][0ex]{$\scriptstyle+$}  \hspace{0.3em} } }  % Additive alpha-convolution
\newcommand{\subsettimes}{  \mathop{\hspace{0.07em}\text{\large$\sqsubset$}\hspace{-0.72em} \raisebox{0.15em}[0ex][0ex]{$\scriptstyle\times$}  \hspace{0.3em} }  }  % Multiplicative alpha-convolution

\begin{document}   
%%%%%%%%%%%

{
%\large %%%%%% for booklet

%\title{An associative product for triplets of states generalizing free, monotone and Boolean products}
\title{A three-state independence in non-commutative probability}
%\title{New associative product for triplets of states generalizing free, monotone, antimonotone, Boolean, conditionally free and conditionally monotone products}
%Unification of  free, monotone, antimonotone and Boolean independences: an approach via triplet of states}
\author{Takahiro Hasebe\footnote{This work is supported by JSPS Grant-in-Aid for Young Scientists 19K14546 and Scientific Research (B) 18H01115. }}

\date{}

\maketitle

\begin{abstract} 
We define a new independence  in non-commutative probability, called $\AP$-freeness, with respect to a triplet of states. This concept unifies several independences in non-commutative probability, in particular,  free, monotone, antimonotone and Boolean ones as well as conditionally free, conditionally monotone and conditionally antimonotone independences.  
Moreover, the associative law of $\AP$-freeness is transferred to the other independences. As a consequence, $\AP$-free cumulants unify the cumulants for free, monotone, antimonotone and Boolean independences. 
The central limit theorem for $\AP$-freeness is computed. The limit distribution turns out to be a triplet of the Kesten distributions.  

\end{abstract}

\tableofcontents

%Keywords: Freeness; conditional freeness; monotone independence; Boolean independence; cumulants \\

%Mathematics Subject Classification: Primary 46L53, 46L54; secondary 06A07 

%%%%%%%%%%%%%%%%%%%%%%%%%%%%%%%%%%%%%%%%%%%%%%%%%%%%%%%%%%%%%%%%%%%%%%%%%%%%%%%%%%%%%%%%%%%%%%%%%%%%%%%%%%%%%%%%%%%%%%%%%%%%%%%%%%%%%%%%%%%%%%%%%%%%%%%%%%%%
\section{Introduction}

%%%%%%%%%%%%%%%%%%%%%%%%%%%%%%%%%%%%%%%%%%%%%%%%%%%%%%%%%%%%%%%%%%%%%
\subsection{Backgrounds}\label{subsec:backgound}

A building block of noncommutative probability is notions of independence for noncommutative random variables. A straightforward extension of the standard independence for (commutative) random variables to noncommutative ones is called tensor independence. It was Voiculescu who discovered another notion of independence called freeness or free independence \cite{V1}, which turned out to have numerous applications to operators algebras and random matrices. The introduction of freeness also motivated further search for other possible notions of independence. Indeed, more examples are then introduced; basic ones include Boolean independence \cite{Boz86,SW97}, monotone independence \cite{Mur2} and antimonotone independence (essentially equivalent to monotone independence).  

Attempts were made to formulate a general notion of independence and classify it \cite{B-S1,GHU,JL20,Leh,Mur3,Mur5, Sch95,Spe1,Var21}.  Among all, universal independence was introduced by Speicher \cite{Spe1}, and later, further investigated by Sch\"{u}rmann and Ben Ghorbal from a category theoretic viewpoint \cite{B-S1}. According to the classification theorem, there are exactly three universal independences: tensor, free and boolean ones. Notably, Muraki \cite{Mur3,Mur5}  formulated a more general notion ``natural independence'' and classified it into five ones; then   monotone independence and antimonotone independence appear in addition to the three universal independences.

On the other hand, Bo\.{z}ejko, Leinert and Speicher \cite{BLS,BS} formulated conditional freeness (or c-freeness for brevity) referring to a pair of states.   
By its definition, c-freeness unifies free and Boolean independences. In addition, Franz pointed out that monotone (and antimonotone) independence is also a special case of c-freeness \cite{Fra}. Along similar spirits, conditionally monotone (or c-monotone for brevity) independence was introduced by the author \cite{Has3} again referring to pairs of states, which unifies monotone and Boolean independences.  By reversing the order structure,  c-antimonotone independence can also be naturally defined. 

The main objective of this paper is to unify the above mentioned independences referring to a \emph{triplet of states}.  As preliminaries, we give detailed definitions of the above mentioned independences in Subsection \ref{subsec:pre} and explain the aim and motivations in Subsection \ref{sec:aim_motivation}. The main results  are briefly summarized in Subsection \ref{subsec:main}.

\subsection{Preliminaries: products of algebraic probability spaces}\label{subsec:pre}

In the present subsection we give precise definitions of several notions of independence mentioned in Subsection \ref{subsec:backgound} to provide sufficient preliminaries of this work. For later use, it is more convenient to work on products of (unital) algebraic probability spaces, which is a kind of canonical realizations of independences.  For the relationship between products of algebraic probability spaces and independences, see e.g.\ \cite{GHU}.

 Before giving definitions of key products, some supplementary definitions and notations are listed below.  
 
 \begin{Notation}\label{notation} Let $I$ be a set and $\{S_i\}_{i \in I}$ a family of sets. 
 \begin{enumerate}[label=\rm(\alph*)]
\item $[n]$ denotes the set $\{1,2,\dots, n\} \subseteq\N$ for all $n\in\N$,

% \item For $i_1,i_2,\dots, i_n \in I$, the symbol $i_1 \neq \cdots \neq i_n$ means that the neighboring elements are different, i.e.\ $i_k \ne i_{k+1}$ for all $1\le k\le n-1$.  

\item We set
\[
\fin{I} :=\{(i_1, \cdots, i_n): n \in \N, i_1,i_2,\dots, i_n \in I,  i_k \neq i_{k+1} ~\text{for all}~k \in [n-1] \}.
\] 
The length of an element $\bfi \in \fin{I}$ is denoted $|\bfi|$.

\item For $\bfi =(i_1,i_2,\dots, i_n) \in \fin{I}$, we denote 
\[
S_\bfi = S_{i_1} \times  S_{i_2} \times  \cdots \times S_{i_n}. 
\]

\item If $I$ is a toset (totally ordered set), then we denote by $I^{\rm op}$ the toset which is equal to $I$ as a set but with the reversed order, i.e.\ $i \le_{I^{\rm op}} j$ if and only if $i \ge_I j$. 

\item For an algebra $\cA$, its unitization is denoted $\wh \cA = \comp \oplus \cA$. For a linear functional $\vp$ on $\cA$, the unitized linear functional on $\wh \cA $ is denoted $\wh \vp$, i.e.\ $\wh \vp(\lambda + a):= \lambda + \vp(a)$ for all $\lambda \in \comp$ and $a \in \cA$.  
\end{enumerate}
\end{Notation}

\begin{rem} The field for algebras is always assumed to be $\comp$ in this paper.  
\end{rem}

In most part of this paper, the positivity structure of algebraic probability spaces is not crucial. When positivity is needed or is of interest, additional remarks will be provided. Our basic setup is as follows. 
 
\begin{defi} \label{def:basic}

A \emph{unital $m$-algebraic probability space} (resp.\ an \emph{$m$-algebraic probability space}) is a tuple $(\cA, \vp^1,\vp^2,\dots, \vp^m)$, where $\cA$ is a unital algebra (resp.\ an algebra) and $\vp^i$ are unital linear functionals on $\cA$ (resp.\ linear functionals on $\cA$).  When $m=1$ we call it a unital algebraic probability space (resp.\ an algebraic probability space). 
\end{defi}

The nonunital and unital free products of algebras will play key roles. The definition is summarized below. For further details see e.g.\ \cite{NS06,VDN92}. 

\begin{defi}
\begin{enumerate}[label=\rm(\roman*)]

\item The \emph{free product} of a family of algebras $\{\cA_i\}_{i\in I}$ is an algebra $\cA$ equipped with a family of homomorphisms $\{\iota_i\colon \cA_i \to \cA\}_{i\in I}$ such that  for any algebra $\cB$ and homomorphisms $f_i\colon \cA_i \to \cB, i\in I$,  there exists a unique homomorphism $f\colon\cA\to\cB$ such that $f_i = f \circ \iota_{i}$ for all $i\in I$.  This universality uniquely determines the algebra $\cA$ up to isomorphism. Unless otherwise specified, in this paper we call the following construction the free product: 
\[
\sstat{i\in I} \cA_i := \bigoplus_{n\in\N} \bigoplus_{\substack{(i_1,i_2,\dots, i_n) \in I^n \\ i_1 \ne \cdots \ne i_n} } \cA_{i_1} \otimes \cdots \otimes \cA_{i_n}
\]
with multiplication determined by, for $a_i\in A_{k_i}, b_j\in A_{\ell_j}, i\in[n]$, $j\in[m]$,  
\begin{align*}
    (a_1 \otimes \cdots \otimes a_n)(b_1\otimes\cdots \otimes b_m)&=
    \begin{cases}
        a_1 \otimes \cdots \otimes a_{n-1}\otimes (a_nb_1)\otimes b_2 \otimes \cdots \otimes b_m& \text{if }k_n=\ell_1,\\
        a_1 \otimes \cdots \otimes a_n\otimes b_1\otimes\cdots \otimes b_m & \text{if }k_n\neq \ell_1, 
    \end{cases} 
\end{align*} 
and with the canonical embeddings $\iota_i\colon \cA_i \xhookrightarrow{} \cA$.

\item\label{item:unital_free_product} The \emph{unital free product} of a family of unital algebras $\{\cA_i\}_{i\in I}$ is a unital algebra $\cA$ equipped with a family of unital homomorphisms $\{\iota_i\colon \cA_i \to \cA\}_{i\in I}$ such that for any unital algebra $\cB$ and unital homomorphisms $f_i\colon \cA_i \to \cB, i\in I$,  there exists a unique unital homomorphism $f\colon  \cA \to\cB$ such that $f_i = f \circ \iota_{i}$ for all $i\in I$. 
This universality uniquely determines the unital algebra $\cA$ up to isomorphism. A useful construction of the unital free product is as follows. For each $i\in I$ we take and fix a  direct sum decomposition (as vector spaces) $\cA_i = \comp1_{\cA_i} \oplus \cA_i^0$. Then the vector space 
\begin{equation}\label{eq:free_product}
\freeproe{i\in I} \cA_i:=  \comp1 \oplus \bigoplus_{n\in\N} \bigoplus_{\substack{(i_1,i_2,\dots, i_n) \in I^n \\ i_1 \ne \cdots \ne i_n} } \cA_{i_1}^0 \otimes \cdots \otimes \cA_{i_n}^0, 
\end{equation}
together with a suitable multiplication and the canonical embeddings, satisfies the universality. In Subsection \ref{subsec:alternative} a more general construction will be given. 
%where $\cI$ is the two-sided ideal generated by $\{1 - 1_{\cA_i}\}_{i\in I}$. The unital free product has

 \end{enumerate}
\end{defi} 

\begin{rem} 
There is a canonical unital algebra isomorphism $
\freeproe{i\in I}\wh{\cA}_i  \simeq \wh{\sstat{i\in I}\cA_i}.$ 
%induced by the map $f\colon \wh{\sstat{i\in I}\wh{\cA}_i } \to \wh{\sstat{i\in I}\cA_i}$ which deletes the unit elements of $\wh{\cA}_i$'s from tensor components. For example, if $I=\{1,2,3\}$, $a \in \cA_2, b,c \in \cA_1$, then the definition $f$ includes the rule   
%\[
%f(a \otimes 1_{\wh{\cA}_3} \otimes b \otimes 1_{\wh{\cA}_2} \otimes c) := a \otimes (bc). 
%\]
\end{rem}

\begin{rem}  When discussing the positivity structure on unital algebraic probability spaces,  we will be concerned with the following setup: 

\begin{enumerate}[label=\rm(\alph*)]
\item A unital linear functional $\vp$ on a unital $\ast$-algebra $\cA$ is called a \emph{state} if $\vp$ is positive, i.e.\ $\vp(a^*a) \ge0$ for all $a\in\cA$;    

\item A \emph{unital $m$-$\ast$-algebraic probability space} is a tuple $(\cA, \vp^1,\vp^2,\dots, \vp^m)$, where $\cA$ is a unital $\ast$-algebra and $\vp^i$ are states on $\cA$.  When $m=1$ we call it a unital $\ast$-algebraic probability space;  

\item For a family of unital $\ast$-algebras $\{\cA_i\}_{i\in I}$ the free product $\freeproe{i\in I} \cA_i$ can be made a $\ast$-algebra with the canonical involution  
\[
(a_1 \otimes \cdots \otimes a_n)^* := a_n^* \otimes \cdots \otimes a_1^*.   
\]

\end{enumerate}
We can also consider a positivity structure on (not necessarily unital) $m$-algebraic probability spaces, which is of less importance in this paper. In that case, for algebras we simply require $\cA$ to be $\ast$-algebras. For linear functionals we require $\vp^i$ to be \emph{restricted states}, i.e.\  their unital extensions $\wh{\vp}^i$ to the unitized $\ast$-algebra $\wh{\cA}= \comp \oplus \cA$ to be positive. For further information on restricted states the reader is referred to \cite{GHU} or \cite{Lac15}. 
\end{rem}

\begin{defi}[\cite{V1}] Let $I$ be an index set with $\#I \ge2$ and $\{(\cA_i, \vp_i)\}_{i \in I}$ be a family of  unital algebraic probability spaces. 
The \emph{free product} of $\{(\cA_i, \vp_i)\}_{i\in I}$ is a unital algebraic probability space $(\cA, \vp)$, denoted $\freeproe{i\in I} (\cA_i, \vp_i)$, where $\cA = \freeproe{i\in I}\cA_i$ is the unital free product of $\{\cA_i\}_{i\in I}$ and $\vp$ is a linear functional on $\cA$ determined by the following rule: for every $\bfi=(i_1,i_2,\dots, i_n) \in \fin{I}$ and every $(a_1,a_2,\dots, a_n) \in  \cA_{\bfi}$ with $\vp_{i_k}(a_k)=0$ for all $k \in [n]$, we have
 \[
\vp(a_1 a_2 \dotsm a_n) = 0. 
 \]
For brevity, we sometimes omit the algebras, writing $\varphi = \freeproe{i\in I} \vp_i$.
 \end{defi}

\begin{defi}[\cite{Boz86,SW97}] \label{def:BP} Let $I$ be an index set with $\#I \ge2$ and $\{(\cA_i, \vp_i)\}_{i \in I}$ be a family of algebraic probability spaces. 
The \emph{Boolean product} of $\{(\cA_i, \vp_i)\}_{i\in I}$ is an algebraic probability space, denoted $(\cA, \vp)=\diamond_{i\in I}(\cA_i, \vp_i)$, where $\cA = \sstat{i\in I} \cA_i$ is the  free product of $\{\cA_i\}_{i\in I}$ and $\vp$ is a linear functional on $\cA$ determined by the following rule: for every $\bfi=(i_1,i_2,\dots, i_n) \in \fin{I}$ and every $(a_1,a_2,\dots, a_n) \in  \cA_{\bfi}$, we have
 \[
\vp(a_1 a_2 \dotsm a_n) = \prod_{k=1}^n \vp_{i_k}(a_k). 
 \]
  We also employ the notation $\varphi = \diamond_{i\in I} \vp_i$. 
 \end{defi}
 
\begin{defi}[\cite{Mur3,Mur2}]  \label{def:MP-AMP} Let $I$ be a toset with $\#I \ge2$ and $\{(\cA_i, \vp_i)\}_{i \in I}$ be a family of algebraic probability spaces. 

\begin{enumerate}[label=\rm(\roman*)]

\item\label{item:MP} The \emph{monotone product} of $\{(\cA_i, \vp_i)\}_{i\in I}$ is an algebraic probability space, denoted $(\cA, \vp)= \trr_{i\in I} (\cA_i, \vp_i)$, where $\cA = \sstat{i\in I} \cA_i$ is the free product of $\{\cA_i\}_{i\in I}$ and $\vp$ is a linear functional on $\cA$ uniquely determined by the following rule: 
for every $\bfi = (i_1,\dots, i_n) \in \fin{I}$ and every $(a_1,a_2,\dots, a_n) \in  \cA_\bfi$, we have
 \begin{equation}\label{eq:monotone}
  \vp(a_1 a_2 \dotsm a_n) = \vp(a_p)\vp(a_1 a_2 \dotsm a_{p-1} a_{p+1}  \dotsm a_n) 
\end{equation}
for any $p \in [n]$ at which the function $[n]\ni k \mapsto i_k \in I$ attains a local maximum, meaning that $p \in \{2,3,\dots, n-1\}, i_{p-1} < i_p > i_{p+1}$, or $p=1, i_1>i_2$, or $p=n, i_{n-1}<i_n$.   We also employ the notation $\varphi = \trr_{i\in I} \vp_i$. 
 
 \item\label{item:AMP} The \emph{antimonotone product} of $\{(\cA_i, \vp_i)\}_{i\in I}$, denoted $\trl_{i\in I} (\cA_i, \vp_i)$, is the algebraic probability space $\trr_{i\in I^{\rm op}} (\cA_i, \vp_i)$. In other words, the defining rule for $(\cA,\vp)=\trl_{i\in I} (\cA_i, \vp_i)$ is the same as \eqref{eq:monotone} except that  ``local maximum'' is replaced with ``local minimum''.  We also employ the notation $\varphi = \trl_{i\in I} \vp_i$. 
  \end{enumerate}
\end{defi}

Below we collect products for (unital or not) 2-algebraic probability spaces. 
 
\begin{defi}[\cite{BLS,BS}]  \label{def:c-free}%\label{eq1} 
Let $I$ be an index set with $\#I \ge2$ and $\{(\cA_i, \vp_i, \psi_i)\}_{i \in I}$ be a family of unital 2-algebraic probability spaces. The \emph{c-free product} of $\{(\cA_i, \vp_i, \psi_i)\}_{i\in I}$, denoted  
$(\cA, \vp, \psi) = \freeproe{i \in I}(\cA_i, \vp_i, \psi_i)$, consists of  
 the unital free product $\cA := \freeproe{i \in I}\cA_i$,  the free product $\psi:= \freeproe{i \in I} \psi_i$ and a unital linear functional $\vp$ determined by the following property:  
\begin{enumerate}[label=\rm(CF),leftmargin=40pt]
\item\label{item:def-c-free} for every $\bfi=(i_1,i_2,\dots, i_n) \in \fin{I}$ and every $(a_1,a_2,\dots, a_n) \in  \cA_{\bfi}$ with $\psi_{i_k}(a_k) = 0$ for all $k\in[n]$, we have 
\begin{equation}\label{123}
\vp (a_1 \cdots a_n) = \prod_{k = 1}^{n} \vp_{i_k}(a_k).
\end{equation} 
\end{enumerate}

We also use the notation $(\vp, \psi) = \freeproe{i \in I}(\vp_i, \psi_i)$, and especially  when $I=\{1,2\}$, we use the notation $\vp= \cf{\vp_1}{\psi_1}{\psi_2}{\vp_2}$. Note that $\vp$ depends on $\vp_1, \vp_2,\psi_1,\psi_2$. 
\end{defi} 

It can be proved recursively that the c-free product $\cf{\vp_1}{\psi_1}{\psi_2}{\vp_2}$ is determined only by $\vp_1, \vp_2, \psi_1, \psi_2$. A few examples will be helpful to see this fact. 

\begin{exa}\label{exa:c-free} In the setting of Definition \ref{def:c-free} for $I=\{1,2\}$ one can prove 
by the standard centering technique and symmetry that for all $a_1, a_1' \in \cA_1$ and $a_2, a_2' \in \cA_2$ 
\begin{align}
&\vp(a_1) = \vp_1(a_1), \quad \vp(a_2) = \vp_2(a_2),    \label{exa:c-free1} \\
&\vp(a_1a_2) = \vp_1(a_1)\vp_2(a_2) = \vp(a_2a_1),   \label{exa:c-free2} \\
&\vp(a_1 a_2 a_1') = \psi_2(a_2)[\vp_1(a_1a_1') - \vp_1(a_1) \vp_1(a_1')] + \vp_1(a_1)\vp_2(a_2) \vp_1(a_1')\quad \text{and}    \label{exa:c-free3}  \\
&\vp(a_2a_1a_2') = \psi_1(a_1)[\vp_2(a_2a_2') - \vp_2(a_2) \vp_2(a_2')] + \vp_2(a_2)\vp_1(a_1) \vp_2(a_2').\label{exa:c-free4} 
\end{align}
For example, formula $\vp(a_1a_2) = \vp_1(a_1)\vp_2(a_2)$ is deduced from the identity 
\[
\vp[(a_1-\psi_1(a_1)1)(a_2-\psi_2(a_2)1)] = [\vp_1(a_1) - \psi_1(a_1)] [\vp_2(a_2) -\psi_2(a_2)]
\]
and \eqref{exa:c-free1}. It is a consequence of the symmetry of c-free product that  formulas \eqref{exa:c-free3} and \eqref{exa:c-free4} show a certain resemblance. 
\end{exa}

\begin{defi}[\cite{Has3}]\label{def:c-monotone} 
Let $I$ be a toset with $\#I \ge2$  and $\{(\cA_i, \vp_i, \psi_i)\}_{i \in I}$ be a family of 2-algebraic probability spaces. The \emph{c-monotone product} of $\{(\cA_i, \vp_i, \psi_i)\}_{i\in I}$, denoted  $(\cA, \vp, \psi) = \trr_{i \in I}(\cA_i, \vp_i, \psi_i)$, consists of  
 the nonunital free product $\cA := \sstat{i \in I}\cA_i$,  the monotone product $\psi:= \trr_{i \in I} \psi_i$ and the unital linear functional $\vp$ determined by the following recursive condition:  
\begin{enumerate}[label=\rm(CM\arabic*),leftmargin=45pt]
\item $\Restr{\vp}{\cA_p} = \vp_{i_p}$ for all $p\in [n]$;  

\item for every $\bfi=(i_1,i_2, \dots, i_n) \in \fin{I}$ with $n\ge2$, every $(a_1,a_2,\dots, a_n) \in  \cA_{\bfi} $ we have
\begin{equation*}
\vp (a_1 \cdots a_n) = \vp(a_1 \cdots a_{p-1}) [\vp(a_p) - \psi(a_p)]\vp(a_{p+1} \cdots a_{n}) + \psi(a_p) \vp(a_1 \cdots a_{p-1}a_{p+1} \cdots a_n)  
\end{equation*} 
for any $p \in [n]$ at which the function $[n]\ni k \mapsto i_k \in I$ attains a local maximum.  Note that we set $ \vp(a_1 \cdots a_{p-1})= 1$ and $\vp(a_{p+1} \cdots a_{n})=1$  if $p=1$ and $p=n$, respectively. 

\end{enumerate}
In case $I=\{1,2\}$,  $\vp$ depends on $\vp_1, \vp_2, \psi_2$ and will be denoted $\cm{\vp_1}{\psi_2}{\vp_2}$. 
 
The \emph{c-antimonotone product} of $\{(\cA_i, \vp_i, \psi_i)\}_{i\in I}$, denoted $\trl_{i \in I}(\cA_i, \vp_i, \psi_i)$, is defined to be the c-monotone product $\trr_{i \in I^{\rm op}}(\cA_i, \vp_i, \psi_i)$ over the opposite toset $I^{\rm op}$ defined in Definition \ref{def:MP-AMP} \ref{item:AMP}. In case $I=\{1,2\}$, we employ  the notation $(\cA_1, \vp_1,\psi_1) \trl (\cA_2,\vp_2,\psi_2) = (\cA, \cam{\vp_1}{\psi_1}{\vp_2}, \psi_1 \trl\psi_2)$. 
\end{defi}

\begin{rem} 
The positivity structure is preserved by all the products above (free, monotone, antimonotone, Boolean, c-free, c-monotone and c-antimonotone products), i.e.\ if the given spaces are (unital or not) $m$-$\ast$-algebraic probability space then so is the resulting space. 
\end{rem}

The above products for the special index set $I=\{1,2\}$ provides a binary operation on (unital) $m$-algebraic probability spaces. 
In general, a binary operation that associates with any pair  $((\cA_1, \vp_1^1, \dots, \vp^m_1 )$, $(\cA_2, \vp_2^1, \dots, \vp^m_2))$ of  $m$-algebraic probability spaces a new $m$-algebraic probability space $(\cA, \varphi^1, \dots, \varphi^m)$, where $\cA = \cA_1 \sstar \cA_2$, is called a \emph{product for  $m$-algebraic probability spaces}. Similarly, the notion of product for unital $m$-probability spaces is defined, where the product algebra $\cA$ is then taken to be $\cA_1\freeprod \cA_2$, the unital free product.  The notion of product is one of the building blocks in \cite{MS17}. 

A product $\square$ for $m$-algebraic probability spaces  is called \emph{associative} if for any $m$-algebraic probability spaces $(\cA_i, \vp_i^1, \dots, \vp_i^m), i=1,2,3$, the identity 
\begin{align*}
& [(\cA_1, \vp_1^1, \dots, \vp^m_1 ) \square (\cA_2, \vp_2^1, \dots, \vp_2^m) ] \square (\cA_3, \vp_3^1, \dots, \vp_3^m) \\
&\qquad\qquad= (\cA_1, \vp_1^1, \dots, \vp^m_1 ) \square [(\cA_2, \vp_2^1, \dots, \vp_2^m)  \square (\cA_3, \vp_3^1, \dots, \vp_3^m) ] 
\end{align*}
holds under the natural isomorphism $(\cA_1 \sstar \cA_2) \sstar \cA_3 \simeq \cA_1 \sstar (\cA_2 \sstar \cA_3)$. The associativity for products for unital $m$-algebraic probability spaces is defined similarly. 
Associativity is satisfied by free, Boolean, monotone, antimonotone, c-free, c-monotone and c-antimonotone products. The associativity of monotone and c-monotone products are  proved in \cite[Theorem 2.4]{Fra4} and \cite[Theorem 3.7]{Has3}, respectively. 

A product $\square$ for $m$-algebraic probability spaces  is called \emph{symmetric} if for any $m$-algebraic probability spaces $(\cA_i, \vp_i^1, \dots, \vp_i^m), i=1,2$, the identity 
\[
(\cA_1, \vp_1^1, \dots, \vp^m_1 ) \square (\cA_2, \vp_2^1, \dots, \vp_2^m)  =  (\cA_2, \vp_2^1, \dots, \vp^m_2 ) \square (\cA_1, \vp_1^1, \dots, \vp_1^m) 
\]
holds under the natural identification $\cA_1 \sstar \cA_2 \simeq \cA_2 \sstar \cA_1$. A similar definition can be given for products for unital $m$-algebraic probability spaces. 
Free, Boolean and c-free products are symmetric, although monotone, antimonotone, c-monotone and c-antimonotone products are not.

%%%%%%%%%%%%%%%%%%%%%%%%%%%%%%%%%%%%%%%%%%%%%%%%%%%%%%%%%%%%%%%
\subsection{Aims and motivations} \label{sec:aim_motivation}

%In the present paper, we propose a notion of independence referring to a \emph{triplet of linear functionals} or a product for unital 3-algebraic probability spaces. 
%The new product, called the $\AP$-free product, unifies free, Boolean, monotone and antimonotone products in a certain way. In some sense this product is not completely new because it is a combination of c-free products. However, when looking at the \emph{associativity of products}, the $\AP$-free product is not a special case of c-free product. Especially, difference between the c-free product and the $\AP$-free product will be clear when one formulates cumulants. 

Aims and motivations come from the notion of associativity of products for (unital or nonunital) $m$-algebraic probability spaces. It is known that the c-free product of unital 2-algebraic probability spaces generalizes the Boolean product and free product. Also, the monotone product of two algebraic probability spaces is a special case of c-free product. However, unlike the Boolean and free products,  the associativity of monotone product is invisible in the associativity of the c-free product. The present paper appeared from attempts to construct a new product that unifies the mentioned products including associativity. 

To see further details, let us start from a connection of  the c-free product and other products.  
 By definition,  the free product of unital algebraic probability spaces is a special case of the c-free product: 
\[
\cf{\vp_1}{\vp_1}{\vp_2}{\vp_2} = \vp_1 \ast \vp_2.
\]
 Boolean and monotone products appear from c-free products with suitable unitization of algebraic probability spaces. Let $(\cA_i, \vp_i)$ be  algebraic probability spaces for $i=1,2$. Recall that $(\wh{\cA}_i,\wh{\vp}_i)$ denotes the unitization of $(\cA_i, \vp_i)$, see Notation \ref{notation}. The \emph{delta functional} $\delta_i$ on $\wh{\cA}_i$ is then defined by $\delta_i := \wh{0}$, i.e.\  $\delta_i (\lambda + a) = \lambda$ for $\lambda \in \comp$ and $a \in \cA_i$.  From now on, when  referring to delta functionals, we always assume that the algebras are constructed in the above way. Under the  natural isomorphism  $\wh{\cA}_1 \freeprod \wh{\cA}_2\simeq \wh{\cA_1 \sstar \cA_2}$ it holds that 
\begin{equation} \label{eq:c-free-boole}
\cf{\wh{\vp}_1}{\delta_1}{\delta_2}{\wh{\vp}_2} = \wh{\vp_1 \diamond \vp_2},    \quad \text{or equivalently} \quad (\wh\vp_1,\delta_1) \freeprod (\wh\vp_2,\delta_2) = (\wh{\vp_1 \diamond \vp_2}, \delta_1 \ast \delta_2). 
\end{equation}
Note that $\delta_1 \ast \delta_2$ is the delta functional on $\wh{\cA}_1 \freeprod \wh{\cA}_2$ for the decomposition $\wh{\cA}_1 \freeprod \wh{\cA}_2 \simeq \comp \oplus (\cA_1 \sstar \cA_2)$.  

 Moreover, 
Franz proved in \cite{Fra} that the monotone product (resp.\ antimonotone product) appears as 
\begin{equation}\label{eq:c-free-monotone}
\cf{\wh{\vp}_1}{\delta_1}{\wh{\vp}_2}{\wh{\vp}_2} = \wh{\vp_1 \trr \vp_2}\qquad \text{(resp. $\cf{\wh{\vp}_1}{\wh\vp_1}{\delta_2}{\wh\vp_2} = \wh{\vp_1 \trl \vp_2}$).}
\end{equation}

The associativity of the c-free product is transferred to the Boolean product, i.e.\ the associativity of c-free product implies
\begin{equation} \label{eq:associative1}
[(\wh{\vp}_1, \delta_1) \freeprod (\wh{\vp}_2, \delta_2) ]\freeprod (\wh\vp_3,\delta_3) = (\wh{\vp}_1, \delta_1) \freeprod [(\wh{\vp}_2, \delta_2) \freeprod (\wh\vp_3,\delta_3)],   
\end{equation}
which, combined with \eqref{eq:c-free-boole}, further implies 
\begin{equation}\label{eq:associative2}
\wh{(\vp_1 \diamond \vp_2)\diamond \vp_3} = \wh{\vp_1 \diamond (\vp_2\diamond \vp_3)} \quad \text{or equivalently} \quad (\vp_1 \diamond \vp_2)\diamond \vp_3 = \vp_1 \diamond (\vp_2\diamond \vp_3). 
\end{equation}
Also, the associativity of the c-free product is trivially transferred to the free product because the second linear functional in the c-free product is just the free product.  
%\[
%[(\vp_1, \vp_1) \freeprod (\vp_2, \vp_2) ]\freeprod (\vp_3,\vp_3) = (\vp_1, \vp_1) \freeprod [(\vp_2, \vp_2) \freeprod (\vp_3,\vp_3)]. 
%\]

By contrast, the associativity of c-free product does not seem to be transferred to the monotone product because of the asymmetry of $\delta_1$ and $\wh{\vp}_2$ in \eqref{eq:c-free-monotone}. 
 Consequently, the central limit theorem or cumulants for c-freeness do not seem to be transferred to monotone independence, either.   A similar observation can be made for the 
antimonotone product.  

This observation on the associativity lead to the introduction of the c-monotone product in \cite{Has3} for $2$-algebraic probability spaces.  It can be defined in terms of c-free product in the way
\[
(\cA_1, \vp_1, \psi_1) \trr (\cA_2, \vp_2, \psi_2) = (\cA_1\sstar  \cA_2, \Restr{(\cf{\wh\vp_1}{\delta_1}{\wh\psi_2}{\wh\vp_2})}{\cA_1\sstar \cA_2}, \psi_1 \trr \psi_2) 
\]
in the category of 2-algebraic probability spaces. 

%The left component $\vp_1 {}_{\delta_1}\!\! \freeproe{\psi_2}\!\vp_2$ is also denoted by $\vp_1 \trr_{\psi_2} \vp_2$. 
The c-monotone product is associative as demonstrated in \cite{Has3}, although this does not come as a consequence of the associativity of the c-free product. 
The c-monotone product generalizes the monotone and Boolean products in the ways% (recall that the zero linear functional is a restricted state)
\begin{align*}
&(\vp_1, \vp_1) \trr (\vp_2, \vp_2) = (\vp_1 \trr \vp_2, \vp_1 \trr \vp_2), \\
&(\vp_1, 0) \trr (\vp_2, 0) = (\vp_1 \diamond \vp_2, 0). 
\end{align*}  
Therefore, the associativity of c-monotone product is transferred to the monotone and Boolean products along the lines \eqref{eq:associative1} -- \eqref{eq:associative2}. 
In this sense, the c-monotone product well unifies the monotone and Boolean products, but it does not unify the free product. The main objective of the present paper is to unify c-monotone and c-free products including their associativity, thus well unifying the free, monotone and Boolean products. 
% and the paragraph following  \eqref{eq:associative2}. 

%Similarly, the c-antimonotone product can be formulated in terms of c-free product in the way
%\[
%(\cA_1, \vp_1, \psi_1) \trl (\cA_2, \vp_2, \psi_2) := (\cA_1\sstar \cA_2, \Restr{\cf{\wh\vp_1}{\wh\psi_1}{\delta_2}{\wh\vp_2}}{\cA_1 \sstar \cA_2}, \psi_1 \trl \psi_2).  
%\]
%The c-antimonotone product generalizes the Boolean and antimonotone products: 
%\begin{align*}
%&(\vp_1, \vp_1) \trl (\vp_2, \vp_2) = (\vp_1 \trl \vp_2, \vp_1 \trl \vp_2), \\
%&(\vp_1, 0) \trl (\vp_2, 0) = (\vp_1 \diamond \vp_2, 0). 
%\end{align*}   

%%%%%%%%%%%%%%%%%%%%%%%%%%%%%%%%%%%%%%%%%%%%%%%%%%%%%%%%%%%%%%%%%%%
\subsection{Main results} \label{subsec:main}
In this paper we construct an associative product for unital 3-algebraic probability spaces so that its associativity is transferred to any of free, monotone, antimonotone, Boolean, c-free, c-monotone and c-antimonotone products. Such a product provides a more complete unification of the existing products; for example,  the notions of cumulants can be unified. The  new product can be described as 
\begin{equation}\label{eq:alpha_product}
(\vp_1, \psi_1, \theta_1) \fmp (\vp_2, \psi_2, \theta_2) = (\cf{\vp_1}{\theta_1}{\psi_2}{\vp_2}, \cf{\psi_1}{\theta_1}{\psi_2}{\psi_2}, \cf{\theta_1}{\theta_1}{\psi_2}{\theta_2}), 
\end{equation}
which will be called the \textit{$\AP$-free product}.\footnote{This product was formerly named indented product.} Although this product is made out of c-free products, the associativity of $\AP$-free product is not a consequence of the associativity of the c-free product.  In particular, the product for unital 2-algebraic probability spaces 
\begin{equation}\label{eq:a-freeproduct}
(\psi_1, \theta_1) \fmp (\psi_2, \theta_2) = (\cf{\psi_1}{\theta_1}{\psi_2}{\psi_2}, \cf{\theta_1}{\theta_1}{\psi_2}{\theta_2}), 
\end{equation} 
that appears in the second and third components of $\AP$-free product, is also associative. This will be referred to as the \textit{$\BGP$-free product}\footnote{This product was formerly named ordered free product with abbreviation ``o-free product''.} and denoted by the same symbol $\fmp$. 

%Furthermore, $\BGP$-free and $\AP$-free products are naturally expected to have connections with the concept of \textit{matricial freeness} introduced by R. Lenczewski recently \cite{Len0}. We however leave this direction to a future research. 
%Note that the ``opposite $\AP$-free product'' 
% \[
%(\vp_1, \psi_1, \theta_1) \rightthreetimes (\vp_2, \psi_2, \theta_2) := (\cf{\vp_1}{\psi_1}{\theta_2}{\vp_2}, \cf{\psi_1}{\psi_1}{\theta_2}{\psi_2}, \cf{\theta_1}{\psi_1}{\theta_2}{\theta_2})
%\]   
%and the ``opposite $\BGP$-free product'' 
%\[ (\psi_1, \theta_1) \rightthreetimes (\psi_2, \theta_2) := (\cf{\psi_1}{\psi_1}{\theta_2}{\psi_2}, \cf{\theta_1}{\psi_1}{\theta_2}{\theta_2}) \]
%are also associative. The structures of these products are essentially equal to the $\AP$-free product and $\BGP$-free product, respectively, and hence will not be mentioned later. 

By definition, the $\AP$-free product unifies other associative products as follows:  
\begin{align*}
(\vp_1, \psi_1, \psi_1) \fmp (\vp_2, \psi_2, \psi_2) &= (\cf{\vp_1}{\psi_1}{\psi_2}{\vp_2}, \psi_1 \ast \psi_2, \psi_1 \ast \psi_2),   
&& (\text{c-free})   \\
(\vp_1, \vp_1, \vp_1) \fmp (\vp_2, \vp_2, \vp_2)  &= (\vp_1 \ast \vp_2, \vp_1 \ast \vp_2, \vp_1 \ast \vp_2),    && (\text{free})    \\
(\wh\vp_1,\wh \psi_1, \delta_1) \fmp (\wh\vp_2, \psi_2, \wh\delta_2) &= (\wh{\cm{\vp_1}{\psi_2}{\vp_2} }, \wh{\psi_1 \trr \psi_2}, \delta_1 \ast \delta_2),   &&(\text{c-monotone}) \\
(\wh\vp_1, \wh\vp_1, \delta_1) \fmp (\wh\vp_2, \wh\vp_2, \delta_2) &= (\wh{\vp_1 \trr \vp_2}, \wh{\vp_1\trr \vp_2}, \delta_1 \ast \delta_2),   && (\text{monotone}) \\
(\wh\vp_1, \delta_1, \wh\psi_1) \fmp (\wh\vp_2, \delta_2, \wh\psi_2) &= (\wh{  \cam{\vp_1}{\psi_1}{\vp_2} }, \delta_1 \ast \delta_2, \wh{\psi_1 \trl \psi_2} ),     && (\text{c-antimonotone})   \\
(\wh\vp_1, \delta_1, \wh\vp_1) \fmp (\wh\vp_2, \delta_2, \wh\vp_2) &= (\wh{\vp_1  \trl \vp_2}, \delta_1 \ast \delta_2, \wh{\vp_1  \trl \vp_2}),    
&& (\text{antimonotone})  \\ 
(\wh\vp_1, \delta_1, \delta_1) \fmp (\wh\vp_2, \delta_2, \delta_2) &= (\wh{\vp_1 \diamond \vp_2}, \delta_1 \ast \delta_2, \delta_1 \ast \delta_2). && (\text{Boolean})  
\end{align*}
Consequently, the associativity of $\AP$-free product is transferred to the seven products above.  
The way the associativity is transferred from a product to another product is visualized in Fig. \ref{dia21}.  

\begin{figure}[t]
\begin{center}
\begin{tikzpicture}
\draw[->] (0.1,3.7) -- (0.1,2.3); 
\draw[->] (-0.3,3.7) -- (-3+0.3,2+0.3); 
\draw[->] (0.3,3.7) -- (3,2+0.3); 
\draw[->,dashed] (-0.1,2.3) -- (-0.1,3.7); 
\draw[->,dashed] (-0.7,1.95) -- (-3+0.7,1.95); 
\draw[->,dashed] (0.7,1.95) -- (3-0.8,1.95); 
\draw[->,dashed] (0,1.7) -- (0, 0.3); 
\draw[->] (0.3,1.7) -- (3, 0.3); 
\draw[->] (-0.3,1.7) -- (-3+0.3, 0.3); 
\draw[->] (-3+0.3,1.7) -- (-0.3, 0.3); 
\draw[->] (-3,1.7) -- (-3, 0.3); 
\draw[->] (3-0.3,1.7) -- (0.3, 0.3); 
\draw[->] (3.3,1.7) -- (3.3, 0.3); 
\node at (0,4) {$\AP$-free};
\node at (0,2) {c-free};
\node at (0,0) {monotone};
\node at (3.3,0) {Boolean};
\node at (-3,0) {free};
\node at (3.4,2)  {c-monotone}; 
\node at (-3.1,2)  {$\BGP$-free}; 
\node at (8,4) {three linear functionals};
\node at (8,2) {two linear functionals};
\node at (8,0) {single linear functional};
\end{tikzpicture}
\end{center}
\caption{Each arrow means that the initial object is reduced to the terminal one if appropriate linear functionals are selected. Undashed arrows preserve the associativity, although dashed arrows do not. Note that one could incorporate antimonotone and c-antimonotone into the above figure but they are omitted for simplicity.} \label{dia21}
\end{figure}

%The definition (\ref{eq:alpha_product}) will be understood as a natural extension of (\ref{eq:a-freeproduct}).  

The rest of the paper is organized in the following way.  In Section \ref{prod} the associativity of the $\AP$-free product is established, which is the fundamental part of this paper. Also, the definition of $\AP$-free product is extended to an arbitrary family of unital 3-algebraic probability spaces indexed by a totally ordered set. 

In Section \ref{rep}, assuming the positivity structure, we construct $\ast$-representations of the unital free product of unital $\ast$-algebras which provide $\AP$-free operators on the free product of pre-Hilbert spaces. %Motivations for this section come from similar constructions for freeness by Avitzour \cite{Avi} and Voiculescu \cite{V1}, for c-freeness by Bo\.zejko and Speicher \cite{BS} and for c-monotone independence by Popa \cite{Pop}. 

%apers \cite{Avi,BS,Mur3,Pop,V1}. 

The remaining contents are mainly devoted to cumulants. In free probability theory, there have been many researches on combinatorial aspects of cumulants since Speicher introduced noncrossing set partitions in \cite{Spe2}. 
In Section \ref{cum}, we define ``$\AP$-free cumulants'' along the lines of the general theory of cumulants for spreadability systems \cite{HL17}. The moment-cumulant formula is established, where ``ordered noncrossing set partitions'' play a critical role.  Because the $\AP$-free product is asymmetric, i.e.\ the $\AP$-freeness of $(a,b)$ does not imply the $\AP$-freeness of $(b,a)$, cumulants do not have the vanishing property as discussed in \cite{HL17}. 
Since the associativity of $\AP$-free product is transferred to the seven products above,  the moment-cumulant formula for $\AP$-free cumulants induces the moment-cumulant formulas for the other cumulants. 
%In particular, we obtain moment-cumulant formula for c-monotone independence. 

Section \ref{sec:convolution} characterizes additive and multiplicative $\AP$-free convolutions in terms of reciprocal of Cauchy transforms. This section contains some heuristics on how the $\AP$-free product (\ref{eq:alpha_product}) and $\BGP$-free product (\ref{eq:a-freeproduct}) were discovered.

Section \ref{sec:generating_function} treats generating functions of $\AP$-free cumulants (for single variables) and their relationships to the generating functions of moments. 
It turns out that a certain set of differential equations relate these generating functions. This generalizes known differential equations for the Cauchy transforms in free probability and monotone probability. 
As an application of $\AP$-free cumulants, we prove the central limit theorem. Because of the transfer of associativity, the central limit theorem for $\AP$-freeness generalizes those for the seven independences mentioned above. The limit distribution consists of three Kesten distributions, which are known to appear already in the c-free and c-monotone central limit theorems.

%%%%%%%%%%%%%%%%%%%%%%%%%%%%%%%%%%%%%%%%%%%%%%%%%%%%%%%%%%%%%%%%%%%%%%%
\subsection*{Note} This work was first posted to the preprint server arXiv in 2010. 
%Since then, some people mentioned this work in their papers, e.g.\ in \cite{CG}. %Rather recently this old preprint seems to begin getting attention of people working in non-commutative probability, see e.g.\ \cite{CG}. 
Because the old version was not easy to read and contained insufficient arguments, the author has thought of making a revision. 
As of 2022, the author has finally finished revising the manuscript thoroughly. Especially,  the treatment of cumulants is revised greatly on the basis of a general theory of cumulants \cite{HL17} which was not available as of the first version of this manuscript.  The author believes that this approach makes the structure of proof clearer. 
The title and some terminologies have also been changed; the previous title was ``New associative product for triplets of states generalizing free, monotone, antimonotone, Boolean, conditionally free and conditionally monotone products''.

%%%%%%%%%%%%%%%%%%%%%%%%%%%%%%%%%%%%%
%%%%%%%%%%%%%%%%%%%%%%%%%%%%%%%%%%%%%%
\section{The $\AP$-freeness}\label{prod}

%%%%%%%%%%%%%%%%%%%%%%%%%%%%%%%%%%%%%%%%%%%%%%%%%%%%%%%%%%%%%%%%%%%%%%%%%%%%%%%%
\subsection{A characterization of the c-free product}
To prove that the $\AP$-free product is associative, we need a better understanding of mixed moments with respect to the c-free product of pairs of unital linear functionals. 
This subsection offers an alternative characterization of the c-free product for the index set $I=\{1,2\}$.

\begin{lem}\label{lem:c-free}
Let $(\cA_i, \vp_i,\psi_i), i=1,2$ be unital 2-algebraic probability spaces. A unital linear functional $\vp$ on $\cA_1 \freeprod \cA_2$ coincides with $\cf{\vp_1}{\psi_1}{\psi_2}{\vp_2}$ if and only if the following condition holds: 

\begin{enumerate}[label=\rm(CF'),leftmargin=40pt]

\item\label{item:c-free}  $\vp (a_1 \cdots a_n) = 0$ holds whenever  $n \in \N$,  $i_1,i_2,\dots, i_n \in \{1,2\}$ with $i_1 \neq \cdots \neq i_n$ and $(a_1,a_2,\dots,  a_n) \in \cA_{i_1} \times \cA_{i_2} \times \cdots \times \cA_{i_n}$ with $\vp_{i_1}(a_1) = 0$ and $\psi_{i_k}(a_k) = 0$ for all $k \in\{2,3,\dots, n-1\}$. 
\end{enumerate}
By symmetry, condition \ref{item:c-free} can be replaced by the following: 
\begin{enumerate}[label=\rm(CF''),leftmargin=40pt]
\item\label{item:c-free2}  $\vp (a_1 \cdots a_n) = 0$ holds whenever  $n \in \N$,  $i_1,i_2,\dots, i_n \in \{1,2\}$ with $i_1 \neq \cdots \neq i_n$ and $(a_1,a_2,\dots,  a_n) \in \cA_{i_1} \times \cA_{i_2} \times \cdots \times \cA_{i_n}$ with $\vp_{i_n}(a_n) = 0$ and $\psi_{i_k}(a_k) = 0$ for all $k \in\{2,3,\dots, n-1\}$. \end{enumerate}
\end{lem}
\begin{rem} In case $n=1,2$, the irrelevant assumption ``$\psi_{i_k}(a_k) = 0$ for all $k \in\{2,3,\dots, n-1\}$'' is to be deleted. For example,  condition \ref{item:c-free} is to be understood as follows. 
\begin{itemize}
\item $n=1$: $\vp (a) = 0$ whenever  $a \in \cA_{i}$ and $\vp_{i}(a) = 0$. 
 \item $n=2$: $\vp (a_1 a_2) = 0$ whenever $i_1,i_2 \in \{1,2\}$ with $i_1 \neq i_2$ and $(a_1,a_2) \in \cA_{i_1} \times \cA_{i_2}$ with $\vp_{i_1}(a_1) = 0$. 
 \end{itemize}
 \end{rem}

\begin{proof}[Proof of Lemma \ref{lem:c-free}] Once the condition $\vp=\cf{\vp_1}{\psi_1}{\psi_2}{\vp_2}$ is shown to be equivalent to \ref{item:c-free}, it is also equivalent by symmetry to \ref{item:c-free2}. 

Suppose that condition  \ref{item:c-free} holds. Let us take $(a_1,a_2,\dots, a_n)$ as described in condition \ref{item:def-c-free} in the definition of c-free product (Definition \ref{def:c-free}). 
The goal is to prove 
\begin{equation}\label{eq:goal}
\varphi(a_1 a_2  \cdots a_n) = \prod_{k \in [n]} \vp_{i_k}(a_k). 
\end{equation}
For $n=1$ condition \ref{item:c-free} yields $\vp(a- \vp_{i_1}(a_1)1)=0$, which is exactly \eqref{eq:goal}.  For $n\ge2$, we introduce the centering $a_1 = \mathring{a}_1 + \rho_1 1$, where $\rho_1 = \vp_{i_1}(a_1)$ and then proceed as 
\[
\vp(a_1\cdots a_n) = \vp(\mathring{a}_1 a_2 \cdots a_n) + \rho_1  \vp(a_2 \cdots a_n), 
\]
which equals $\vp_{i_1}(a_1) \vp(a_2 \cdots a_n) $ because of \ref{item:c-free}.  The same arguments can be repeated to arrive at the formula $\vp(a_1\cdots a_n) = \vp_{i_1}(a_1)\vp_{i_2}(a_2) \cdots \vp_{i_n}(a_n)$ as desired.

Conversely,  assume that $\vp = \cf{\vp_1}{\psi_1}{\psi_2}{\vp_2}$.  Condition \ref{item:c-free} for $n=1$ is a consequence of  \eqref{exa:c-free1} in Example \ref{exa:c-free}. Taking $n\ge2$ and $(a_1,a_2,\dots, a_n)$ as described in condition \ref{item:c-free}, we introduce the centering $a_k = \mathring{a}_k + \lambda_k 1$ where $\lambda_k=\psi_{i_k}(a_k)$ for $k \in \{1,n\}$ and apply the definition of c-free product to proceed as 
\begin{equation*}
\begin{split}
\vp (a_1 \cdots a_n) &= \vp(\mathring{a}_1a_2 \cdots a_{n-1}\mathring{a}_n)  + \vp(\mathring{a}_1a_2 \cdots a_{n-1})  \lambda_n
                             +\lambda_1 \vp(a_2 \cdots a_{n-1}\mathring{a}_n) + \lambda_1 \vp(a_2 \cdots a_{n-1}) \lambda_n \\
                         &= \vp_{i_1}(\mathring{a}_1) \left[ \prod_{k=2}^{n-1} \vp_{i_k}(a_k) \right]  \vp_{i_n}(\mathring{a}_n) +  \vp_{i_1}(\mathring{a}_1) \left[ \prod_{k=2}^{n-1} \vp_{i_k}(a_k)  \right] \psi_{i_n}(a_n) \\
                         &\quad +\psi_{i_1}(a_1)  \left[ \prod_{k=2}^{n-1} \vp_{i_k}(a_k)  \right] \vp_{i_n}(\mathring{a}_n)  +   \psi_{i_1} (a_1)  \left[  \prod_{k=2}^{n-1} \vp_{i_k}(a_k)  \right] \psi_{i_n}(a_n).  
\end{split}
\end{equation*}
Because the assumption $\vp_{i_1}(a_1)=0$ implies $\vp_{i_1}(\mathring{a}_1) = - \psi_{i_1} (a_1)$, we deduce from straightforward calculations that 
\[
\vp_{i_1}(\mathring{a}_1)  \vp_{i_n}(\mathring{a}_n) +    \vp_{i_1}(\mathring{a}_1) \psi_{i_n}(a_n) +\psi_{i_1}(a_1)\vp_{i_n}(\mathring{a}_n)  +   \psi_{i_1} (a_1) \psi_{i_n}(a_n) =0
\]
and hence  $\vp (a_1 \cdots a_n)=0$ as desired. 
\end{proof}

The following fact will not be directly used later but is stated here for potential use in future work. The proof is similar to that of Lemma \ref{lem:c-free} and is omitted. 
It roughly says that if we put the assumption $\vp_2=\psi_2$ then the kernel conditions are unnecessary for $a_1$ and $a_n$ when $n\ge4$.

\begin{lem} \label{lem:afree} Let $(\cA_1, \vp_1,\psi_1)$ be a unital 2-algebraic probability space and $(\cA_2, \vp_2)$ be a unital algebraic probability space. A unital linear functional $\vp$  on $\cA_1 \freeprod \cA_2$ coincides with $\cf{\vp_1}{\psi_1}{\vp_2}{\vp_2}$  if and only if the following conditions hold:
\begin{enumerate}[label=\rm(\roman*)]

\item\label{item:afree1}  for all $a_1, a_1' \in \cA_1$ and $a_2, a_2' \in \cA_2$ we have 
\begin{align*}
%&\vp(a_1) = \vp_1(a_1), \quad \vp(a_2) = \vp_2(a_2), \\
%&\vp(a_1a_2) = \vp_1(a_1)\vp_2(a_2) = \vp(a_2a_1), \\
&\vp(a_1 a_2 a_1') = \vp_2(a_2)\vp_1(a_1a_1') \quad \text{and} \\
&\vp(a_2a_1a_2') = \psi_1(a_1)[\vp_2(a_2a_2') - \vp_2(a_2) \vp_2(a_2')] + \vp_2(a_2)\vp_1(a_1) \vp_2(a_2'); 
\end{align*}

\item\label{item:afree2}  $\vp (a_1 \cdots a_n) = 0$  
whenever $n \geq 4$, $i_1,i_2,\dots, i_n \in \{1,2\},i_1 \neq \cdots \neq i_n$ and $(a_1,a_2,\dots,  a_n) \in \cA_{i_1} \times \cA_{i_2} \times \cdots \times \cA_{i_n}$ with $\psi_1(a_k) = 0$ for all $k\in\{2,3,\dots, n-1\}$ such that $i_k = 1$ and $\vp_2(a_k) = 0$ for all $k\in\{2,3,\dots, n-1\}$ such that $i_k = 2$. 
\end{enumerate}

\end{lem}
\begin{rem}
For words of length 3, it holds true that 
$\vp(a_1 a_2 a_1')=0$ whenever $\vp_2(a_2)=0$, although $\vp(a_2 a_1 a_2')$ may not vanish even when $\psi_1(a_1)=0$. In order to get $\vp(a_2 a_1 a_2')=0$, we need to  assume e.g.\ $\psi_1(a_1)=0$ \emph{and} $\vp_2(a_2)=0$ as in Lemma \ref{lem:c-free}. 
\end{rem}

\subsection{The $\AP$-freeness and $\BGP$-freeness} 

We will extend the definition of $\AP$-free product \eqref{eq:alpha_product} to an arbitrary ordered family of unital 3-algebraic probability spaces and then establish the associativity.  Notions of peaks and bottoms are crucial.

\begin{defi} Let $I$ be a toset. For each $\bfi= (i_1, \cdots, i_n) \in \fin{I}$ with $n \in \N$, a \textit{peak} (resp.\ a \textit{bottom}) with respect to $\bfi$ is a number $k \in \{2,3,\dots, n-1\}$ at which the function $[n]\ni p\mapsto i_p \in I$ attains a local maximum (resp.\ local minimum), i.e.\ $i_{k-1}< i_k >i_{k+1}$  (resp.\ $i_{k-1}> i_k < i_{k+1}$).  Let $\Peak(\bfi)$ and $\Bottom(\bfi)$ be the sets of the peaks and the bottoms, respectively.  Note that $ \Peak(\bfi)=\Bottom(\bfi) =\emptyset$ for $n\in\{1,2\}$.  
%For instance, $\Peak(i_1, \cdots, i_{9}) = \{3,7\}$ and $\Bottom(i_1, \cdots, i_{9}) = \{6,8 \}$ in Fig.\ \ref{fig:peaks}.  
\end{defi}
 
\begin{comment}
\begin{figure}[t]
\begin{center}
\begin{tikzpicture}[scale=0.7]
\draw[-] (1,1) -- (2,2) -- (3,4) -- (4,3) -- (5,2) -- (6,1) -- (7,4) -- (8,3) -- (9,5);
\node at (1,0) {1};
\node at (2,0) {2};
\node at (3,0) {3};
\node at (4,0) {4};
\node at (5,0) {5};
\node at (6,0)  {6}; 
\node at (7,0)  {7}; 
\node at (8,0) {8};
\node at (9,0) {9};
\node at (1,1) {$\bullet$};
\node at (2,2) {$\bullet$};
\node at (3,4) {$\bullet$};
\node at (4,3) {$\bullet$};
\node at (5,2) {$\bullet$};
\node at (6,1) {$\bullet$};
\node at (7,4) {$\bullet$};
\node at (8,3) {$\bullet$};
\node at (9,5) {$\bullet$};
\end{tikzpicture}
\end{center}
\caption{$i_1 = i_6 = 1$, $i_2 = i_5 = 2$, $i_4 = i_8 = 3$, $i_3 = i_7 =4$,  $i_9 = 5$.} \label{fig:peaks}
\end{figure}
\end{comment}

On peaks and bottoms and an endpoint 1 or $n$ we consider kernel conditions with respect to suitable linear functionals. 
For each $n\ge1, \bfi =(i_1,i_2,\dots, i_n)\in\fin{I}$ and $k\in[n]$, let $\alpha_{\bfi,k}$ be the linear functionals on $\cA_{i_k}$ defined by 
\begin{equation}\label{eq:alpha}
\alpha_{\bfi,k} = 
\begin{cases} \varphi_{i_1} & \text{if~} k=1, \\
\psi_{i_k} & \text{if~} k \in \Peak(\bfi), \\
\theta_{i_k} & \text{if~} k \in \Bottom(\bfi), \\
0& \text{otherwise} 
\end{cases}
\quad (n\ge2), \quad \text{and} \quad  \alpha_{(i_1),1}=\vp_{i_1}\quad (n=1). 
  \end{equation}
With a slight abuse of notation, we also write $\alpha_{\bfi,k}^{(\vp,\psi,\theta)}$ when stressing the dependence on linear functionals. As special cases we also introduce 
\[
\beta_{\bfi,k} = \beta_{\bfi,k}^{(\psi,\theta)} := \alpha_{\bfi,k}^{(\psi,\psi,\theta)}\quad \text{and} \quad  \gamma_{\bfi,k} = \gamma_{\bfi,k}^{(\psi,\theta)}:= \alpha_{\bfi,k}^{(\theta,\psi,\theta)}. 
\]
Note that if $k \in \{2,3,\dots, n\} \setminus [\Peak(\bfi) \cup \Bottom(\bfi)] $ then $\ker \alpha_{\bfi,k}= \ker \beta_{\bfi,k} = \ker \gamma_{\bfi,k} = \cA_{i_k}$. 

Finally, the following linear functionals $\omega_{\bfi,k}=\omega_{\bfi,k}^{(\vp,\psi,\theta)}$ are useful when discussing operators on the free product pre-Hilbert space in Section \ref{rep}: 
\begin{equation} \label{eq:omega}
\omega_{\bfi,k} = 
\begin{cases} \varphi_{i_n} & \text{if~} k=n, \\
\psi_{i_k} & \text{if~} k \in \Des(\bfi), \\
 \theta_{i_k} & \text{if~} k \in \Asc(\bfi), 
\end{cases}
\quad (n\ge2), \quad \text{and} \quad  \omega_{(i_1),1}=\vp_{i_1}\quad (n=1), 
\end{equation}
where $\Des(\bfi)$ and $\Asc(\bfi)$ are defined by 
\[
\Des(\bfi)=\{k \in [n-1]: i_k > i_{k+1}\} \quad \text{and} \quad \Asc(\bfi)=\{k \in [n-1]: i_k < i_{k+1}\}
\]
and are called  the sets of descents and ascents of $[n-1]$ regarding the sequence $\bfi$, respectively. 
%Note also that $\ker \omega_{\bfi,k} \subseteq \ker\alpha_{\bfi,k}$ for all $\bfi \in \fin{I}$ and $k \in[|\bfi|]$. 

\begin{defi}\label{def:alpha_indep} Let $I$ be a toset with $\# I \ge2$.

\begin{enumerate}[label=\rm(\arabic*)]

\item\label{item:AP_prod} ($\AP$-free product) Let $\{(\cA_i, \vp_i, \psi_i, \theta_i)\}_{i\in I}$ be a family of unital 3-algebraic probability spaces.  
Then its \emph{$\AP$-free product} $\fmp_{i\in I} (\cA_i, \vp_i, \psi_i, \theta_i)$ is a tuple $(\cA, \vp, \psi, \theta)$,  where $\cA = \freeproe{i\in I} \cA_i$ is the unital free product and $(\vp,\psi, \theta)$ is a triplet of unital linear functionals such that the following conditions hold: 

\begin{enumerate}[label=$(\AP)$]
\item\label{item:AP}  $\vp(a_1 \cdots a_n) = 0$ whenever $\bfi =(i_1,i_2,\dots,i_n)\in \fin{I}$ and $(a_1,a_2,\dots, a_n) \in \ker \alpha_{\bfi,1} \times \ker \alpha_{\bfi,2} \times \cdots \times \ker \alpha_{\bfi,n}$; 
\end{enumerate}
\begin{enumerate}[label=$(\beta)$]
%\item\label{item:AF1}   $\Restr{\psi}{\cA_i} = \psi_i$ and $\Restr{\theta}{\cA_i} = \theta_i$ for all $i\in I$; 
\item\label{item:AF}  $\psi(a_1 \cdots a_n) = 0$ whenever $\bfi =(i_1,i_2,\dots,i_n)\in \fin{I}$ and $(a_1,a_2,\dots, a_n) \in \ker \beta_{\bfi,1} \times \ker \beta_{\bfi,2} \times \cdots \times \ker \beta_{\bfi,n}$; 
\end{enumerate}
\begin{enumerate}[label=$(\gamma)$]
%\item\label{item:AF1}   $\Restr{\psi}{\cA_i} = \psi_i$ and $\Restr{\theta}{\cA_i} = \theta_i$ for all $i\in I$; 
\item\label{item:AF2}  $\theta(a_1 \cdots a_n) = 0$ whenever $\bfi =(i_1,i_2,\dots,i_n)\in \fin{I}$ and $(a_1,a_2,\dots, a_n) \in \ker \gamma_{\bfi,1} \times \ker \gamma_{\bfi,2} \times \cdots \times \ker \gamma_{\bfi,n}$. 
\end{enumerate}

The $\AP$-free product $\fmp_{i\in \{1,2\}} (\cA_i, \vp_i,\psi_i,\theta_i)$ over the index toset $I=\{1,2\} \subseteq \N$ is also denoted 
\[
(\cA_1, \vp_1 , \psi_1,\theta_1) \fmp (\cA_2,\vp_2, \psi_2,\theta_2)
\]
and $\fmp$ is regarded as a binary operation on unital 3-algebraic probability spaces.

%Note that condition \ref{item:AP} implies $\Restr{\vp}{\cA_i} = \vp_i$ for all $i\in I$. 

\item ($\AP$-freeness) Let $(\cA, \vp, \psi, \theta)$ be a unital 3-algebraic probability space. 
Let $\{\cA_i\}_{i\in I}$ be a family of subalgebras of $\cA$ containing the unit of $\cA$. Then $\{\cA_i\}_{i\in I}$ is said to be \emph{$\AP$-free}  
if conditions \ref{item:AP}, \ref{item:AF} and \ref{item:AF2} are satisfied, where $(\vp_i,\psi_i,\theta_i)$ is set to be $(\Restr{\vp}{\cA_i}, \Restr{\psi}{\cA_i},\Restr{\theta}{\cA_i})$ for all $i\in I$. Also, a family of subsets $\{S_i\}_{i\in I}$ of $\cA$ is said to be \emph{$\AP$-free} if $\{\alg{S_i}\}_{i\in I}$ is $\AP$-free, where $\alg{S_i}$ is the subalgebra of $\cA$ generated by $S_i$ and $1_\cA$.

\item\label{def:AF} ($\BGP$-free product) Let $\{(\cA_i, \psi_i,\theta_i)\}_{i\in I}$ be a family of unital 2-algebraic probability spaces. 
Then its \emph{$\BGP$-free product} $\fmp_{i\in I} (\cA_i, \psi_i,\theta_i)$ is a unital 2-algebraic probability space $(\cA, \psi,\theta)$, where $\cA = \freeproe{i\in I} \cA_i$ is the unital free product and $\psi,\theta$ are unital linear functionals satisfying conditions \ref{item:AF} and \ref{item:AF2}.

\item ($\BGP$-freeness) Let $(\cA, \psi,\theta)$ be a unital 2-algebraic probability space. 
Let $\{\cA_i\}_{i\in I}$ be a family of subalgebras of $\cA$ containing the unit of $\cA$. Then $\{\cA_i\}_{i\in I}$ is said to be \emph{$\BGP$-free} 
if conditions \ref{item:AF} and \ref{item:AF2} are satisfied, where $(\psi_i,\theta_i)$ is set to be $(\Restr{\psi}{\cA_i},\Restr{\theta}{\cA_i})$ for all $i\in I$.

\end{enumerate}
\end{defi}

\begin{rem} \label{rem:well_definedness}
The existence and uniqueness of $\AP$-free product can be proved with small modifications of the standard methods in free probability. 
For the interested reader, we include proofs in Subsection \ref{subsec:alternative}. In the presence of positivity structures, Theorem \ref{thm:reduced} also provides the proof of existence. 
\end{rem}

\begin{rem} 
If $(\cA_i,\vp_i,\psi,\theta_i)$ are unital $3$-$\ast$-algebraic probability states, so is their $\AP$-free product $(\cA, \vp, \psi, \theta)$. This is because the c-free product of states is a state \cite[Theorem 2.2]{BLS} and the $\AP$-free product is a combination of c-free products, see Proposition \ref{rem:cfree_alpha} below. 
Another proof of positivity can be given from a canonical operator model of $\AP$-freeness on the free product pre-Hilbert space, see Theorem \ref{thm:reduced}. 
%If $\psi = \theta$, then $\AP$-freeness of the subalgebras 
\end{rem}

\begin{rem} 
The $\AP$-free product can be regarded as a special case of the product constructed by Cabanal-Duvillard and Ionescu \cite{CDI97} for unital $\infty$-$\ast$-algebraic probability spaces. 
\end{rem}

An important fact is that the $\AP$-free product is a combination of the c-free product as already mentioned in  \eqref{eq:alpha_product}.

\begin{prop} \label{rem:cfree_alpha}  For unital 3-algebraic probability spaces $(\cA_i, \vp_i, \psi_i,\theta_i), i=1,2,$ we have 
\[
(\cA_1, \vp_1, \psi_1,\theta_1) \fmp (\cA_2,\vp_2, \psi_2,\theta_2)  =   (\cA_1 \freeprod \cA_2, \cf{\vp_1}{\theta_1}{\psi_2}{\vp_2}, \cf{\psi_1}{\theta_1}{\psi_2}{\psi_2}, \cf{\theta_1}{\theta_1}{\psi_2}{\theta_2}). 
\]
 \end{prop}
 \begin{proof}
This is a consequence of Lemma \ref{lem:c-free}. 
 \end{proof}

When proving an assertion for the $\AP$-free product of general unital 3-algebraic probability spaces, the following fact often makes the proof shorter because one only needs to work with $\vp$.  
\begin{prop} \label{prop:useful}
Let $\{(\cA_i, \vp_i, \psi_i, \theta_i)\}_{i\in I}$ be a family of unital 3-algebraic probability spaces with $I$ a toset and let  $(\cA, \vp, \psi, \theta)$ be the $\AP$-free product $\fmp_{i\in I} (\cA_i, \vp_i, \psi_i, \theta_i)$. 

\begin{enumerate}[label=\rm(\roman*)]
\item If $\vp_i = \psi_i$ for all $i\in I$ then $\vp=\psi$. 
\item If $\vp_i = \theta_i$ for all $i\in I$ then $\vp=\theta$. 
\end{enumerate}
\end{prop}
\begin{proof}  
Obvious from the definitions of linear functionals $\alpha_{\bfi,k}, \beta_{\bfi,k}, \gamma_{\bfi,k}$.  
\end{proof}

Some  conditions equivalent to \ref{item:AP} are presented below. 

\begin{prop}\label{prop:omega} Let $(\cA, \vp, \psi, \theta)$ be a unital 3-algebraic probability space and $\{\cA_i\}_{i\in I}$ be subalgebras of $\cA$ containing the unit of $\cA$, where $I$ is a toset with $\# I \ge 2$.  For each $i\in I$, let $\vp_i, \psi_i, \theta_i$ be unital linear functionals on the subalgebra $\cA_i$. Then condition \ref{item:AP} is equivalent to any of the following conditions: 
\begin{enumerate}[label=\rm($\AP$'),leftmargin=35pt]
\item\label{item:AP'}  $\vp(a_1 \cdots a_n) = 0$ whenever $\bfi =(i_1,i_2,\dots,i_n)\in \fin{I}$, $(a_1,a_2,\dots, a_n) \in \cA_\bfi$, $\vp_{i_n}(a_n) = 0$, $\psi_{i_k}(a_k) = 0 $ for all $k \in \Peak(\bfi)$ 
and $\theta_{i_k}(a_k) = 0$ for all $k \in \Bottom(\bfi)$.  
\end{enumerate}
\begin{enumerate}[label=\rm($\omega$),leftmargin=35pt]
\item \label{item:omega}  $\vp(a_1 \cdots a_n) = 0$ whenever $\bfi =(i_1,i_2,\dots,i_n)\in \fin{I}$ and $(a_1,a_2,\dots, a_n) \in \ker \omega_{\bfi,1} \times \ker \omega_{\bfi,2} \times \cdots \times \ker \omega_{\bfi,n}$. 
\end{enumerate}
\end{prop}

\begin{proof}
{\bf \ref{item:AP} $\Leftrightarrow$ \ref{item:AP'}.}  By symmetry it suffices to prove one implication. Assuming \ref{item:AP} we prove that condition \ref{item:AP'} holds which is the same as \ref{item:AP} except that condition $\vp_{i_1}(a_1)=0$ is replaced with condition $\vp_{i_n}(a_n)=0$. 

For $n=1$, condition \ref{item:AP'} and condition \ref{item:AP} are both equivalent to the condition that $\Restr{\vp}{\cA_i} = \vp_i$ for all $i\in I$. 
Let $n\ge2$, $\bfi =(i_1,i_2,\dots,i_n)\in \fin{I}$, $(a_1,a_2,\dots, a_n) \in \cA_\bfi$, $\vp(a_n) = 0$, $\psi(a_k) = 0 $ for all $k \in \Peak(\bfi)$ and $\theta(a_k) = 0$ for all $k \in \Bottom(\bfi)$. Repeatedly applying the decomposition $a_k = \mathring{a}_k + \vp(a_k) 1$ to the word $a_1 \cdots a_n$ from the left side we get 
\begin{equation}\label{eq:decomp_initial}
a_1 \cdots a_n
= \sum_{\ell=1}^{n-1}  \left[\prod_{j=1}^{\ell-1} \vp(a_j) \right] \mathring{a}_{\ell} a_{\ell+1} \cdots a_n + \left[\prod_{j=1}^{n-1} \vp(a_j) \right]  a_n, 
\end{equation} 
where $\prod_{j=1}^{\ell-1} \vp(a_j)$ is set to be one if $\ell=1$. 
The evaluation of the each monomial by $\vp$ is zero because the leftmost element $\mathring{a}_\ell$ has the zero $\vp$-moment, and the other $a_{\ell+1}, \dots, a_n$ satisfy the kernel condition as described in \ref{item:AP}.

\vspace{2mm}
\noindent 
{\bf \ref{item:AP'} $\Leftrightarrow$ \ref{item:omega}.}  
It is easy to see that condition \ref{item:AP'} implies condition \ref{item:omega}. For the converse statement, let us assume condition \ref{item:omega} and take  $\bfi =(i_1,i_2,\dots,i_n)\in \fin{I}$ and $(a_1,a_2,\dots, a_n) \in \cA_\bfi$  as described in condition \ref{item:AP'}. We may assume that $n\ge2$ because  otherwise conditions \ref{item:AP'} and \ref{item:omega} are trivially equivalent. 
Let us introduce the centering $a_k = \mathring{a}_k + \omega_{\bfi,k}(a_k)$ for all $k \in [n]$. In the expansion 
\begin{align}
a_1a_2 \cdots a_n 
&= (\mathring{a}_1 + \omega_{\bfi, 1} (a_1) 1_\cA) (\mathring{a}_2 + \omega_{\bfi, 2} (a_2) 1_\cA)  \cdots (\mathring{a}_n + \omega_{\bfi, n} (a_n) 1_\cA) \notag   \\
&=  \mathring{a}_1  \mathring{a}_2   \cdots \mathring{a}_n  + \omega_{\bfi, 1} (a_1)   \mathring{a}_2\mathring{a}_3   \cdots \mathring{a}_n + \omega_{\bfi, 2} (a_2)   \mathring{a}_1 \mathring{a}_3  \cdots \mathring{a}_n  +  \cdots,   \label{eq:expansion}
\end{align}
each nonzero term satisfies the assumptions described in \ref{item:omega}, and hence belongs to the kernel of $\vp$.  Note here that we already have $\omega_{\bfi, k}(a_k)=0$ for $k \in \Peak(\bfi) \cup \Bottom(\bfi)\cup\{n\}$ from the way we took $(a_1,a_2,\dots, a_n) \in \cA_\bfi$  and therefore, the centering $a_k = \mathring{a}_k + \omega_{\bfi,k}(a_k)$ matters only for $k$ not being a peak or bottom or $n$. For such numbers $k$, $a_{k-1}$ and $a_{k+1}$ belong to different subalgebras, so that the centering procedure never creates new peaks or bottoms which would be an obstruction for condition \ref{item:omega}. 
\end{proof}

%\begin{rem} %The equivalence of \ref{item:AP} and \ref{item:omega} explains the reason why 
%An attentive reader might have noticed that in the definition of $\omega_{\bfi,k}$ one can also choose an arbitrary kernel condition at $k$ which is not being peak or bottom or $n$.  A similar remark applies to $\alpha_{\bfi, k}$. 
%\end{rem}

As the main result of this section, the associativity of the $\AP$-free product and $\BGP$-free product is proved here in a generalized form.  

\begin{thm}\label{thm:main_associativity}
Let $I$ be a toset and $J,K$ be disjoint nonempty  subsets of $I$ such that $I= J \cup K$ and any element of $J$ is smaller than any element of $K$. Let $\{(\cA_i, \vp_i, \psi_i, \theta_i)\}_{i\in I}$ be a family of unital 3-algebraic probability spaces.  Then we have
\[
 \fmp_{i \in I}   (\cA_i, \vp_i, \psi_i,\theta_i) = \left( \fmp_{j \in J}   (\cA_j, \vp_j, \psi_j,\theta_j) \right) \fmp \left( \fmp_{k \in K}   (\cA_k, \vp_k, \psi_k,\theta_k) \right)  
\]
under the natural identification $\freeproe{i\in I} \cA_i \simeq (\freeproe{j\in J} \cA_j) \freeprod  (\freeproe{k\in K} \cA_k) $. 
%\end{enumerate}
\end{thm}
\begin{rem} 
This implies that the binary operation $\fmp$ (for $I=\{1,2\}$) for unital 3-algebraic probability spaces is associative. To see this,  setting $I = \{1,2,3\}$ and $J=\{1,2\}, K=\{3\}$ we get 
\[
\fmp_{i \in \{1,2,3\}}   (\cA_i, \vp_i, \psi_i,\theta_i) = \left[  (\cA_1, \vp_1, \psi_1,\theta_1) \fmp  (\cA_2, \vp_2, \psi_2,\theta_2) \right] \fmp  (\cA_3, \vp_3, \psi_3,\theta_3).   
\]
On the other hand, taking $I = \{1,2,3\}$ and $J=\{1\}, K=\{2,3\}$ yields 
\[
\fmp_{i \in \{1,2,3\}}   (\cA_i, \vp_i, \psi_i,\theta_i) = (\cA_1, \vp_1, \psi_1,\theta_1) \fmp  \left[ (\cA_2, \vp_2, \psi_2,\theta_2) \fmp  (\cA_3, \vp_3, \psi_3,\theta_3) \right].  
\]
Combining them together entails the associativity. 

\end{rem}
\begin{proof}[Proof of Theorem \ref{thm:main_associativity}] Let us denote 
\[
(\cA_X, \vp_X, \psi_X, \theta_X) = \fmp_{i \in X}   (\cA_i, \vp_i, \psi_i,\theta_i)
\]
for $X \in \{I,J,K\}$. In view of Remark \ref{rem:cfree_alpha},  what we need to prove is 
\begin{equation} \label{eq:maingoal}
(\vp_I, \psi_I, \theta_I) = \left( \cf{\vp_J}{\theta_J}{\psi_K}{\vp_K},   \cf{\psi_J}{\theta_J}{\psi_K}{\psi_K},  \cf{\theta_J}{\theta_J}{\psi_K}{\theta_K} \right). 
\end{equation}
It suffices to prove $\vp_I =  \cf{\vp_J}{\theta_J}{\psi_K}{\vp_K}$ because, thanks to Proposition \ref{prop:useful}, the other two identities $\psi_I=   \cf{\psi_J}{\theta_J}{\psi_K}{\psi_K}$ and $ \theta_I =  \cf{\theta_J}{\theta_J}{\psi_K}{\theta_K}$ correspond to the special cases $\vp_i= \psi_i$ and $\vp_i= \theta_i$ for all $i \in I$, respectively. 
It further suffices to prove 
\begin{equation}\label{eq:suffice}
\vp_I (a_1 a_2 \cdots a_n)=  (\cf{\vp_J}{\theta_J}{\psi_K}{\vp_K})(a_1 a_2 \cdots a_n)
\end{equation}
for every $(a_1,a_2,\dots, a_n) \in \cA_{\bfi}$, $n\ge1$, $\bfi =(i_1,i_2,\dots,i_n)\in \fin{I}$. 

It is obvious that \eqref{eq:suffice} holds in case $n=1$ (the common value equals $\vp_{i_1}(a_1)$). For $n\ge2$, it suffices to show that both hands sides of \eqref{eq:suffice} are zero by assuming that $a_k \in \ker \alpha
_{\bfi,k}$ for all $k \in [n]$, i.e.\ $\vp_I(a_1)=0$, $\psi_I(a_k) = 0 $ for all $k \in \Peak(\bfi)$ and $\theta_I(a_k) = 0$ for all $k \in \Bottom(\bfi)$, because of the uniqueness of $\AP$-free product. Then, by definition, $\vp_I (a_1 a_2 \cdots a_n)=0$, so that we need to prove that 
\begin{equation}\label{eq:maingoal2}
(\cf{\vp_J}{\theta_J}{\psi_K}{\vp_K})(a_1 a_2 \cdots a_n) =0. 
\end{equation}

We regroup the numbers in $\bfi$ to get 
\begin{equation}\label{eq:decomp_i}
\bfi = (\bfj_1, \bfk_1, \bfj_2, \bfk_2, \dots) \quad \text{or} \quad    \bfi = (\bfk_1, \bfj_1, \bfk_2, \bfj_2,  \dots), 
\end{equation}
where $\bfj_1,\bfj_2 ,\dots$ are elements in $\fin{J}$ and $\bfk_1, \bfk_2, \dots$ are elements in $\fin{K}$.  We denote further 
\[
\bfj_p = (i_{\alpha(p)}, i_{\alpha(p)+1},\dots, i_{\omega(p)}) \qquad \text{and} \qquad  \bfk_p = (i_{\beta(p)}, i_{\beta(p)+1},\dots, i_{\zeta(p)})
\]
 for each $p$. Let us work with the case $\bfi = (\bfj_1, \bfk_1, \bfj_2, \bfk_2, \dots)$ because the other case is almost the same. 

We prove below that 
\begin{align} 
&\vp_J(a_{\alpha(1)} \cdots a_{\omega(1)}) =0, \label{eq:key1} \\
 &\theta_J(a_{\alpha(p)} \cdots a_{\omega(p)})=0 \quad \text{for all} \quad  p\ge2, \label{eq:key2} \\  
 &\psi_K (a_{\beta(q)} \cdots a_{\zeta(q)})=0 \quad \text{for all} \quad  q\ge1  \label{eq:key3}
\end{align}
except for the last  subsequence. Once \eqref{eq:key1}--\eqref{eq:key3} are proved then by the definition of $\AP$-free product for the binary case, the desired \eqref{eq:maingoal2} holds.

We first work on \eqref{eq:key2}. If $\# J =1$ then $\bfj_p$ is a singleton for every $p\ge2$ and $\alpha(p)=\omega(p)$ and $\theta_J(a_{\alpha(p)})=0$ because $\alpha(p)$ is always a bottom (unless $\bfi$ ends with $\bfj_p$). One may therefore assume that $\# J \ge2$. 
Note that a number in $\{\alpha(p)+1, \dots, \omega(p)-1\}$ is a peak (resp.\ bottom) with respect to the sequence $\bfi$ if and only if it is a peak (resp.\ bottom) with respect to the subsequence $\bfj_p$. If, moreover, $\alpha(p)$ is a bottom with respect to $\bfi$, then our assumption includes $\theta_J(a_{\alpha(p)})=0$, and hence $\theta_J(a_{\alpha(p)} \cdots a_{\omega(p)})=0$ holds directly by the definition of $\BGP$-free product.  
If $\alpha(p)$ is not a bottom with respect to $\bfi$, because $\bfj_p$ is preceded by $\bfk_{p-1}$ in the sequence $\bfi$,  there exists $r \ge 1$ such that 
\begin{align*}
&K \ni i_{\alpha(p)-1} > i_{\alpha(p)} > i_{\alpha(p)+1} > \cdots > i_{\alpha(p)+r} < i_{\alpha(p)+r+1} \in J \cup K \quad \text{or} \quad \\
 &K \ni i_{\alpha(p)-1} > i_{\alpha(p)} > i_{\alpha(p)+1} > \cdots > i_{\alpha(p)+r} =i_{\omega(p)} = i_n,  
\end{align*}
The latter case occurs only when the sequence $\bfi$ ends with its subsequence $\bfj_p$ and hence need not be handled (recall that we need not prove \eqref{eq:key2} for the last subsequence of $\bfi$). In the first case $\alpha(p)+r$ is contained in $\Bottom(\bfi)$, so that we already assumed $\theta_J(a_{\alpha(p)+r})=0$. Repeatedly applying the decomposition $a_j = \mathring a_j + \theta_J(a_j)1$ to the word $a_{\alpha(p)} a_{\alpha(p)+1} \cdots a_{\alpha(p)+r-1}$ from the left side, we arrive at the identity
\begin{align}
a_{\alpha(p)} \cdots a_{\alpha(p)+r-1} 
&= \mathring a_{\alpha(p)}  a_{\alpha(p)+1} \cdots a_{\alpha(p)+r-1}  + \sum_{\ell=1}^{r-2}\left[ \prod_{j=0}^{\ell-1} \theta_J(a_{\alpha(p)+j}) \right] \mathring a_{\alpha(p)+\ell} a_{\alpha(p)+\ell+1} \cdots a_{\alpha(p)+r-1}  \notag \\
&\quad+ \left[ \prod_{j=0}^{r-2} \theta_J(a_{\alpha(p)+j}) \right] \mathring a_{\alpha(p)+r-1} + \left[ \prod_{j=0}^{r-1} \theta_J(a_{\alpha(p)+j}) \right] 1,  \label{eq:decomp_alpha}
\end{align}
which gives a decomposition of $a_{\alpha(p)} \cdots a_{\omega(p)}$ into the sum of words in the kernel of $\theta_J$ because the leftmost element $\mathring a_j$ (or $a_{\alpha(p)+r}$) in each word belongs to the kernel of $\theta_J$. 
Altogether, claim \eqref{eq:key2} holds.  By a similar argument \eqref{eq:key3} also holds true. Claim \eqref{eq:key1} follows by the definition of $\AP$-free product because we already assumed $\vp_J(a_1)=0$. 
\end{proof}

%We note that the \AF-product generalizes free, monotone and antimonotone products. 

%%%%%%%%%%%%%%%%%%%%%%%%%%%%%%%%%%%%%%%%%%%%%%%%%%%%%%%%%%%%%%%%%%
\subsection{The uniqueness and existence of the $\AP$-free product} \label{subsec:alternative}

This subsection shows the uniqueness and existence of the $\AP$-free product and thus complements Definition \ref{def:alpha_indep}, cf.\ Remark \ref{rem:well_definedness}. 
The uniqueness is a consequence of the following decomposition result.

 \begin{lem}\label{lem:decomp} Let $\cA$ be a unital algebra generated by a family of subalgebras $\{\cA_i\}_{i\in I}$ of $\cA$ containing the unit of $\cA$. 
For each $\bfi =(i_1,i_2,\dots, i_n)\in\fin{I}$ and $k\in[n]$, let $\pp_{\bfi,k}$ be a (possibly nonunital) linear functional on $\cA_{i_k}$. 
Then any element of $\cA$ is of the form 
\begin{equation}\label{eq:elementary}
\lambda 1_{\cA}+ \sum_{\bfi \in \fin{I}}a_{\bfi,1} a_{\bfi,2}\cdots a_{\bfi,|\bfi|}, 
\end{equation}
where $\lambda\in \comp, (a_{\bfi,1} ,a_{\bfi,2},\cdots, a_{\bfi,|\bfi|}) \in \ker p_{\bfi,1}  \times \ker p_{\bfi,2} \times \cdots \times \ker p_{\bfi,|\bfi|}, \bfi=(i_1,i_2,\dots, i_n)\in\fin{I}$. 
\end{lem}
\begin{proof}
It suffices to prove the statement for words. There is nothing to do with the words of length $0$. Suppose that the desired decomposition exists for all words of length less than $\ell-1$. Let us take a word $a_1 a_2 \cdots a_\ell$, where $(a_1,\dots, a_\ell) \in \cA_\bfj$, $\bfj \in \fin{I}, |\bfj|=\ell$. 
Setting $\mathring{a}_k= a_k-\lambda_k1_{\cA}$, $\lambda_k = \pp_{\bfj,k} (a_k)$ for $k \in[\ell]$, we have the formula 
\begin{align}
a_1a_2 \cdots a_\ell &= (\mathring{a}_1 +\lambda_1 1_{\cA})(\mathring{a}_2+ \lambda_2 1_{\cA})\cdots (\mathring{a}_n+ \lambda_n 1_{\cA})  \notag  \\
   &= \mathring{a}_1 \mathring{a}_2\cdots \mathring{a}_\ell  
     + \sum_{S \subseteq [\ell], S\neq \emptyset }\left[\prod_{s \in S} \pp_{\bfj,s}(a_{s}) \right] \mathring{a}_{S_0} \mathring{a}_{S_1} \cdots   \mathring{a}_{S_{k}},    \label{eq:decomp0} 
\end{align}
where the following notation is employed above: $S=\{s_1,s_2,\dots, s_k\}$ with $s_1 < \cdots < s_k$, $S_i = \{s_{i}+1,\dots, s_{i+1}-1\}$, $s_0:=0, s_{k+1}:=\ell+1$ and $\mathring{a}_{T}:= \mathring{a}_{t_1} \cdots \mathring{a}_{t_p}$ whenever $T=\{t_1,t_2,\dots, t_p\} \subset \N$ with $t_1 < \cdots < t_p$, and $\mathring{a}_{\emptyset}=1_\cA$. By the induction hypothesis each term $ \mathring{a}_{S_0} \mathring{a}_{S_1} \cdots   \mathring{a}_{S_{k}}$ can be decomposed as stated, so that we are done. 
\end{proof}

Lemma \ref{lem:decomp} with $\pp_{\bfi,k}:=\AP_{\bfi,k}$ guarantees the uniqueness of the $\AP$-free product in Definition \ref{def:alpha_indep}: if $\vp, \vp'$ satisfy condition \ref{item:AP} then choosing a decomposition of each $a \in \cA$ as in  Lemma \ref{lem:decomp} implies $\vp(a)=\lambda = \vp'(a)$.

Lemma \ref{lem:decomp} is not useful for proving the existence because the decomposition might not be unique. A key to existence is the following construction of the unital free product that slightly generalizes $\freeproe{i\in I} \cA_i$ defined in \eqref{eq:free_product}.  
%This can be regarded as a refinement of Lemma \ref{lem:decomp}. 

Let $I$ be a set and let $\{\cA_i\}_{i\in I}$ be a family of unital algebras (over $\comp$). 
For each $\bfi =(i_1,i_2,\dots, i_n)\in\fin{I}$ and $k\in[n]$, let $\pp_{\bfi,k}$ be a unital linear functional on $\cA_{i_k}$. 
We define the vector space 
\[
\cA^{\pp}:= \comp 1 \oplus \bigoplus_{\bfi \in \fin{I}} ( \ker \pp_{\bfi,1}  \otimes \ker \pp_{\bfi,2} \otimes \cdots \otimes \ker \pp_{\bfi, |\bfi|} ). 
\]
This vector space can be made into a unital algebra with unit $1$ and with an associative product defined as follows. For $\bfi = (i_1,i_2,\dots, i_n),\bfj=(j_1,j_2,\dots, j_m) \in \fin{I}$,  $(a_1,a_2,\dots, a_n) \in \ker \pp_{\bfi,1}  \otimes \ker \pp_{\bfi,2} \otimes \cdots \otimes \ker \pp_{\bfi,n}$ and $(b_1,b_2, \dots, b_m) \in \ker \pp_{\bfj,1}  \otimes \ker \pp_{\bfj,2} \otimes \cdots \otimes \ker \pp_{\bfj,m}$ we define 
\begin{align*}
   & (a_1 \otimes \cdots \otimes a_n)(b_1\otimes\cdots \otimes b_m) \\
 &:=   \begin{cases}
        a_1 \otimes \cdots \otimes a_n\otimes b_1\otimes\cdots \otimes b_m & \text{if }i_n\neq j_1, \\
         a_1 \otimes \cdots \otimes a_{n-1}\otimes (a_n b_1 - \lambda 1_{\cA_{i_n}}) \otimes b_2 \otimes \cdots \otimes b_m  + \lambda (a_1 \otimes \cdots \otimes a_{n-1})(b_2 \otimes \cdots \otimes b_m)& \text{if } i_n  = j_1, 
    \end{cases}
\end{align*} 
 where $\lambda$ is selected to be $\pp_{(i_1,i_2,\dots, i_{n}, j_2, j_3,\dots, j_m), n}(a_nb_1)$. This is a recursive definition:  if $i_n=j_1$ and furthermore $i_{n-1} = j_2$, then we need to again apply the centering procedure for $a_{n-1}b_2$ until we get a linear combination of words that fit the direct sum of $\cA^p$. 
 The algebra $\cA^p$ satisfies the universality of the unital free product and hence is isomorphic to $\freeproe{i\in I} \cA_i$ as unital algebras.

We are ready to construct the $\AP$-free product. Let $\{(\cA_i, \vp_i, \psi_i,\theta_i)\}_{i\in I}$ be a family of unital 3-algebraic probability spaces with $I$ a toset. Let us select the unital linear functional $\pp_{\bfi,k}:=\omega_{\bfi,k}$. On the unital free product $\cA^\omega$ a linear functional $\vp$ can be simply defined as 
\[
\varphi(\lambda 1 + a):= \lambda \quad \text{for all $\lambda \in \comp$ and~} a \in \bigoplus_{\bfi \in \fin{I}}( \ker \omega_{\bfi,1}  \otimes \ker \omega_{\bfi,2} \otimes \cdots \otimes \ker \omega_{\bfi, |\bfi|} ). 
\]
Then it satisfies condition \ref{item:omega} which is equivalent to \ref{item:AP}. The $\BGP$-free product $(\psi, \theta)$ can be defined as special cases of $\vp$ according to Proposition \ref{prop:useful}. The existence of $\AP$-free product is thus justified. 

One can also change the definition of the linear functionals $p_{\bfi,k}$. For example, the linear functionals 
\begin{equation}\label{eq:alpha2}
\pp_{\bfi,k} = 
\begin{cases} \varphi_{i_1} & \text{if~} k=1, \\
\psi_{i_k} & \text{if~} k \in \Peak(\bfi), \\
\theta_{i_k} & \text{if~} k \in \Bottom(\bfi), \\
\vp_{i_k}& \text{otherwise} 
\end{cases}
\quad (n\ge2), \quad \text{and} \quad  \pp_{(i_1),1}=\vp_{i_1}\quad (n=1),  
  \end{equation}
 give rise to the same linear functional $\vp$ on the unital free product algebra. Basically one needs to define $\pp_{\bfi,k}$ equal to $\AP_{\bfi,k}$ for $k \in \Peak(\bfi) \cup \Bottom(\bfi)$ and needs to define $\pp_{\bfi,1} = \vp_{i_1}$ or $\pp_{\bfi,n}= \vp_{i_n}$, but for the other $k$'s any unital linear functional can be chosen.

%%%%%%%%%%%%%%%%%%%%%%%%%%%%%%%%%%%%%%%%%%%%%%%%%%%%%%%%%%%%%%%%%%%%%%%%%%%%%%%%%%%%%%%%%%%%%%%%%%%%%%%%%%%%%%%%%%%%%%%%%%%%%%%%%%%%%%%%%%%%%%%%%%%%%%%%%%%%
\section{The $\AP$-freeness and the free product Hilbert space}\label{rep}

We propose a canonical construction of $\AP$-free operators on the free product of pre-Hilbert spaces. This construction is based on \cite{CDI97} and generalizes constructions in \cite{Avi,BS,Mur3,Pop,V1}. A crucial role is played by certain invariant subspaces of the free product  of pre-Hilbert spaces for a canonical operator acting from the left side.

%%%%%%%%%%%%%%%%%%%%%%%%%%%%%%%%%%%%%%%%%%%%%%%%%%%%%%%%%%%%%%%%%%%%%
\subsection{The free product pre-Hilbert space and its subspaces}\label{subsec:Hilbert}

A pair $(H,\xi)$ of a pre-Hilbert space $H$ and its unit vector $\xi$ is called a pointed pre-Hilbert space. As a convention, an inner product is linear on the right component and antilinear on the left. 
The set of adjointable linear operators on a pre-Hilbert space $H$ is denoted $\Lin(H)$.   
The orthogonal projection from $H$ onto its closed subspace $K$ is denoted $P_K$ if it exists. Note that the orthogonal projection onto a finite-dimensional subspace always exists. For further details, see \cite{GHU}.   

Let $I$ be a toset and  $\{(H_i,\xi_i)\}_{i\in I}$ be a family of pointed pre-Hilbert spaces. 
Let $H_i ^0 := H_i \ominus \comp\xi_i$, i.e.\ $H_i$ has the orthogonal decomposition $H_i = \comp\xi_i \oplus H_i^0$. 
Let $(H, \xi)$ be the free product of $\{(H_i, \xi_i)\}_{i\in I}$, denoted $\freeproe{i\in I}(H_i,\xi_i)$, 
\[
\begin{split}
&H= \comp \xi \oplus \bigoplus_{n = 1}^\infty \bigoplus _{i_1 \neq \cdots \neq  i_n} H_{i_1} ^0 \otimes \cdots \otimes H_{i_n}^0, \\ 
\end{split}
\]
%Let $\pi, \sigma: \freeprod \cA_k \to  \Lin(H)$ be the $\freeprod$-representations defined as the extensions of $\pi_k$ and $\sigma_k$, respectively.  
and $H(i), i \in I,$ be complementary subspaces 
\[
H(i)=\comp \xi \oplus \bigoplus_{n = 1}^\infty \bigoplus _{\substack{i_1 \ne \cdots \ne  i_n, \\ i_1 \neq i}} H_{i_1} ^0 \otimes \cdots \otimes H_{i_n}^0. 
\]  
For each $i \in I$ a canonical isomorphism $U_i \colon H \to H_i \otimes H(i)$ is defined by
\begin{align*}
&U_i(\xi) = \xi_i \otimes \xi, \\
&U_i(f) = f \otimes \xi,  \qquad f \in H_i^0,  \\
&U_i(f_1\otimes \cdots \otimes f_n) = 
\begin{cases}
f_1 \otimes (f_2\otimes \cdots \otimes f_n), ~~~ i_1 = i, ~n \geq 2, \\
\xi_i \otimes  (f_1\otimes \cdots \otimes f_n), ~~~ i_1 \neq i,  
\end{cases}
\end{align*}
where $f_1 \otimes \cdots \otimes f_n \in  H_{i_1} ^0 \otimes \cdots \otimes H_{i_n}^0$ with $\bfi =(i_1,i_2,\dots, i_n) \in \fin{I}$. 

Let $\lambda_i\colon \Lin(H_i) \to \Lin(H)$ be the $\ast$-homomorphism  defined by 
\[
\lambda_i(A_i):=U_i^\ast (A_i \otimes \id_{H(i)})U_i.
\]
The operator $\lambda_i(A_i)$ ($A_i \in \Lin(H_i)$) acting on elementary vectors are given by 
\[
\begin{split}
\lambda_i(A_i)\xi  
&=  \langle  \xi_i, A_i \xi_i \rangle \xi + P_{H_i^0}(A_i \xi_i), \\
\lambda_i(A_i)(f_1\otimes \cdots \otimes f_n)  
&= 
\begin{cases}
 \langle \xi_i, A_i f_1 \rangle f_2\otimes \cdots \otimes f_n + P_{H_i^0}(A_i f_1)\otimes f_2 \otimes \cdots \otimes f_n, & \text{if~} i_1 = i, ~n \geq 2, \\
 \langle \xi_i, A_i f_1 \rangle \xi + P_{H_i^0}(A_i f_1), & \text{if~}   i_1 = i, ~n =1, \\ 
 \langle \xi_i, A_i \xi_i \rangle f_1\otimes \cdots \otimes f_n + P_{H_i^0}(A_i \xi_i)\otimes f_1 \otimes \cdots \otimes f_n, &\text{if~}  i_1 \neq  i,  
 \end{cases}
\end{split}
\]
where $f_1 \otimes \cdots \otimes f_n \in  H_{i_1} ^0 \otimes \cdots \otimes H_{i_n}^0$. 
Then the $\ast$-subalgebras $\{\lambda_i(\Lin(H_i))\}_{i\in I}$ of $\Lin(H)$ is free in the unital algebraic probability space $(\Lin(H), \langle \xi, \cdot~\xi\rangle)$.

To construct the $\AP$-free product, we introduce the two subspaces of $H$
\[
\nA(i) =  H_<(i) \oplus H_{=, <}(i)\quad \text{and} \quad \nB(i) = H_>(i) \oplus H_{=, >}(i), 
\]
where 
\[
\begin{split}
&H_<(i) =  \bigoplus_{n = 1}^\infty \bigoplus _{\substack{i_1 \neq \cdots \neq  i_n, \\ i_1 < i}} H_{i_1} ^0 \otimes \cdots \otimes H_{i_n}^0, \\
&H_{=,<}(i) =  \bigoplus_{n = 2}^\infty \bigoplus _{\substack{i_1 \neq \cdots \neq  i_n, \\ i_1 = i, ~i_2 < i}} H_{i_1} ^0 \otimes \cdots \otimes H_{i_n}^0, \\ 
&H_>(i) =  \bigoplus_{n = 1}^\infty \bigoplus _{\substack{i_1 \neq \cdots \neq  i_n, \\ i_1 > i}} H_{i_1} ^0 \otimes \cdots \otimes H_{i_n}^0, \\
&H_{=, >}(i) =  \bigoplus_{n = 2}^\infty \bigoplus _{\substack{i_1 \neq \cdots \neq  i_n, \\ i_1 = i, ~i_2 > i}} H_{i_1} ^0 \otimes \cdots \otimes H_{i_n}^0. 
\end{split}
\]
Observe that $\nA(i)$ and $\nB(i)$ are invariant regarding the actions of $\lambda_i$, and  that the decomposition 
\[
H = H_i \oplus \nA(i) \oplus \nB(i)
\]
holds  under the identification of $\comp \xi \oplus H_i^0$ with $H_i$.

%%%%%%%%%%%%%%%%%%%%%%%%%%%%%%%%%%%%%%%%%%%%%%%%%%%%%%%%%%%%%%%%%%%%%%%%%%%%
\subsection{The existence of $\AP$-free product: the second approach}

Throughout this subsection, let $I$ be a toset and $\{(\cA_i, \vp_i, \psi_i, \theta_i)\}_{i\in I}$ be a family of \emph{unital 3-$\ast$-algebraic probability spaces} because we take GNS representations.  
Let $\pi_i, \sigma_i, \rho_i\colon \cA_i \to \Lin(H_i)$ be unital $\ast$-representations on a pointed pre-Hilbert space $(H_i,\xi_i)$ for each $i\in I$ satisfying 
$\pi_i(a_i) = \langle\xi_i,  \pi_i(a_i)\xi_i\rangle$,  $\psi_i(a_i) = \langle \xi_i, \sigma_i(a_i)\xi_i\rangle$ and 
$\theta_i(a_i) = \langle \xi_i, \rho_i(a_i)\xi_i\rangle$.  Note that such $\ast$-representations exist. For example, take GNS representations 
$(\mathring{\pi}_i, H_i^\pi, \xi_i^\pi)$, $(\mathring{\sigma}_i, H_i^\sigma, \xi_i^\sigma)$ and $(\mathring{\rho}_i,H_i^\rho,\xi_i^\rho)$ for the unital algebraic probability spaces $(\cA_i, \vp_i)$, $(\cA_i,\psi_i)$ and $(\cA_i, \theta_i)$, respectively. On the larger space $(H_i,\xi_i) := (H_i^\pi\otimes H_i^\sigma \otimes H_i^\rho, \xi_i^\pi \otimes \xi_i^\sigma \otimes \xi_i^\rho)$, unital $\ast$-representations $\pi_i,\sigma_i,\rho_i\colon \cA_i\to \Lin(H_i)$ can be constructed by 
\begin{align*}
&\pi_i(a_i):= \mathring{\pi}_i(a_i) \otimes \id \otimes \id, \\
&\sigma_i(a_i) := \id\otimes \mathring{\sigma}_i(a_i)\otimes \id,  \\ 
&\rho_i(a_i) := \id \otimes \id \otimes \mathring{\rho}_i(a_i).  
\end{align*}
Let $H$ be the free product of $\{(H_i,\xi_i)\}_{i\in I}$ and  $\nA(i), \nB(i)$ be its invariant subspaces regarding the $\ast$-homomorphism $\lambda_i\colon \Lin(H_i)\to \Lin(H)$ as constructed in Subsection \ref{subsec:Hilbert}. 
The main role is played by the unital $\ast$-homomorphisms $J_i^\AP,J_i^{\BP},J_i^\GP \colon \cA_i \to \Lin(H)$ defined by  
\begin{align}
& J^{{\AP}}_i(a_i) = \lambda_i(\pi_i(a_i))P_{H_i} \oplus \lambda_i(\sigma_i(a_i))P_{\nA(i)} \oplus \lambda_i(\rho_i(a_i))P_{\nB(i)},  \\
& J^{\BP}_i(a_i) =\lambda_i(\sigma_i(a_i))P_{H_i \oplus \nA(i)} \oplus \lambda_i(\rho_i(a_i))P_{\nB(i)}, \\
& J^{\GP}_i(a_i) = \lambda_i(\sigma_i(a_i))P_{\nA(i)} \oplus  \lambda_i(\rho_i(a_i))P_{H_i \oplus \nB(i)}.
\end{align}
Note that these are $\ast$-homomorphisms because $\lambda_i(\Lin(H_i))$ preserves the subspaces $H_i, \nA(i)$ and $\nB(i)$ for all $i\in I$. Apparently,  $J_i^{{\AP}}$ coincides with $J_i^{\BP}$ (resp.\ $J_i^{\GP}$) in case $(\vp_i, \pi_i) = (\psi_i, \sigma_i)$ (resp.\ $(\vp_i, \pi_i)=(\theta_i,\rho_i)$).

By the universality of the unital free product, a 
unique $\ast$-homomorphism 
\[
J^{{\AP}}  \colon \freeproe{i\in I} \cA_i \to \Lin(H)  \quad \text{with} \quad \Restr{J^{{\AP}}}{\cA_i} = J_i^\AP \quad \text{for all} \quad i \in I 
\]
exists. Similarly, we set $\ast$-homomorphisms $J^{\BP}, J^{\GP}\colon \freeproe{i\in I} \cA_i \to \Lin(H)$.

\begin{thm} \label{thm:reduced}
%Let $(\freeproe{i\in I} \cA_i,\vp, \psi, \theta)$ be the $\AP$-free product of $\{(\cA_i,\vp_i, \psi_i, \theta_i)\}_{i\in I}$.  
In the setting above, the triplet of states $(\vp, \psi, \theta)$ defined by 
\[
\vp(a) := \langle \xi, J^{{\AP}}(a)\xi \rangle, \quad  \psi(a)  := \langle \xi,  J^{\BP}(a)\xi \rangle  \quad \text{and} \quad 
\theta(a)  := \langle \xi, J^{\GP}(a)\xi \rangle   \quad \text{for all}\quad a \in \freeproe{i\in I} \cA_i
\]
is the $\AP$-free product. 
\end{thm}     
\begin{proof}

Let $\omega_{\bfi,k}$ be as defined in \eqref{eq:omega}. 
We check condition \ref{item:omega}.  Suppose that $(a_1, a_2,\dots, a_n) \in \ker \omega_{\bfi,1}  \times \ker \omega_{\bfi,2} \times \cdots \times \ker \omega_{\bfi,n}$. Straightforward calculations yield 
\[
 J^{{\AP}}(a_1 \cdots a_n)\xi = (\tau_{\bfi,1}(a_1)\xi_{i_1}) \otimes \cdots \otimes (\tau_{\bfi,n}(a_n)\xi_{i_n}) \in H_{i_1}^0 \otimes \cdots \otimes H_{i_n}^0,   
\]
where $\tau_{\bfi,k}$ is defined by 
\begin{equation}\label{eq:tau}
\tau_{\bfi,k} = 
\begin{cases} 
\pi_{i_n} & \text{if~} k=n, \\
\sigma_{i_k} & \text{if~} k \in \Des(\bfi), \\
 \rho_{i_k} & \text{if~} k \in \Asc(\bfi), 
\end{cases}
\quad (n\ge1) \quad \text{and} \quad \tau_{(i_1),1} = \pi_{i_1} \quad (n=1). 
\end{equation}
Therefore, $\langle \xi, J^{\AP}(a_1a_2 \cdots a_n)\xi \rangle = 0$ and hence condition \ref{item:omega} is satisfied as desired. Condition \ref{item:AP} is also the case according to Proposition \ref{prop:omega}. 

Conditions \ref{item:AF} and \ref{item:AF2} for the states $\psi$ and $\theta$ follow by setting $(\vp_i,\pi_i)= (\psi_i,\sigma_i)$ for all $i \in I$ and $(\vp_i,\pi_i)= (\theta_i,\rho_i)$ for all $i \in I$, respectively. 
\end{proof}

%\begin{rem} 
%Because $J^{{\AP}}, J^{\BP}, J^{\GP}$ are unital $\ast$-homomorphisms, Theorem \ref{thm:reduced} implies that $\vp, \psi, \theta$ are all states on $ \freeproe{i \in I} \cA$.  
%\end{rem}

%%%%%%%%%%%%%%%%%%%%%%%%%%%%%%%%%%%%%%%%%%%%%%%%%%%%%%%%%%%%%%%%%%%%%%%%%%%%%%%%%%%%%%%%%%%%%%%%%%%%%%%%%%%%%%%%%%%%%%%%%%%%%%%%%%%%%%%%%%%%%%%%%%%%%%%%%%%%
\section{The $\AP$-free cumulants}\label{cum}

In this section we define $\AP$-free cumulants for $\AP$-freeness. The main result of this section is the moment-cumulant formula for $\AP$-free cumulants  (Theorem \ref{thm:m-c-alpha}). 
The proof is based on the general theory of cumulants for spreadability systems \cite{HL17}. The free, monotone, antimonotone and Boolean cumulants are then unified by the $\AP$-free cumulants as a consequence of transfer of associativity. 

Subsections \ref{subsec:partition} and \ref{subsec:orderedpartition} collect several notions on unordered and ordered set partitions, and then Subsection \ref{subsec:spreaded_cumulants} reminds us of spreadability systems and the associated cumulants developed in \cite{HL17}. 
With these preparatory materials, we work on the highest coefficients for $\AP$-freeness in Subsection \ref{subsec:highest}, which is a key ingredient in the proof of the moment-cumulant formula. 

Throughout this section the following notation is useful for  simplifying formulas.  

\begin{Notation}
Let $\cA$ be an algebra and $\vp\colon \cA\to \comp$ be a linear functional and $M=(M_n)_{n\in\N}$ be a sequence of multilinear functionals $M_n\colon \cA^n \to\comp$.  Then for every $P =\{p_1, p_2,\dots, p_k\} \subseteq [n]$ with $p_1 < p_2 < \cdots < p_k$ and every $a_1, \dots, a_n \in \cA$ we set 
\[
\vp(a_P) := \vp (a_{p_1} a_{p_2} \cdots a_{p_k}) \quad \text{and } \quad M_{k}(a_P) := M_k(a_{p_1} ,a_{p_2}, \dots, a_{p_k}). 
\]
\end{Notation}

%%%%%%%%%%%%%%%%%%%%%%%%%%%%%%%%%%%%%%%%%%%%%%%%%%%%%%%%%%%%
\subsection{Set partitions}\label{subsec:partition}

Let $L$ be a totally ordered finite set below. 
A \emph{set partition} of $L$ is a set $\pi=\{P_1,P_2,\dots, P_k\}$, whose elements $P_1,\dots, P_k$ are disjoint nonempty subsets of $L$ such that $\cup_{i=1}^k P_i =L$. The number $k$ is denoted $|\pi|$. Each $P_i$ is called a block of $\pi$. The set of the set partitions of $L$ is denoted $\SP(L)$. In particular $\SP([n])$ is abbreviated to $\SP(n)$. We also introduce the notation 
\[
\SP = \bigcup_{n\in\N} \SP(n). 
\]  
Similar notation will be used later for other classes without explicitly defining, e.g.\ $\NC(n)$, $\IP(n),$ $\OSP,$ $\ONCP$. 

For set partitions $\pi,\sigma$ of $L$ a binary relation $\pi \le \sigma$ means that for every $P\in \pi$ there exists $S \in \sigma$ such that $P \subseteq S$. This makes $\SP(L)$ a poset. Regarding this partial order, the set partition $1_L := \{L\}$ is a maximum. The notation $1_{[n]}$ is simplified to $1_n$. 

A set partition is said to be crossing if it has two blocks $B, B'$ such that some elements $i,j\in B$ and $k,\ell \in B'$ satisfy $i < k < j <\ell$. Otherwise, a set partition is said to be \emph{noncrossing}. The set of the noncrossing set partitions of $L$ is denoted $\NC(L)$. A set partition of $L$ is said to be \emph{interval} if every block is an interval, i.e.\ of the form $\{i, i+1, \dots, j\}$, where $ i \le j$. The set of the interval set partitions of $L$ is denoted $\IP(L)$.

Let $\pi$ be a noncrossing partition. For blocks $P,P' \in \pi$, the notation $P \gtrdot P'$ means that  $P'$ \emph{covers} $P$, i.e.\  there exist $i,j \in P'$ such that $i < k < j$ for all $k \in P$. A block of $\pi$ which is covered by another block of $\pi$ is called an \emph{inner block}; otherwise it is called an \emph{outer block}. The set of the inner blocks of $\pi$ is denoted $\Inn(\pi)$ and the set of the outer blocks of $\pi$ is denoted $\Out(\pi)$. For an inner block $P$ of $\pi$, there is a unique block, denoted $\Cov(P)$ and called the \emph{nearest cover}, that covers $P$ and there is no block $P'$ such that $P \gtrdot P' \gtrdot C(P)$.

\subsection{Ordered set partitions and colored set partitions}\label{subsec:orderedpartition}
Let $L$ be a totally ordered finite set.  An \emph{ordered set partition} 
 of $L$ 
is a sequence 
$\pi=(P_1, P_2, \ldots,P_p)$ of subsets of $L$ such that $\bar\pi:= \{P_1, P_2, \ldots,P_p\}$ 
is a set partition of $L$. 
The set of ordered set partitions of $L$ is denote by $\OSP(L)$. The particular case $\OSP([n])$ is abbreviated to $\OSP(n)$. 

As a more general notion, an \emph{$m$-colored set partition} of $L$ is a pair $(\pi, f)$ of a set partition $\pi$ of $L$ and a function $f\colon \pi \to [m]$.  This can be regarded as an ordered set partition when $f$ is injective. The set of the $m$-colored set partitions of $[n]$ is denoted $\CSP^m(n)$.

Some terms and notations for ordinary set partitions are also used for ordered set partitions when no confusion can arise.
For example, the notation $P\in \pi$  indicates that $P$ is a block of $\bar{\pi}$ and we also call $P$ a block of $\pi$. 

For an ordered set partition $\pi=(P_1, P_2, \ldots,P_p)$, we say that $P_i$ \emph{precedes} $P_j$ if $i<j$. For  $\pi=(P_1, P_2,\dots,P_p), \sigma=(S_1, S_2,\dots,S_s) \in \OSP(L)$, a binary relation $\pi\leq \sigma$ means that: 
\begin{itemize}
\item 
 $\bar{\pi}\leq\bar{\sigma}$ as set partitions;   
\item if $i<j$, $P_i \subseteq S_k$ and $P_j \subseteq S_l$ then $k \leq l$. 
\end{itemize}
%In other words, $\pi\leq\sigma$ if every block of $\sigma$ is a union of
This binary relation makes $(\OSP(L),\leq)$ a poset. With a slight abuse of notation, we denote by $1_L =(L)$ the maximum element of $\OSP(L)$ regarding the partial order $\le$. We also abbreviate $1_{[n]}$ to $1_n$.

\begin{figure}[t]
\begin{center}
\begin{tikzpicture}[scale=0.7]
\draw[-] (3,0.5) -- (3,1.5) -- (4,1.5) --(4,0.5);
\draw[-] (5,0.5) -- (5,1.5);
\draw[-] (8,0.5) -- (8,1.5) -- (9,1.5) --(9,0.5);
\draw[-] (2,0.5) -- (2,2.5) -- (6,2.5) --(6,0.5);
\draw[-] (6,2.5) -- (11,2.5) --(11,0.5);
\draw[-] (7,0.5) -- (7,2) -- (10,2) --(10,0.5);
\draw[-] (1,0.5) -- (1,3) -- (12,3) --(12,0.5);
\draw[-] (13,0.5) -- (13,3) -- (14,3) --(14,0.5);
\node at (1,0) {1};
\node at (2,0) {2};
\node at (3,0) {3};
\node at (4,0) {4};
\node at (5,0) {5};
\node at (6,0)  {6}; 
\node at (7,0)  {7}; 
\node at (8,0) {8};
\node at (9,0) {9};
\node at (10,0) {10};
\node at (11,0) {11};
\node at (12,0)  {12}; 
\node at (13,0)  {13}; 
\node at (14,0) {14};
%%%%%
\node at (8.5,2.2) {\tiny1};
\node at (3.5,1.7) {\tiny 2};
\node at (13.5,3.2)  {\tiny 3}; 
\node at (6.5,3.2)  {\tiny4}; 
\node at (8.5,1.7) {\tiny 5};
\node at (5,1.7)  {\tiny 6}; 
\node at (6.5,2.7) {\tiny 7};
\end{tikzpicture}
\end{center}
\caption{The diagram of $\pi = (P_1, \cdots, P_7) \in \ONCP(14)$, where $P_1 = \{7,10 \}, P_2=\{3, 4 \}, P_3 = \{13, 14 \}, P_4 = \{1, 12 \}, P_5 = \{8, 9 \}, 
P_6 = \{5 \}$ and $P_7 = \{2, 6, 11 \}$. The set of outer blocks is $\Out(\pi) = \{P_3,P_4\}$ and the set of inner blocks is $\Inn(\pi) = \{P_1,P_2,P_5,P_6,P_7\}$. }
\label{fig:ordered_nc}

\end{figure}

For an ordered set partition $\sigma=(S_1, S_2, \dots,S_s) \in \OSP(L)$ 
  and a nonempty subset
  $P\subseteq [n]$, the \emph{restriction} $\Restr{\sigma}{P}\in\OSP(P)$ 
  is the ordered set partition $(P \cap S_1,P\cap S_2, \dots,P\cap S_s)\in
  \OSP(P)$, where the empty sets, if appear, are removed. 

 The \emph{quasi-meet} of $\pi =(P_1, P_2,\ldots,P_p)$ and $\sigma=(S_1,S_2,\ldots,S_s) \in \OSP(L)$ is the ordered set partition of $L$
 \[
\pi \curlywedge \sigma  := (P_1 \cap S_1, P_1\cap S_2, \ldots, P_1 \cap S_s, P_2 \cap S_1, P_2\cap S_2,\ldots, P_2\cap S_s, \dots, P_p \cap S_s),
 \]
  where empty sets are to be deleted. 

We consider the following subclasses of ordered set partitions on $[n]$.   
An ordered set partition $\pi\in\OSP(n)$ is called \textit{noncrossing}
    if the underlying set partition $\bar{\pi}$ is noncrossing.
   The set of the ordered noncrossing set partitions of $[n]$ is denoted $\ONCP(n)$. 
 For an ordered noncrossing partition $\pi=(P_1,P_2,\dots, P_p)$, we introduce the following subsets of $\Inn(\pi)$: 
\[
\AA(\pi) = \{P \in \Inn(\pi): \text{$\Cov(P)$ precedes $P$}  \} \quad \text{and} \quad \BB(\pi) = \{P \in \Inn(\pi): \text{$P$ precedes $\Cov(P)$}  \}. 
\]
Note that $\Inn(\pi) = \AA(\pi)\cup \BB(\pi).$  An example is shown in Fig.~\ref{fig:ordered_nc}, in which $\AA(\pi)=\{P_5, P_7 \}$ and $\BB(\pi) =\{P_1, P_2, P_6 \}$.

 A \emph{monotone partition} is a noncrossing ordered set partition such that every inner block is preceded by its covers.  
  The set of the monotone partitions is denoted $\MP(n)$. 
  
 An \emph{antimonotone partition} is a noncrossing ordered set partition such that every inner block precedes its covers.   
  The set of the antimonotone partitions is denoted $\AMP(n)$.

Let $I$ be a toset. The \emph{ordered kernel set partition} $\kappa(\bfi)$ of a sequence 
    $\bfi= (i_1,i_2,\ldots,i_n) \in \fin{I}$ is defined as follows.  
  First, pick the smallest value, say $p_1$,
  from $i_1, i_2,\ldots, i_n$  and define the block $P_1=\{k \in[n]
  \mid i_k=p_1\}$.  
  Then pick the second smallest value $p_2$ from $i_1, i_2,\ldots,
  i_n$ and define the block $P_2=\{k \in[n]\mid i_k=p_2 \}$.  Repeating this procedure leads us to an ordered set partition $(P_1, P_2, \ldots)$,
  which we denote by  $\kappa(\bfi)$.  Forgetting the order structure we set $\overline{\kappa} (\bfi) := \overline{\kappa(\bfi)}$.

%%%%%%%%%%%%%%%%%%%%%%%%%%%%%%%%%%%%%%%%%%%%%%%%%%%%%%%%%%%%%%%%%%
\subsection{General theory of cumulants for spreadability systems} \label{subsec:spreaded_cumulants}

The notion of spreadability system is a generalized notion of ``independent random variables''. We modify the definition a bit from \cite{HL17} because  it is natural in this paper to work on the category of \emph{unital} algebras.

\begin{defi} \label{def:sprsyst}
Let $(\mathcal{A}, \vp)$ be a unital algebraic probability space.   
 A \emph{unital spreadability system} for $(\mathcal{A}, \vp)$ is a triplet
  $\cS=(\cU, \widetilde{\vp}, (\iota^{(i)})_{i=1}^\infty)$ satisfying the
  following properties:

\begin{enumerate}[label=\rm(\roman*)]
 \item $(\cU, \widetilde{\vp})$ is a unital algebraic probability space, 

 \item $\iota^{(j)}\colon \cA \to \cU$ is a unital homomorphism
 such that $\vp = \widetilde{\vp} \circ \iota^{(j)}$ for
  each $j \geq 1$. For simplicity, $\iota^{(j)}(a)$ is denoted by $a^{(j)}$, $a \in
  \mathcal{A}$. 
  %, and we denote by $\mathcal{A}^{(j)}$ the image of $\mathcal{A}$ under $\iota^{(j)}$.
 \item The identity 
  \begin{equation}\label{eq:spreadability}
    \widetilde{\varphi}(a_1^{(i_1)}a_2^{(i_2)}\cdots a_n^{(i_n)}) = \widetilde{\varphi}(a_1^{(j_1)}a_2^{(j_2)}\cdots a_n^{(j_n)}) 
  \end{equation}
  holds for any $a_1, a_2,\ldots, a_n \in \mathcal{A}$ and any $i_1, i_2, \ldots, i_n, j_1,j_2,\dots, j_n \in
  \N$ such that $\kappa(i_1,i_2, \ldots, i_n) = \kappa(j_1,j_2,\ldots,j_n)$. 
\end{enumerate}
\end{defi}

\begin{defi}\label{def:dot}
 Let $\cS=(\mathcal{U}, (\iota^{(j)})_{j \geq 1}, \widetilde{\varphi})$ be a unital spreadability system for a unital algebraic probability space $(\mathcal{A}, \varphi)$. 
\begin{enumerate}[label=\rm(\roman*)]
\item  \label{eq:deltaAX}
 
The \emph{dot operation} for $a \in \cA$ is the sum of copies of $a$: 
$$
N.a:=a^{(1)} + a^{(2)} + \cdots + a^{(N)} \in \cU,   \quad N\in \N;  \qquad 0.a :=0.  
$$ 

\item For every $\pi \in \OSP(n)$ 
we define a multilinear functional
$\varphi_\pi\colon \cA^n \to \comp$ by choosing any representative sequence
$(i_1, i_2, \ldots, i_n)$ with $\kappa(i_1, i_2, \ldots, i_n)=\pi$ and setting
\begin{equation}
  \label{eq:phipi}
  \vp_{\pi}(a_1, a_2, \ldots,a_n)=\wt\vp(a_1^{(i_1)} a_2^{(i_2)} \cdots a_n^{(i_n)}). 
\end{equation}
Invariance \eqref{eq:spreadability} ensures that
this definition does not depend on the choice of the representative.

\item We define
 \begin{equation}
   \label{eq:phipixjij}
 \wt\vp_{\pi}(a_1^{(i_1)}, a_2^{(i_2)}, \ldots,a_n^{(i_n)}):=\vp_{\pi \curlywedge \kappa(i_1,i_2,\ldots,i_n)}(a_1,a_2,\ldots,a_n),    
 \end{equation}
and 
 $$
 \wt\vp_\pi(N.a_1,N.a_2,\ldots,N.a_n):=
\sum_{(i_1,i_2,\dots, i_n) \in [N]^n} 
\wt \vp_{\pi}(a_1^{(i_1)}, a_2^{(i_2)}, \ldots,a_n^{(i_n)}). 
 $$
\end{enumerate}
It holds that $\wt \varphi_{\pi}(N.a_1,N.a_2,\ldots,N.a_n)$ is a polynomial on $N$ whose coefficients of $N^k$ are zero for every $0\le k < |\pi|$. 
\end{defi}

\begin{rem}
Roughly, $\{a^{(i)}\}_{i \in \N}$ is meant to be ``i.i.d.\ copies of $a$'' and $ \vp_{\pi}(a_1, a_2, \ldots,a_n)$ is meant to be the value $\vp(a_1a_2 \cdots a_n)$ computed under the assumption that $(\{a_i\}_{i\in P_1}, \{a_i\}_{i\in P_2}, \dots, \{a_i\}_{i\in P_{|\pi|}})$ is ``independent''. 
\end{rem}

\begin{defi}\label{def:cumulants} In the setting of Definition \ref{def:dot}, for $\pi \in \OSP(n)$ the coefficient of $N^{|\pi|}$ in the polynomial $\wt \varphi_{\pi}(N.a_1,N.a_2,\ldots,N.a_n)$ is called the \emph{$\cS$-cumulant (over $\pi$)} and is denoted $K_\pi^\cS(a_1,a_2,\dots, a_n)$. In particular, $K_{1_n}^\cS$ is denoted $K_n^\cS$.  For a single element $a \in \cA$, the evaluation $K_\pi^\cS(a):= K_\pi^\cS(a,a,\dots,a)$ on the diagonal is called the \emph{$\cS$-cumulant of $a$} (over $\pi$). 
  \end{defi}
  Note that $K_\pi^\cS$ in this general definition does not guarantee any kind of multiplicativity.  Typically some $K_\pi^\cS$'s are multiplicative and some $K_\pi^\cS$'s are zero for specific kinds of partitions $\pi$. For example, for the free spreadability system $\cS$ (see Example \ref{exa:c-free}), $K_\pi^\cS=0$ for all crossing ordered set partitions $\pi$, while $K_\pi^\cS(a_1,\dots, a_n) = \prod_{P\in \pi} K_{|P|}^\cS(a_P)$ for all noncrossing ordered set partitions $\pi$. This fact together with the following general moment-cumulant formula entails the usual moment-cumulant formula for free cumulants. 

\begin{thm}[{\cite[Theorem 4.8]{HL17}}] \label{mc2} For any $\pi \in \OSP(n)$, we have 
  \[
  \varphi_{\pi}(a_1, a_2,\ldots, a_n) = \sum_{\substack{ \sigma \in \OSP(n) \\ \sigma \leq \pi } } \frac1{[\sigma:\pi]!}K_{\sigma}^\cS(a_1,a_2,\ldots,a_n) ,  
  \]
  where $[\sigma:\pi]! := \prod_{P \in \pi} (\#\{S \in \sigma: S \subseteq P\})!$. 
\end{thm}  

%{\color{red} Extensivity here and maybe uniqueness here} 

\begin{prop} 
The $\cS$-cumulants satisfy the following \emph{extensivity}: 
\[
K_\pi^\cS (N.a_1,\dots, N.a_n) = N^{|\pi|} K_\pi^\cS(a_1,\dots, a_n)
\] 
for all $a_1,\dots, a_n \in\cA, $ $\pi \in \OSP(n)$ and $N \in \N$. 
\end{prop}

\begin{rem}
A spreadability system $\cS$ is called an \emph{exchangeability system} if $\vp_\pi = \vp_{\pi'}$ whenever $\overline{\pi} = \overline{\pi'}$, i.e.\ $\vp_\pi$ is independent of orders on the set partition $\overline{\pi}$. In this case the cumulants $K_\pi^\cS$ are also independent of orders on $\pi$, so that it suffices to consider set partitions $\pi$. If $\cS$ is an exchangeability system then the cumulants satisfy the vanishing property that is stronger than the extensivity \cite[Proposition 7.1]{HL17}: if $\pi$ contains a block $P$ which can be partitioned into nonempty two subsets, $P=P_1 \sqcup P_2$, such that the subsets $\{a_i: i \in P_1\}$ and $\{a_i: i \in P_2\}$ are ``$\cS$-independent''  then $K_\pi^\cS(a_1,\dots, a_n)=0$.  This vanishing property further implies the additivity of cumulants of single elements: if $a_1$ and $a_2$ are $\cS$-independent then 
\[
K_n^\cS (a_1+a_2) = K_n^\cS (a_1) + K_n^\cS (a_2). 
\]
The exchangeability holds for the standard spreadability systems for free, Boolean and c-free independences, see \cite[Subsection 3.3]{HL17}.  On the other hand, the exchangeability fails to hold for monotone, antimonotone and c-monotone spreadability systems. 
\end{rem}

\begin{exa}[C-free cumulants] \label{exa:c-free} C-free cumulants, originally introduced in \cite{BLS}, can be formulated in terms of exchangeability systems \cite{Leh}.
Let $(\cA,\vp,\psi,\theta)$ be a unital 2-algebraic probability space. We take copies $\{(\cA_n,\vp_n,\psi_n)\}_{n\in\N}$ of $(\cA,\vp,\psi)$ and set the c-free product 
\[
(\cU^\Free, \wt\vp^\CF, \wt\psi^\Free) := \freeproe{n\in \N} (\cA_n,\vp_n,\psi_n)
\] 
together with the natural embeddings 
$
\iota^{(j)}\colon \cA \xhookrightarrow{}  \cU^\Free 
$
into the $j$-th component of $\cU^\Free$. 
Then both 
\[
\cS^{\CF}:= (\cU^\Free, (\iota^{(j)})_{j\in\N}, \wt\vp^\CF) \quad \text{and} \quad \cS^{\Free}:= (\cU^\Free, (\iota^{(j)})_{j\in\N}, \wt\psi^\Free)
\]
are unital spreadability systems for the unital algebraic probability spaces $(\cA, \vp)$ and $(\cA,\psi)$, respectively. Because of symmetry of c-free product, it turns out that they are exchangeability systems. Hence, we can associate the c-free cumulants  $(K_\pi^{{\CF}(\vp, \psi)})_{\pi \in \SP}$ and the free cumulants $(K_\pi^{{\Free}(\psi)})_{\pi \in \SP}$.  Note that the relation
\begin{equation} %\label{eq:special_cases}
K_\pi^{{\Free}(\psi)}  = K_\pi^{{\CF}(\psi, \psi)}
\end{equation}
holds because $\wt\vp^\CF = \wt\psi^\Free$ if $\vp = \psi$. It is demonstrated in \cite[Subsection 4.7]{Leh}  that 
\[
K^{\CF(\vp,\psi)}_\pi(a_1,a_2,\dots,a_n) = \begin{cases} 
\displaystyle\left[\prod_{P \in \Out(\pi)} K^{\CF(\vp,\psi)}_{|P|}(a_P) \right] \left[\prod_{P \in \Inn(\pi)} K^{\Free(\psi)}_{|P|}(a_P) \right], & \text{if~} \pi \in \NC(n), \\[8mm] 
0,& \text{if~} \pi \notin \NC(n). 
\end{cases}
\]
This multiplicativity and Theorem \ref{mc2} recover the moment-cumulant formula in \cite{BLS}: 
 \begin{equation}\label{eq:cfree_mc}
  \varphi(a_1 a_2\cdots a_n) = \sum_{\substack{ \pi \in \NC(n)} } \left[\prod_{P \in \Out(\pi)} K^{\CF(\vp,\psi)}_{|P|}(a_P) \right] \left[\prod_{P \in \Inn(\pi)} K^{\Free(\psi)}_{|P|}(a_P) \right]  
  \end{equation}
  and, in the special case $\vp=\psi$, the formula in  \cite{Spe2} 
 \begin{equation}\label{eq:free_mc}
  \psi(a_1 a_2\cdots a_n) = \sum_{\substack{ \pi \in \NC(n)} }  \left[\prod_{P \in \pi} K^{\Free(\psi)}_{|P|}(a_P) \right].   
  \end{equation}
 Formulas  \eqref{eq:cfree_mc} and \eqref{eq:free_mc} provide a recursive way to express cumulants $K_n^{\CF(\vp,\psi)}, K_n^{\Free(\psi)}$ in terms of $\vp,\psi$. For $n=1,2$ the cumulants are given by 
   \begin{align}
 K^{\CF(\vp,\psi)}_{n}(a) &=   \varphi(a), & K^{\CF(\vp,\psi)}_{n}(a_1,a_2) &= \vp(a_1 a_2) -\vp(a_1)\vp(a_2), \\
  K^{\Free(\psi)}_{n}(a) &=   \psi(a), &    K^{\Free(\psi)}_{n}(a_1,a_2) &=   \psi(a_1a_2) - \psi(a_1) \psi(a_2).   
 \end{align}
 In general, cumulants of order $n$ are of the forms
 \begin{align}
 K^{\CF(\vp,\psi)}_{n}(a_1,\dots, a_n) &=   \varphi(a_1 a_2\cdots a_n)  +  (\text{polynomial on $\{\vp(a_P), \psi(a_P): P \subsetneq [n]\}$}), \label{eq:cfree_cm}  \\
  K^{\Free(\psi)}_{n}(a_1,\dots, a_n) &=   \psi(a_1 a_2\cdots a_n)  +  (\text{polynomial on $\{\psi(a_P): P \subsetneq [n]\}$}).  \label{eq:free_cm}
 \end{align}
 Note that the polynomials above are actually multilinear functionals on $\cA^n$. For example, by induction on $n$ the polynomial part of \eqref{eq:cfree_cm} is of the form
 \[
 \sum_{\substack{(\pi, f)\in \CSP^2(n) \\ \pi \in \NC(n), \pi < 1_n}}  \delta(\pi,f) \left[ \prod_{\substack{P \in \pi \\ f(P)=1}} \vp(a_P) \right]  \left[ \prod_{\substack{ P \in \pi \\  f(P)=2}} \psi(a_P) \right], 
 \]
 where $\delta(\pi,f) \in \real$ is a universal coefficient. 
    
Finally, because of exchangeability, c-free and free cumulants have the vanishing property: given $a_1,a_2,\dots,a_n \in \cA$ for which there is a subset $P \subseteq [n]$ such that $P$ and $[n]\setminus P$ are nonempty and $\{a_i: i \in P\}$ and $\{a_i: i \in [n] \setminus P\}$ are c-free then we have 
\begin{equation}\label{eq:cfree_vanish}
K^{\CF(\vp,\psi)}_n(a_1,a_2,\dots,a_n) = K^{\Free(\psi)}_n(a_1,a_2,\dots,a_n) =0. 
\end{equation}
Note that it suffices to consider $K_n^{\CF(\vp,\psi)}$ and $K_n^{\Free(\psi)}$ because the other cumulants over $\pi$ are just the product of them or zero. 
\end{exa}

%%%%%%%%%%%%%%%%%%%%%%%%%%%%%%%%%%%%%%%%%%%%%%%%%%%%%%%%%%%%%%%%%%%%
\subsection{The highest coefficients for $\AP$-freeness} \label{subsec:highest}

The highest coefficients for each notion of independence are crucial objects to prove a moment-cumulant formula. This concept uncovers the way each notion of independence associates a specific class of (ordered or unordered) set partitions. Here we formulate the highest coefficients for $\AP$-freeness and determine them.

\begin{prop}  \begin{enumerate}[label=\rm(\arabic*)] 

\item There exists an  integer $c^{\AP}((\pi,f), \sigma)$ for every $(\pi,f) \in \CSP^3(n), \sigma \in \OSP(n), \pi \le \overline{\sigma},  n\in\N$, called a \emph{universal coefficient for $\AP$-freeness}, such that for every family of unital $3$-algebraic probability spaces $\{(\cA_i, \vp_i, \psi_i,\theta_i)\}_{i\in I}$ with $I$ a toset, every $\bfi =(i_1,i_2,\dots, i_n)\in \fin{I}$ and $(a_1,a_2,\dots, a_n) \in \cA_\bfi$ we have 
\begin{equation}\label{eq:universal_coefficients}
\vp (a_1 a_2 \cdots a_n) = \sum_{\substack{(\pi,f) \in \CSP^3(n) \\ \pi \le \overline{\kappa}(\bfi)  }}  c^{\AP}((\pi,f), \kappa(\bfi)) \left[ \prod_{\substack{P \in \pi \\ f(P)=1}} \vp(a_P) \right]  \left[ \prod_{\substack{ P \in \pi \\  f(P)=2}} \psi(a_P) \right]     \left[ \prod_{\substack{ P \in \pi \\  f(P)=3}} \theta(a_P) \right],  
\end{equation}
where $(\cA, \vp, \psi,\theta) = \fmp_{i\in I}(\cA_i, \vp_i, \psi_i,\theta_i)$ is the $\AP$-free product. 

\item There exist  integers $c^\BP((\pi,f), \sigma), c^\GP((\pi,f), \sigma)$ for every $(\pi,f) \in \CSP^2(n)$, $\sigma \in \OSP(n)$, $\pi \le \overline{\sigma},  n\in\N$, called the \emph{universal coefficients for $\BGP$-freeness}, such that for every family of unital $2$-algebraic probability spaces $\{(\cA_i, \psi_i,\theta_i)\}_{i\in I}$ with $I$ a toset, every $\bfi =(i_1,i_2,\dots, i_n)\in \fin{I}$ and $(a_1,a_2,\dots, a_n) \in \cA_\bfi$ we have 
\begin{align}
\psi (a_1 a_2 \cdots a_n) 
&= \sum_{\substack{(\pi,f) \in \CSP^2(n) \\ \pi \le \overline{\kappa}(\bfi)  }}  c^\BP((\pi,f), \kappa(\bfi))   \left[ \prod_{\substack{ P \in \pi \\  f(P)=1}} \psi(a_P) \right]     \left[ \prod_{\substack{ P \in \pi \\  f(P)=2}} \theta(a_P) \right], \label{eq:universal_coefficients_AF}    \\
\theta (a_1 a_2 \cdots a_n) 
&= \sum_{\substack{(\pi,f) \in \CSP^2(n) \\ \pi \le \overline{\kappa}(\bfi)  }}  c^\GP((\pi,f), \kappa(\bfi)) \left[ \prod_{\substack{P \in \pi \\ f(P)=1}} \psi (a_P) \right]      \left[ \prod_{\substack{ P \in \pi \\  f(P)=2}} \theta(a_P) \right], \label{eq:universal_coefficients_OAF}
\end{align}
where $(\cA, \psi,\theta) = \fmp_{i\in I}(\cA_i, \psi_i,\theta_i)$ is the $\BGP$-free product. 
\end{enumerate}
\end{prop}

\begin{proof} The existence of universal coefficients for $\AP$-freeness can be verified by induction on $n$. Then we apply the  induction hypothesis to the recursive formula 
\begin{align}
\vp(a_1a_2 \cdots a_n) &=  \sum_{S \subseteq [n], S\neq \emptyset }\left[\prod_{s \in S} \alpha_{\bfi,s}(a_{s}) \right] \vp(\mathring{a}_{S_0} \mathring{a}_{S_1} \cdots   \mathring{a}_{S_{k}}),
\end{align}
where the notation is taken from \eqref{eq:decomp0}, to arrive at the desired conclusion. 

Because the moment formulas for $\psi$ and $\theta$ are special cases of $\AP$-freeness with $\vp =\psi$ and $\vp =\theta$, respectively, \eqref{eq:universal_coefficients_AF} and \eqref{eq:universal_coefficients_OAF} follow from \eqref{eq:universal_coefficients}. 
\end{proof}

Recall from Proposition \ref{prop:useful} that selecting the special combination $(\vp,\psi,\psi)$ makes $\AP$-free product the c-free product. Therefore we are able to obtain the following more or less known result (note that using the moment-cumulant formula yields a more concrete coefficients).

\begin{cor}\label{cor:universal_coefficients_c-free}  There exists an integer $c^{\CF}((\pi,f), \sigma)$ for every $(\pi,f) \in \CSP^2(n), \sigma \in \SP(n), \pi \le \sigma,  n\in\N$ such that for every family of unital $2$-algebraic probability spaces $\{(\cA_i, \vp_i, \psi_i)\}_{i\in I}$ with $I$ a set, every $\bfi =(i_1,i_2,\dots, i_n)\in \fin{I}$ and $(a_1,a_2,\dots, a_n) \in \cA_\bfi$ we have 
\begin{equation}\label{eq:universal_coefficients_c-free}
\vp (a_1 a_2 \cdots a_n) = \sum_{\substack{(\pi,f) \in \CSP^2(n) \\ \pi \le \overline{\kappa}(\bfi)  }}  c^{\CF}((\pi,f),  \overline{\kappa}(\bfi)) \left[ \prod_{\substack{P \in \pi \\ f(P)=1}} \vp(a_P) \right]  \left[ \prod_{\substack{ P \in \pi \\  f(P)=2}} \psi(a_P) \right], 
\end{equation}
where $(\cA, \vp, \psi) = \freeproe{i\in I}(\cA_i, \vp_i, \psi_i)$ is the c-free product. 
\end{cor}

\begin{rem} 
The fact that c-freeness is symmetric implies that $c^{\CF}((\pi,f), \sigma)$ does not depend on the order on $\sigma$. Therefore, we formulated Corollary \ref{cor:universal_coefficients_c-free} in terms of unordered set partitions. 
\end{rem}

The goal is to determine the numbers $c^{\AP}((\pi,f),\sigma)$ when $\overline{\sigma} = \pi$. For this we prepare some lemmas.

\begin{lem}\label{lem:c-free_highest}
There exists an integer $d((\pi,f), (\sigma,g))$ for every $(\pi,f) \in \CSP^2(n+1), (\sigma,g)\in \CSP^2(n), n \in \N$ such that for any unital 2-algebraic probability spaces $\{(\cA_i, \vp_i,\psi_i)\}_{i \in \{1,2\}}$ and every $n\in \N, x_1,x_2,\dots, x_{n+1} \in \cA_1$ and $y_1,y_2,\dots, y_n \in \cA_2$ we have   
\begin{align}%\label{eq:universal_coefficients_c-free_special}
&\vp (x_1 y_1  x_2 y_2 \cdots y_{n} x_{n+1})  \label{eq:goal1} \\
&= [ \vp(x_1 x_{n+1}) - \vp(x_1)\vp(x_{n+1})] \left[\prod_{i =2}^{n}  \psi(x_i)\right] \psi(y_1 \cdots y_{n})  \label{eq:goal2}\\ 
&\quad + \vp(x_1)\vp(x_{n+1}) \left[\prod_{i =2}^{n}  \psi(x_i)\right]  \vp(y_1 \cdots y_{n}) \label{eq:goal3} \\
&\quad + \sum_{(\pi, f), (\sigma, f)}  d((\pi,f), (\sigma, g)) \left[ \prod_{\substack{P \in \pi \\ f(P)=1}} \vp(x_P) \right]  \left[ \prod_{\substack{ P \in \pi \\  f(P)=2}} \psi(x_P) \right] \left[ \prod_{\substack{S \in \sigma \\ g(P)=1}} \vp(y_S) \right]  \left[ \prod_{\substack{S \in \sigma \\  g(P)=2}} \psi(y_S) \right],  \label{eq:goal4}
\end{align}
where  $(\cA, \vp, \psi) = (\cA_1, \vp_1,\psi_1) \freeprod (\cA_2, \vp_2,\psi_2)$  is the c-free product and $(\pi,f), (\sigma,g)$ run over  $\CSP^2(n+1)$ and $\CSP^2(n)$, respectively, so that $\sigma \ne 1_n$. 
\end{lem}

\begin{proof} 
Because of the vanishing of c-free cumulants \eqref{eq:cfree_vanish}, only noncrossing set partitions $\pi$ of $[2n+1]$ contribute to the moment-cumulant formula \eqref{eq:cfree_mc}, where $a_1=x_1, a_2= y_1, a_3=  x_2, \dots, a_{2n+1} = x_{n+1}$. As a consequence, the factor $\psi(y_1 \cdots y_{n})$ appears in \eqref{eq:cfree_mc} if and only if $\pi=\{\{1,2n+1\}, \{2,4,6,\dots, 2n\},$ $\{3\}, \{5\}, \dots, \{2n-1\}\}$. For this $\pi$ we have
\begin{align*}
&\left[\prod_{P \in \Out(\pi)} K^{\CF(\vp,\psi)}_{|P|}(a_P) \right] \left[\prod_{P \in \Inn(\pi)} K^{\Free(\psi)}_{|P|}(a_P) \right]    \\
&\qquad = K_{2}^{\CF(\vp,\psi)}(x_1,x_{n+1}) \left[\prod_{i =2}^{n}  K_1^{\Free(\psi)}(x_i)\right] K_{n}^{\Free(\psi)}(y_1,y_2,\dots, y_n),  
\end{align*}
which consists of the term \eqref{eq:goal2} and terms contained in \eqref{eq:goal4}. In a similar manner \eqref{eq:goal3} appears only from $\pi=\{\{1\}, \{2n+1\}, \{2,4,6,\dots, 2n\}, \{3\}, \{5\}, \dots, \{2n-1\}\}$. The other noncrossing partitions $\pi$ only yield terms to be included in \eqref{eq:goal4}. 
\end{proof}

\begin{lem}\label{lem:c-free_highest2}
In the setting of Corollary \ref{cor:universal_coefficients_c-free}, let $I=\{1,2\}$. Suppose that there exist $p < s < q < t$ such that $i_p = i_q =1$  and $i_s = i_t =2$. Suppose also that a set partition $\pi \le \overline{\kappa}(\bfi)$ contains two blocks $P_1, P_2$ such that $p, q \in P_1$ and $r,s \in P_2$. Then $c^{\CF}((\pi,f),  \overline{\kappa}(\bfi)) =0$ for every coloring $f$.  
\end{lem}
\begin{proof}
In the moment-cumulant formula \eqref{eq:cfree_mc}, a nonzero contribution appears in the sum only if the noncrossing set partition $\pi$ satisfies $\pi \le \overline{\kappa}(\bfi)$ because of vanishing property of c-free and free cumulants. Such a set partition $\pi$ never contains two blocks $P_1$ and $P_2$ such that $p, q \in P_1$ and $r,s \in P_2$. Substituting the cumulant-moment formulas \eqref{eq:free_cm} and \eqref{eq:cfree_cm} into the moment-cumulant formula, we get the desired conclusion. 
\end{proof}

\begin{thm}\label{thm:highest_coefficients}
The integers $c^{\AP}(\sigma; g):= c^{\AP}((\overline{\sigma},g),\sigma),  \sigma \in \OSP(n), g\colon \overline{\sigma} \to [3], n\in \N$, called the \emph{highest coefficients for $\AP$-freeness}, are determined as follows: 
\begin{equation}\label{eq:highest_coefficients}
c^{\AP}(\sigma; g) = 
\begin{cases}
1,  & \text{if $\sigma \in \ONCP(n)$,  $g(S) = 1$ for $S \in \Out(\sigma)$,  $g(S) = 2$ for $S \in \AA(\sigma)$} \\ 
& \text{and $g(S)=3$ for $S \in \BB(\sigma)$}, \\
0 , & \text{otherwise}. 
\end{cases}
\end{equation}
Similarly, the highest coefficients for $\BGP$-freeness $c^\BP(\sigma; g):= c^\BP((\overline{\sigma},g),\sigma)$ and $c^\GP(\sigma; g):= c^\GP((\overline{\sigma},g),\sigma)$,  $\sigma \in \OSP(n),  g\colon \overline{\sigma} \to [2], n\in \N$ are given by  
\begin{align}
c^\BP(\sigma; g) 
&= 
\begin{cases}
1,  & \text{if $\sigma \in \ONCP(n)$,  $g(S) = 1$ for $S \in \Out(\sigma) \cup \AA(\sigma)$} \\ 
&\text{and  $g(S) = 2$ for $S \in \BB(\sigma)$}, \\ 
0 , & \text{otherwise},
\end{cases} 
\label{eq:highest_coefficients_AF}   \\
c^\GP(\sigma; g) &= 
\begin{cases}
1,  & \text{if $\sigma \in \ONCP(n)$,  $g(S) = 1$ for $S \in \Out(\sigma) \cup \BB(\sigma)$} \\ 
&\text{and  $g(S) = 2$ for $S \in \AA(\sigma)$}, \\ 
0 , & \text{otherwise}. 
\end{cases}
\label{eq:highest_coefficients_OAF}
\end{align}
\end{thm}

\begin{proof}

The proof is given by induction on $n$. We take a family  of unital 3-algebraic probability spaces $\{(\cA_i, \vp_i,\psi_i,\theta_i)\}_{i\in I}$ with $I$ a toset, set $(\cA,\vp,\psi, \theta)$ to be their $\AP$-free product,  take $\bfi =(i_1,\dots, i_n) \in \fin{I}$ of length $n$ and $(a_1,a_2,\dots, a_n) \in \cA_\bfi$ and set $\sigma:= \kappa(\bfi) = (S_1,S_2,\dots, S_\ell)$.

\begin{figure}[b]
\begin{center}
\begin{tikzpicture}[scale=0.7]
\draw[-] (3,0.5) -- (3,1.5) -- (4,1.5) --(4,0.5);
\draw[-] (5,0.5) -- (5,1.5);
\draw[-] (8,0.5) -- (8,1.5) -- (9,1.5) --(9,0.5);
\draw[line width=1.5pt] (2,0.5) -- (2,2.5) -- (6,2.5) --(6,0.5);
\draw[line width=1.5pt] (6,2.5) -- (11,2.5) --(11,0.5);
\draw[-] (7,0.5) -- (7,2) -- (10,2) --(10,0.5);
%\draw[-] (1,0.5) -- (1,3) -- (12,3) --(12,0.5);
\draw[-] (12,0.5) -- (12,2.5) -- (13,2.5) --(13,0.5);
%\node at (1,0) {1};
\node at (2,0) {1};
\node at (3,0) {2};
\node at (4,0) {3};
\node at (5,0) {4};
\node at (6,0)  {5}; 
\node at (7,0)  {6}; 
\node at (8,0) {7};
\node at (9,0) {8};
\node at (10,0) {9};
\node at (11,0) {10};
\node at (12,0)  {11}; 
\node at (13,0)  {12}; 
%\node at (14,0) {14};
%%%%%
\node at (8.5,2.2) {\tiny1};
\node at (3.5,1.7) {\tiny 2};
\node at (12.5,2.7)  {\tiny 3}; 
%\node at (6.5,3.2)  {\tiny4}; 
\node at (8.5,1.7) {\tiny 4};
\node at (5,1.7)  {\tiny 5}; 
\node at (6.5,2.75) {\tiny 6};
\end{tikzpicture}
\end{center}
\caption{The diagram of $\sigma = (S_1, \cdots, S_6) \in \ONCP(12)$ for which $S_6 = \{1, 5, 10 \}$. Then $B_1 = \emptyset, B_2= \{2,3,4\}, B_3=\{6,7,8,9\}, B_4=\{11,12\}$. }
\label{fig:interval}

\end{figure}

\vspace{2mm}
\noindent
{\bf Case a.}  First we treat the case $\sigma \in \ONCP(n)$. 
Let $j = \max\{i_1,\dots, i_n\}$, i.e.\ $S_\ell = \{i_\ell: \ell \in [n], i_\ell = j \}$.  The set $S_\ell$ divides $[n]$ into (possibly empty) intervals $B_1,B_2,\dots, B_{\ell+1}$, such that $[n]= S_\ell \sqcup B_1 \sqcup \cdots \sqcup B_{\ell+1}$ and $a \in B_i, b \in B_j, i<j$ implies $a<b$, see Figure \ref{fig:interval}. Let $\cA_{<j}$ and $\cA_{\ge j}$ be the subalgebra generated by $\{\cA_i\}_{i \in I, i< j}$ and $\{\cA_i\}_{i \in I, i\ge  j}$, respectively. By the associativity of $\AP$-free product and Proposition \ref{rem:cfree_alpha}, we have 
\[
\vp =  \cf{(\Restr{\vp}{\cA_{<j}})} {(\Restr{\theta}{\cA_{<j}})} {(\Restr{\psi}{\cA_{\ge j}} )} {(\Restr{\vp}{\cA_{\ge j}})}, 
\]
\[ \psi =  \cf{(\Restr{\psi}{\cA_{<j}})} {(\Restr{\theta}{\cA_{<j}})} {(\Restr{\psi}{\cA_{\ge j}})} {(\Restr{\psi}{\cA_{\ge j}})} \quad  \text{and} \quad \theta =  \cf{(\Restr{\theta}{\cA_{<j}})} {(\Restr{\theta}{\cA_{<j}})} {(\Restr{\psi}{\cA_{\ge j}})} {(\Restr{\theta}{\cA_{\ge j}})}
\]
under the natural identification $\cA \simeq \cA_{<j} \freeprod \cA_{\ge j}$. 

Combining consecutive elements $a_i$'s that belong to $\cA_{<j}$, we get a reduced form 
\[
a_1a_2 \cdots a_n = x_1 y_1 x_2 y_2 \cdots y_\ell x_{\ell+1},
\]
 where $x_1,\dots, x_{\ell+1} \in \cA_{<j}$ and $y_1,\dots, y_\ell \in \cA_j$. We take a convention here that $x_1=1$ if $B_1 =\emptyset$, and similarly for $x_{\ell+1}$. Lemma \ref{lem:c-free_highest} then yields
\begin{align}%\label{eq:universal_coefficients_c-free_special}
&\vp (x_1 y_1  x_2 y_2 \cdots y_{\ell} x_{\ell+1})  \label{eq:rgoal1} \\
&= [ \vp(x_1 x_{\ell+1}) - \vp(x_1)\vp(x_{\ell+1})] \left[\prod_{i =2}^{\ell}  \theta(x_i)\right] \psi(y_1 \cdots y_{\ell})  \label{eq:rgoal2}\\ 
&\quad + \vp(x_1)\vp(x_{\ell+1}) \left[\prod_{i =2}^{\ell}  \theta(x_i)\right]  \vp(y_1 \cdots y_{\ell}) \label{eq:rgoal3} \\
&\quad + \text{[terms not containing $ \vp(y_1 \cdots y_{\ell}),  \psi(y_1 \cdots y_{\ell})$ or $ \theta(y_1 \cdots y_{\ell})$ as factors].} \notag
\end{align}

\noindent
{\bf Case a1:} $S_\ell \in \Out(\sigma)$. Because then $x_1$ and $x_{\ell+1}$ are elements of different algebras, we have $\vp(x_1 x_{\ell+1}) - \vp(x_1)\vp(x_{\ell+1})=0$. Therefore, the highest coefficient can appear only from part \eqref{eq:rgoal3} and hence for any $g$
\[
c^{\AP}(\sigma; g) =  
\begin{cases}
0, & \text{if}~ g (S_\ell) \in \{2,3\}, \\
c^{\AP}(\Restr{(\sigma;g)}{B_1})c^{\AP}(\Restr{(\sigma;g)}{B_{\ell+1}})\prod_{i=2}^\ell c^{\GP}(\Restr{(\sigma;g)}{B_i}), & \text{if}~ g(S_\ell)=1.  
\end{cases}
\]
 By the induction hypothesis applied to $c^{\AP}(\Restr{(\sigma;g)}{B_1}), c^{\AP}(\Restr{(\sigma;g)}{B_{\ell+1}}), c^{\GP}(\Restr{(\sigma;g)}{B_i})$, we can see that \eqref{eq:highest_coefficients} is satisfied. Setting the specialization $\vp_i = \psi_i$ and $\vp_i =\theta_i$ implies \eqref{eq:highest_coefficients_AF} and \eqref{eq:highest_coefficients_OAF}, respectively. 

\vspace{2mm}
\noindent
{\bf Case a2:} $S_\ell \in \AA(\sigma)$. Because $S_\ell$ is covered by another block, an alphabet in the word $x_1$ and an alphabet in the word $x_{\ell+1}$ belong to the same subalgebra, so that the highest coefficient never appears from part \eqref{eq:rgoal3} and it should arise from part \eqref{eq:rgoal2}, more precisely, from part 
\[
\vp(x_1 x_{\ell+1})  \left[\prod_{i =2}^{n}  \theta(x_i)\right] \psi(y_1 \cdots y_{n}). 
\]
This yields the recursive formula
\[
c^{\AP}(\sigma; g) =  
\begin{cases}
0, & \text{if}~ g (S_\ell) \in \{1,3\}, \\
 c^{\AP}(\Restr{(\sigma; g)}{B_1 \cup B_{\ell+1}})\prod_{i=2}^n c^{\GP}(\Restr{(\sigma;g)}{B_i}), & \text{if}~ g(S_\ell)=2.  
\end{cases}
\]
By the induction hypothesis applied to  $c^{\AP}(\Restr{(\sigma; g)}{B_1 \cup B_{\ell+1}})$ and $c^{\GP}(\Restr{(\sigma;g)}{B_i})$, we can see that \eqref{eq:highest_coefficients} is satisfied. Setting the specialization $\vp_i = \psi_i$ and $\vp_i =\theta_i$ implies \eqref{eq:highest_coefficients_AF} and \eqref{eq:highest_coefficients_OAF}, respectively.

\vspace{2mm}
\noindent
{\bf Case a3:} $S_\ell \in \BB(\sigma)$. This never occurs because $S_\ell$ precedes no other blocks of $\sigma$. 

\vspace{2mm}

The above arguments justify formulas \eqref{eq:highest_coefficients} --  \eqref{eq:highest_coefficients_OAF}  when $\sigma \in\ONCP(n)$.  

\vspace{2mm}
\noindent
{\bf Case b.} It remains to treat the case $\sigma \notin\ONCP(n)$. Then there are two blocks of $\sigma$, say $S_i$ and $S_j$ with $i<j$, and $p,q \in S_i$ and $r,s \in S_j$ such that $p< r< q< s$ or $p > r > q > s$. We set $\cA_{\le i}:= \freeproe{k \in I, k \le i}\cA_i$ and $\cA_{> i}:= \freeproe{k\in I, k >i} \cA_k$. Then we have 
\[
\vp = \cf{(\Restr{\vp}{\cA_{\le i}})} {\Restr{\theta}{\cA_{\le i}}} {\Restr{\psi}{\cA_{> i}}} {(\Restr{\vp}{\cA_{> i}})}.
\] 
By Lemma \ref{lem:c-free_highest2},  it never happens in formula \eqref{eq:universal_coefficients} that $a_p, a_q$ appear inside a common  linear functional and also $a_r, a_s$ appear inside a common linear functional. This implies that $c^{\AP}(\sigma; g) =0 $ for any $g$. 
\end{proof}

%We reduce $a_1a_2 \cdots a_n$ to $x_1 y_1 x_2 y_2 \cdots y_r x_{r+1}$, where $x_1, \dots, x_{r+1} \in \cA_{\le i}$  and $y_1,\dots, y_r \in \cA_{> i}$. 

%%%%%%%%%%%%%%%%%%%%%%%%%%%%%%%%%%%%%%%%%%%%%%%%%%%%%%%%%%%%%%%%%%%
\subsection{The $\AP$-free spreadability system and the $\AP$-free cumulants} \label{subsec:alpha_spread}

Let $(\cA,\vp,\psi,\theta)$ be a unital 3-algebraic probability space. We take countable copies $(\cA_n,\vp_n,\psi_n,\theta_n):=(\cA,\vp,\psi,\theta)$ for $n\in\N$ and set 
\[
(\cU^\Free, \wt\vp^\AP, \wt\psi^\BP, \wt\theta^\GP) := \fmp_{n\in \N} (\cA_n,\vp_n,\psi_n,\theta_n)
\] 
together with the natural embeddings 
$
\iota^{(j)}\colon \cA_j \xhookrightarrow{}  \cU^\Free. 
$
Then 
\[
\cS^{\AP}:= (\cU^\Free, (\iota^{(j)})_{j\in\N}, \wt\vp^\AP), \quad \cS^{\BP}:= (\cU^\Free, (\iota^{(j)})_{j\in\N}, \wt\psi^\BP) \quad \text{and} \quad \cS^{\GP}:= (\cU^\Free, (\iota^{(j)})_{j\in\N}, \wt\theta^\GP)
\]
are unital spreadability systems for the unital algebraic probability spaces $(\cA, \vp), (\cA,\psi)$ and $(\cA,\theta)$, respectively. Hence, along the lines of Subsection \ref{subsec:spreaded_cumulants} we can define the \emph{$\AP$-free cumulants} $(K_\pi^{{\AP}})_{\pi \in \OSP}$, the \emph{$\BP$-free cumulants} $(K_\pi^{{\BP}})_{\pi \in \OSP}$ and the \emph{$\GP$-free cumulants} $(K_\pi^{{\GP}})_{\pi \in \OSP}$.   When we emphasize the dependence on linear functionals, we will also use the notation $K_\pi^{{\AP}(\vp,\psi,\theta)}, K_\pi^{{\BP}(\psi,\theta)}, K_\pi^{{\GP}(\psi,\theta)}$. 

Note that, by Proposition \ref{prop:useful}, the relations 
\begin{equation} \label{eq:special_cases}
K_\pi^{{\BP}(\psi, \theta)}  = K_\pi^{{\AP}(\psi, \psi, \theta)}  \quad \text{and} \quad   K_\pi^{{\GP}(\psi, \theta)}  = K_\pi^{{\AP}(\theta, \psi, \theta)}
\end{equation}
hold, which are useful for shortening proofs below.

\begin{thm}\label{thm:m-c-alpha}
Let $(\cA, \vp, \psi, \theta)$ be a unital 3-algebraic probability space. Then 
\begin{align}
\vp(a_1 \cdots a_n)  
&= \sum_{\pi \in \ONCP (n)} \frac{1}{|\pi|!} \left[ \prod_{P \in \Out(\pi)   } K^{{\AP}}_{|P|}(a_P)\right] \left[ \prod_{P \in \AA(\pi)} K^{\BP}_{|P|}(a_P) \right]\left[  \prod_{P \in \BB(\pi)} K^{\GP}_{|P|}(a_P)\right],  \label{eq503}  \\ 
\psi(a_1 \cdots a_n) 
&= \sum_{\pi \in \ONCP (n)} \frac{1}{|\pi|!}   \left[ \prod_{P \in \Out(\pi)\cup \AA(\pi)} K^{\BP}_{|P|}(a_P) \right]\left[  \prod_{P \in \BB(\pi)} K^{\GP}_{|P|}(a_P)\right],  \label{eq501} \\
 \theta(a_1 \cdots a_n) 
 &= \sum_{\pi \in \ONCP (n)} \frac{1}{|\pi|!}  \left[ \prod_{P \in \AA(\pi)} K^{\BP}_{|P|}(a_P) \right]\left[  \prod_{P \in \Out(\pi) \cup \BB(\pi)} K^{\GP}_{|P|}(a_P)\right].     \label{eq502}
\end{align}
\end{thm}

\begin{proof} It suffices to prove \eqref{eq503} because formulas \eqref{eq501} and \eqref{eq502} follow from \eqref{eq503}  in view of \eqref{eq:special_cases}. 
For $\pi \in \OSP(n)$ and $u_1, \dots, u_n \in \cU^\Free$ let 
\[
\wt\vp^\AP_{(\pi)} (u_1,\dots, u_n) :=  \left[ \prod_{P \in \Out(\pi)} \wt\vp^\AP(u_P) \right]  \left[ \prod_{P \in \AA(\pi)} \wt\psi^\BP(u_P) \right]     \left[ \prod_{P \in \BB(\pi)} \wt\theta^\GP(u_P) \right]  
\]
and 
\[
\wt\vp^\AP_{<\pi}(u_1, \dots, u_n) :=  \sum_{\substack{(\rho,h) \in \CSP^3(n) \\ \rho < \overline{\pi}  }}  c^{\AP}((\rho,h), \pi ) \left[ \prod_{\substack{R \in \rho \\ h(R)=1}} \vp(u_R) \right]  \left[ \prod_{\substack{ R \in \rho \\  h(R)=2}} \psi(u_R) \right]     \left[ \prod_{\substack{ R \in \rho \\  h(R)=3}} \theta(u_R) \right]. 
\]
Below we  regard $\cA$ as a subalgebra of $\cU^\Free$ via the embedding $\iota^{(1)}$ and make use of the simplified notation $\wt\vp^\AP_\pi(a_1,\dots, a_n)$ and $\wt\vp^\AP_{<\pi}(a_1, \dots, a_n)$ for $a_1,a_2,\dots, a_n \in \cA$. 
As a consequence of Theorem \ref{thm:highest_coefficients}, formula \eqref{eq:universal_coefficients} can be expressed in the form  
\begin{align}\label{eq:universal_coefficients2}
 \vp_\pi (a_1, \cdots, a_n) 
&= \begin{cases}
\wt\vp^\AP_{(\pi)} (a_1,\dots, a_n)    +   \wt\vp^\AP_{< \pi}(a_1, \dots, a_n) ,     &  \text{if}~ \pi\in \ONCP(n), \\
 \vp_{< \pi}(a_1, \dots, a_n),  & \text{if} ~  \pi \notin \ONCP(n),  
\end{cases}
\end{align}
for all  $(a_1, \dots, a_n) \in \cA^n$.

The value
\[
\wt\vp^\AP_\pi (a_1^{(i_1)}\cdots a_n^{(i_n)}) = \wt\vp^\AP_{\pi \curlywedge \kappa(\bfi)} (a_1, \dots, a_n)
\]
is, by its construction, nothing but the value $\vp (a_1 a_2 \dots a_n)$ under the additional assumption that $(a_1, a_2,\dots, a_n)$ is $\AP$-free according to the ordered set partition $\pi \curlywedge \kappa(\bfi)$.\footnote{When $(\{a_i\}_{i\in P_1}, \dots, \{a_i\}_{i\in P_k})$ is $\AP$-free with $\pi=(P_1,P_2,\dots, P_k) \in \OSP(n)$, we say that $(a_1,\dots, a_n)$ is $\AP$-free according to $\pi$.} By the associativity, this is equivalent to computing $\vp (a_1 a_2 \dots a_n)$ in \emph{two steps}:  first we compute $\vp (a_1 a_2 \dots a_n)$ by assuming that $(a_1,a_2,\dots, a_n)$ is $\AP$-free according to $\pi$; then continue the computation by assuming that for every block $P=\{p_1,p_2,\dots, p_k\} \in \pi$ $(p_1 < \cdots < p_k)$ the tuple $(a_{p_1}, a_{p_2}, \dots, a_{p_k})$ is $\AP$-free according to the restricted partition $\Restr{\kappa(\bfi)}{P}.$ This reasoning extends formula \eqref{eq:universal_coefficients2} to 
\begin{align}
\wt\vp^\AP_\pi (a_1^{(i_1)},\cdots, a_n^{(i_n)}) 
&= \begin{cases}
\wt\vp^\AP_{(\pi)} (a_1^{(i_1)},\dots, a_n^{(i_n)})    +   \wt\vp^\AP_{< \pi}(a_1^{(i_1)}, \dots, a_n^{(i_n)}) ,     &  \text{if}~ \pi \in \ONCP(n), \\
 \wt\vp^\AP_{<\pi}(a_1^{(i_1)}, \dots, a_n^{(i_n)}) ,  & \text{if} ~  \pi \notin \ONCP(n). 
\end{cases}
\end{align}
Taking the sum over all the tuples $(i_1,\dots, i_n) \in [N]^n$ we obtain 
\begin{align} \label{eq:moments}
\wt\vp^\AP_\pi (N.a_1,\cdots, N.a_n) 
= \begin{cases}
\wt\vp^\AP_{(\pi)} (N.a_1,\dots, N.a_n)    +   \wt\vp^\AP_{< \pi}(N.a_1, \dots, N.a_n) ,     &  \text{if}~ \pi \in \ONCP(n), \\
 \wt\vp^\AP_{< \pi}(N.a_1, \dots, N.a_n) ,  & \text{if} ~ \pi \notin \ONCP(n).
\end{cases}
\end{align}
Because the coefficient of $N^{|\pi|}$ in $\wt\vp^\AP_{< \pi}(N.a_1, \dots, N.a_n)$ is zero (because each of $\wt\vp^\AP(N.a_P)$, $\wt\psi^\BP(N.a_P),$ $\wt\theta^\GP(N.a_P)$ has factor $N$), comparing  the coefficient of $N^{|\pi|}$ in  \eqref{eq:moments} yields 
\begin{align*}
K_\pi^{\AP} (a_1,\cdots, a_n) = \begin{cases}
\displaystyle \left[ \prod_{P \in \Out(\pi)   } K^{{\AP}}_{|P|}(a_P)\right] \left[ \prod_{P \in \AA(\pi)} K^{\BP}_{|P|}(a_P) \right]\left[  \prod_{P \in \BB(\pi)} K^{\GP}_{|P|}(a_P)\right] ,     &  \text{if}~ \pi \in \ONCP(n), \\[10mm]
 0,  & \text{if} ~ \pi \notin \ONCP(n).   
\end{cases}
\end{align*}
Combined with the general moment-cumulant formula in Theorem \ref{mc2} this entails the desired formula. 
\end{proof}

In the literature, moment-cumulant formulas are established for free \cite{Spe2}, c-free \cite{BLS}, monotone \cite{H-S,H-S2}, antimonotone, c-monotone \cite{Has3}  and Boolean independences \cite{SW97}. The antimonotone case is essentially the same as the monotone case. 
These formulas are all unified by Theorem \ref{thm:m-c-alpha}. 

\begin{cor}  Let $(\cA, \vp,\psi)$ be a 2-algebraic probability space. 
The multilinear functionals $K_{n}^{\CM(\vp, \psi)}:= \Restr{K_n^{\AP(\wh\vp,\wh\psi,\delta)}}{\cA^n}$ and $K_{n}^{\Mon(\psi)}:= \Restr{K_n^{\AP(\wh\psi,\wh\psi,\delta)}}{\cA^n}$ conincide with the c-monotone cumulants and monotone cumulants, respectively, and 
\begin{align}
&\vp(a_1 \cdots a_n) = \sum_{\pi \in \MP(n)}\frac{1}{|\pi|!}  \left[  \prod_{P \in \Out(\pi)}K_{|P|}^{\CM(\vp, \psi)}(a_P) \right] \left[  \prod_{P \in \Inn(\pi) }  K_{|P|}^{\Mon (\psi)}(a_P) \right], \label{eq:CM-mc}\\
&\psi(a_1 \cdots a_n) = \sum_{\pi \in \MP(n)}\frac{1}{|\pi|!} \prod_{P \in \pi}K_{|P|}^{\Mon(\psi)}(a_P). \label{eq:M-mc}
\end{align}
Similarly, $K_n^{\CAM(\vp, \psi)} := \Restr{K_n^{{\AP}(\wh\vp, \delta, \wh\psi)}}{\cA^n}$ and $K_n^{\AM(\psi)} := 
\Restr{K_n^{{\AP}(\wh\psi, \delta, \wh\psi)}}{\cA^n}$ are the $n$-th c-antimonotone cumulant and the $n$-th antimonotone cumulant, respectively, and the following moment-cumulant formulas hold:  
\begin{align}
&\vp(a_1 \cdots a_n) = \sum_{\pi \in \AMP(n)}\frac{1}{|\pi|!}      \left[   \prod_{P \in \Out(\pi)}   K_{|P|}^{\CAM(\vp, \psi)}(a_P) \right] \left[   \prod_{  P \in \Inn(\pi)  } K_{|P|}^{\AM(\psi)}(a_P) \right],   \label{eq:CAM-mc} \\
&\psi(a_1 \cdots a_n) = \sum_{\pi \in \AMP(n)}  \frac{1}{|\pi|!}    \prod_{P \in \pi }  K_{|P|}^{\AM(\psi)}(a_P). \label{eq:AM-mc}
\end{align}
Finally, $K_n^{\Boole(\vp)}:=K_n^{{\AP}(\wh\vp,\delta,\delta)}$ is the $n$-th Boolean cumulant and  
\begin{equation}
\vp(a_1 \cdots a_n) = \sum_{\pi \in \IP(n)} \prod_{P \in \pi}  K_{|P|}^{\Boole(\vp)}(a_P). 
\end{equation}
\end{cor}
\begin{proof}
Because $\Restr{K_n^{\GP(\wh\psi,\delta)}}{\cA^n}=0$, only the ordered noncrossing partitions $\pi$ with $\BB(\pi)=\emptyset$ contribute to the RHS of \eqref{eq503}. Such a class of ordered set partitions is exactly the set of the monotone partitions, and thus we conclude  \eqref{eq:CM-mc}. Formula \eqref{eq:M-mc} follows by setting $\vp = \psi$ in  \eqref{eq:CM-mc}. A recursive argument demonstrates that formulas  \eqref{eq:CM-mc} and \eqref{eq:M-mc} uniquely determine the multilinear functionals $K_{n}^{\CM(\vp, \psi)}$ and $K_{n}^{\Mon(\psi)}$ and hence they are identified with the c-monotone cumulants and monotone cumulants. The other formulas  can be proved in a similar manner. In the Boolean case, the summands do not depend on the order structure of $\ONCP(n)$, and therefore, 
the factor $\frac{1}{|\pi|!}$ vanishes after taking the partial sum over the possible orders on each set partition $\pi$.
\end{proof}

\begin{rem} 
Although $\MP(n)$ and $\AMP(n)$ are not identical as subsets of $\OSP(n)$, we actually have $K_{n}^{\CM(\vp, \psi)} = K_{n}^{\CAM(\vp, \psi)}$ and $K_{n}^{\Mon(\psi)} = K_{n}^{\AM(\psi)}$. This is because c-antimonotone independence with reversed order on the index set is c-monotone independence, and therefore the two notions of dot operations, i.e.\ sums of i.i.d.\ random variables,  are identical in distribution. 
\end{rem}

\begin{rem} 
C-monotone cumulants are only defined for single variables in \cite{Has3}. 
\end{rem}

The c-free and free moment-cumulant formulas \eqref{eq:cfree_mc} -- \eqref{eq:free_mc} can be retrieved by selecting the triplet of linear functionals $(\vp,\psi,\psi)$.

%%%%%%%%%%%%%%%%%%%%%%%%%%%%%%%%%%%%%%%%%%%%%%%%%%%%%%%%%%%%%%%%%%%%%%%%%%%%%%%%%%%%%%%%%%%%%%%%%%%%%%%%%%%%%%%%%%%%%%%%%%%%%%%%%%%%%%%%%%%%%%%%%%%%%%%%%%%%
\section{Convolutions}\label{sec:convolution}

%%%%%%%%%%%%%%%%%%%%%%%%%%%%%%%%%%%%%%%%%%%%%%%%%%%%%%%%%%%%%%%%%%%
\subsection{Additive and multiplicative convolutions of distributions}

%We start from the description of additive convolutions of probability measures. 
%This section will be useful for the reader to understand the idea of Section \ref{prod}.

Let $\comp[x]$ be the unital algebra generated from one indeterminate $x$. A \emph{distribution} is a unital linear functional on $\comp[x]$.  
Let $(\cA,\vp^1,\dots, \vp^m)$ be a (unital or not) $m$-algebraic probability space. 
The \emph{distribution} of $a\in \cA$ with respect to $(\cA, \vp^1,\dots, \vp^m)$ is the tuple $(\mu_a^1, \dots, \mu_a^m )$ of distributions determined by $\mu_a^i(x^n) = \vp^i(a^n)$ for all $n\in\N$ and $i\in[m]$.  Elements $a_i \in \cA, i \in I,$ are said to be \emph{identically distributed} if the distributions of them are the same. 

%If $a = a^\freeprod$ the there exists a probability measure $\rho$ such that $\int_\R t^n\, \rho_a(\di t) = \mu_a(x^n)$ for all $n \in \N$ but $\rho$ need not be unique.  when 

Suppose that a product $\square$ for (unital or not) $m$-algebraic probability spaces is given. For two $m$-algebraic probability spaces $(\cA_i,\vp^1_i, \dots, \vp^m_i), i=1,2$, let 
\[
(\cA, \vp^1,\dots, \vp^m) := (\cA_1, \vp_1^1, \dots, \vp_1^m) \square  (\cA_2, \vp_2^1, \dots, \vp_2^m),
\] where $\cA$ is the unital or nonunital free product of $\cA_1$ and $\cA_2$.  
Let $(\mu_i^1,\dots, \mu_i^m)$ be the distribution of $a_i \in \cA_i$ for $i=1,2$. For typical products $\square$ such as free product and Boolean product, the distribution of $a_1 + a_2$ with respect to $(\cA, \vp^1,\dots, \vp^m)$ is determined only by $(\mu_i^1,\dots, \mu_i^m), i=1,2$ and do not depend on other information on the elements  $a_1$ and $a_2$.  In this case  the distribution of $a_1 + a_2$ is called the \emph{additive convolution} of the distributions of $a_1$ and of $a_2$. 

In a similar manner the distribution of $a_1 a_2$ is called the \emph{multiplicative convolution} of the distributions of $a_1$ and of $a_2$ if it is determined by $(\mu_i^1,\dots, \mu_i^m), i=1,2$ only.  

Note that the convolutions defined as above are associative whenever the product for $m$-algebraic probability spaces is associative. 

\begin{rem} 
In the nonunital case, the above definition of multiplicative convolution often results in a trivial convolution, see e.g.\ \cite{Ber1,Fra}. As a remedy, the following modified definition is widely adopted in the literature (see e.g.\ \cite{Fra08,FHS20}):  the distribution of $x_1x_2$ is called the multiplicative convolution of the distributions of $x_1$ and $x_2$ if $(x_1-1, x_2-1)$ is independent. Note that the independence of $(x_1 -1,x_2-1)$ is not equivalent to the independence of $(x_1,x_2)$ in general (but still equivalent to the independence of $(x_1-1,x_2)$ in case of monotone independence). An interpretation of the  ``$-1$'' can be given for monotone  and Boolean independences  from the viewpoint of c-freeness, see \cite[Section 2]{Has4}. 
\end{rem}
 
In what follows the free additive and multiplicative convolutions of $\mu_1$ and $\mu_2$ are respectively denoted $\mu_1 \boxplus \mu_2$ and $\mu_1 \boxtimes \mu_2$.  The c-free additive and multiplicative convolutions of $(\mu_1, \nu_1)$ and $(\mu_2,\nu_2)$ are respectively denoted $(\mu_1,\nu_1)  \boxplus (\mu_2,\nu_2)$ and $(\mu_1,\nu_1)  \boxtimes (\mu_2,\nu_2)$. 
Note that the left  component of $(\mu_1,\nu_1) \boxplus (\mu_2,\nu_2)$ depends on $\mu_1,\mu_2,\nu_1,\nu_2$ and the right component of   $(\mu_1,\nu_1) \boxplus (\mu_2,\nu_2)$ is the free additive convolution $\nu_1 \boxplus \nu_2$. The left component will be denoted $\cfconv{\mu_1}{\nu_1}{\nu_2}{\mu_2}$.  
Similarly, $(\mu_1,\nu_1)  \boxtimes (\mu_2,\nu_2)$ is denoted $(\mcfconv{\mu_1}{\nu_1}{\nu_2}{\mu_2}, \nu_1 \boxtimes \nu_2)$.

The \emph{additive $\AP$-free convolution} and the \emph{multiplicative $\AP$-free convolution} of $(\lambda_1,\mu_1,\nu_1)$ and $(\lambda_2,\mu_2,\nu_2)$  are denoted $(\lambda_1,\mu_1,\nu_1)\subsetplus(\lambda_2,\mu_2,\nu_2)$ and $(\lambda_1,\mu_1,\nu_1)\subsettimes(\lambda_2,\mu_2,\nu_2)$, respectively. By the associativity of $\AP$-free product, these convolutions are associative. 
Also, by Proposition \ref{rem:cfree_alpha}, it holds that 
\[
(\lambda_1,\mu_1,\nu_1)\subsetplus(\lambda_2,\mu_2,\nu_2) = (\cfconv{\lambda_1}{\nu_1}{\mu_2}{\lambda_2}, \cfconv{\mu_1}{\nu_1}{\mu_2}{\mu_2}, \cfconv{\nu_1}{\nu_1}{\mu_2}{\nu_2})
\]
and 
\[
(\lambda_1,\mu_1,\nu_1)\subsettimes(\lambda_2,\mu_2,\nu_2) = (\mcfconv{\lambda_1}{\nu_1}{\mu_2}{\lambda_2}, \mcfconv{\mu_1}{\nu_1}{\mu_2}{\mu_2}, \mcfconv{\nu_1}{\nu_1}{\mu_2}{\nu_2}).  
\]
The binary operations $\subsetplus$ and $\subsettimes$ restricted to the second and third components are called the additive and multiplicative \emph{$\BGP$-free convolutions}, respectively, and will be denoted by the same symbols, e.g.\  
$
(\mu_1,\nu_1)\subsetplus(\mu_2,\nu_2) = (\cfconv{\mu_1}{\nu_1}{\mu_2}{\mu_2}, \cfconv{\nu_1}{\nu_1}{\mu_2}{\nu_2}). 
$
 
Because the $\AP$-free product unifies other products of algebraic probability spaces as listed in Subsection \ref{subsec:main}, the $\AP$-free convolution  unifies the corresponding convolutions; for instance,  c-free convolution and monotone convolution are included in the ways 
\[
(\lambda_1,\mu_1,\mu_1)\subsetplus(\lambda_2,\mu_2,\mu_2) = (\cfconv{\lambda_1}{\mu_1}{\mu_2}{\lambda_2}, \mu_1 \boxplus \mu_2, \mu_1 \boxplus \mu_2) = ((\lambda_1,\mu_1) \boxplus (\lambda_2,\mu_2), \mu_1 \boxplus \mu_2)
\]
and 
\[
(\mu_1,\mu_1,\delta_0)\subsetplus(\mu_2,\mu_2,\delta_0) = (\mu_1\trr\mu_2, \mu_1\trr \mu_2, \delta_0).  
\]
These formulas show that the associativity of $\AP$-free convolution is transferred to the associativity of c-free convolution and monotone convolution.  

In the following subsections we summarize generating functional machineries  for calculating c-free convolutions and hence $\AP$-free convolutions. 
Results will be used later when we analyze $\AP$-free cumulants and compute the central limit distribution for $\AP$-freeness, see Section \ref{sec:generating_function}.

%%%%%%%%%%%%%%%%%%%%%%%%%%%%%%%%%%%%%%%%%%%%%%%%%%%%%%%%%%%%%%%%%%%%%%%%%%%
\subsection{Analytic tools for describing convolutions}

The results in this section are based on \cite{Bel,BLS}. Let $\mu, \nu$ be distributions. Note that $(\comp[x],\mu,\nu)$ can be regarded as a unital 2-algebraic probability space. The \emph{Cauchy transform} of $\mu$ is the formal Laurent series 
\begin{equation}
G_\mu(z) = \sum_{n\ge0} \frac{\mu(x^n)}{z^{n+1}}
\end{equation}
and the \emph{reciprocal Cauchy transform} of $\mu$ is the formal Laurent series 
\begin{equation}
F_\mu(z) = \frac{1}{G_\mu(z)} = z  -b_1(\mu) - \frac{b_2(\mu)}{z} - \frac{b_3(\mu)}{z^2} - \cdots, 
\end{equation}
where $b_n(\mu)$ are polynomials on $\mu(x), \dots, \mu(x^n)$ for each $n\in \N$. These coefficients $b_n(\mu)$ coincide with the Boolean cumulants $K_n^{\Boole(\mu)}(x)$ of $x$.  The reciprocal Cauchy transform has a formal compositional inverse Laurent series of the form 
\[
F^{-1}_\mu(z) = z  +  \sum_{n\ge1}\frac{r_n(\mu)}{z^{n-1}}. 
\]
The coefficients $r_n(\mu)$ coincide with the free cumulants $K_n^{\Free(\mu)}(x)$ of $x$. The formal Laurent series $\phi_\nu(z)$ (called the \emph{Voiculescu transform} of $\nu$) and $\phi_{(\mu, \nu)}(z)$ determined by 
\begin{gather}
F_{\nu}(z) = z - \phi_{\nu}(F_{\nu}(z)),  \label{rel13}\\
F_{\mu}(z) = z - \phi_{(\mu, \nu)}(F_{\nu}(z)) \label{rel14}   
\end{gather}
are useful for computing additive c-free convolution because of they are generating functions of free cumulants and c-free cumulants: 
\[
\phi_{\nu}(z) = \sum_{n\ge1}\frac{K_n^{\Free(\nu)}(x)}{z^{n-1}} \quad \text{and} \quad \phi_{(\mu,\nu)}(z) = \sum_{n\ge1}\frac{K_n^{\CF(\mu,\nu)}(x)}{z^{n-1}}. 
\]
Consequently the additivity properties hold: 
\begin{gather}
\phi_{\nu_1 \boxplus \nu_2}(z) = \phi_{\nu_1}(z) + \phi_{\nu_2} (z), \label{eq11111}\\
\phi_{(\mu_1, \nu_1) \boxplus (\mu_2, \nu_2)}(z) = \phi_{(\mu_1, \nu_1)}(z) + \phi_{(\mu_2, \nu_2)}(z) \label{cum45}. 
\end{gather}
Note that (\ref{eq11111}) has the alternative form 
\begin{equation}\label{rel15}
F_{\nu_1 \boxplus \nu_2} ^{-1}(z) = F_{\nu_1}^{-1}(z) + F_{\nu_2}^{-1}(z) - z.   
\end{equation} 
%Because \eqref{eq11111} and \eqref{cum45} are symmetric, i

The formal series 
\[
R_{(\mu,\nu)} (z) := \phi_{(\mu, \nu)}\left(\frac1{z}\right) \quad \text{and} \quad R_\nu(z) : =\phi_{\nu}\left(\frac1{z}\right)
\]
 (the latter called the \emph{$R$-transform} of $\nu$) can be replacements for $\phi_{(\mu,\nu)}(z), \phi_\nu(z)$. These formal series are determined uniquely by 
\begin{gather}
\frac{1}{G_{\nu}(z)} = z - R_{\nu}(G_{\nu}(z)), \label{rel1} \\
\frac{1}{G_{\mu}(z)} = z - R_{(\mu, \nu)}(G_{\nu}(z)). \label{rel2}  
\end{gather}

The following formula was proved in \cite[Proposition 3]{Bel} for probability measures possibly with unbounded supports. 
We state the result  in the sense of formal Laurent series. 
\begin{prop}\label{eq124} 
Let $\mu_i,\nu_i~(i=1,2)$ be distributions. Then the distribution $\cfconv{\mu_1}{\nu_1}{\nu_2}{\mu_2}$ is characterized by 
\[
F_{\subcfconv{\mu_1}{\nu_1}{\nu_2}{\mu_2}}(z) = F_{\mu_1} \circ F_{\nu_1} ^{-1} \circ F_{\nu_1 \boxplus \nu_2}(z) + F_{\mu_2} \circ F_{\nu_2} ^{-1} \circ F_{\nu_1 \boxplus \nu_2} - F_{\nu_1 \boxplus \nu_2}(z).  
\]

\end{prop}
\begin{cor}\label{eq125} 
Let $\mu_1, \mu_2, \nu_1$ be distributions. Then the distribution $\cfconv{\mu_1}{\nu_1}{\mu_2}{\mu_2}$ is characterized by 
\[
F_{\subcfconv{\mu_1}{\nu_1}{\mu_2}{\mu_2}}(z) = F_{\mu_1} \circ F_{\nu_1} ^{-1} \circ F_{\nu_1 \boxplus \mu_2}(z). 
\]
\end{cor}

%%%%%%%%%%%%%%%%%%%%%%%%%%%%%%%%%%%%%%%%%%%%%%%%%%%%%%
\subsection{Analytic tools for multiplicative c-free convolution}
Results in this subsection will not be used later but are summarized here as a reference for future work. Let $\mu,\nu$ be distributions. 
The formal series 
\[
\eta_{\mu}(z)= 1 - \frac{z}{G_\mu (\frac{1}{z})} = \sum_{n\ge1} K_n^{\Boole(\mu)}(x) z^n
\]
 is useful for describing multiplicative c-free convolution. 
 %As was used effectively in \cite{Ber1,Has4}, t
% to characterize $\mu_1 {}_{\nu_1}\!\! \boxtimes_{\mu_2}\mu_2$. 
We define the formal series 
\[
\widetilde{R}_{(\mu, \nu)}(z):=z R_{(\mu, \nu)}(z) \quad \text{and} \quad \widetilde{R}_{\mu}(z):=z R_{\mu}(z)
\]
which were used in \cite{P-W} with different notation. 
Formulas (\ref{rel1}) and (\ref{rel2}) are translated into 
\begin{align}
&\widetilde{R}_\nu \left( \frac{z}{1 - \eta_\nu(z)}\right) = \frac{\eta_\nu(z)}{1 - \eta_\nu(z)}, \label{eq01}\\
&\widetilde{R}_{(\mu, \nu)}\left(\frac{z}{1 - \eta_\nu(z)}\right) = \frac{\eta_\mu(z)}{1 - \eta_\nu(z)}. \label{eq02}   
\end{align}

Multiplicative c-free convolution has been characterized in \cite{P-W} in the following way. 
Assume that $\nu(x)\ne0$ which guarantees that the compositional formal inverse series $\widetilde{R}^{-1}_\nu(z)$ exists.  Let $T_{(\mu, \nu)}(z), T_\nu(z)$ be defined by 
\[
T_{(\mu, \nu)}(z)= \frac{\widetilde{R}_{(\mu, \nu)}(\widetilde{R}^{-1}_\nu (z))}{\widetilde{R}^{-1} _\nu (z)}  \quad \text{and}  \quad T_\nu(z)= T_{(\nu, \nu)}(z) = \frac{z}{\widetilde{R}^{-1}_\nu (z)}. 
\] 
The multiplicative c-free convolution of $(\mu_1, \nu_1)$ and $(\mu_2, \nu_2)$ with $\nu_1(x), \nu_2(x) \neq 0$ is characterized by 
\begin{align}
&T_{(\submcfconv{\mu_1}{\nu_1}{\nu_2}{\mu_2},  \nu_1 \boxtimes \nu_2)}(z) = T_{(\mu_1, \nu_1)}(z) T_{(\mu_2, \nu_2)}(z), \label{eq:c-freeT}\\ 
&T_{\nu_1 \boxtimes \nu_2}(z) = T_{\nu_1}(z) T_{\nu_2}(z). \label{eq04}
\end{align}

Formula \eqref{eq:c-freeT} can be expressed only in terms of $\eta_\mu$ in the following way. 

\begin{prop}\label{prop002} 
Let $\mu_1,\mu_2,\nu_1,\nu_2$ be distributions such that $\nu_1(x), \nu_2(x) \neq 0$. Then 
\begin{equation}\label{eq003} 
\eta_{\submcfconv{\mu_1}{\nu_1}{\nu_2}{\mu_2}}(z) = \frac{\eta_{\mu_1} \circ \eta_{\nu_1}^{-1} \circ \eta_{\nu_1 \boxtimes \nu_2}(z) \eta_{\mu_2} \circ \eta_{\nu_2}^{-1} \circ \eta_{\nu_1 \boxtimes \nu_2}(z)}{\eta_{\nu_1 \boxtimes \nu_2}(z)}. 
\end{equation}
\end{prop}
\begin{proof}
Note first that $\widetilde{R}_{\nu_i}(z), \eta_{\nu_i}(z)$ ($i=1,2$), $\widetilde{R}_{\nu_1 \boxtimes \nu_2}(z)$ and $\eta_{\nu_1 \boxtimes \nu_2}(z)$ all have formal compositional inverse series because $\nu_i(x) \neq 0$ and $(\nu_1 \boxtimes \nu_2)(x)= \nu_1(x)\nu_2(x) \ne0$. 
From \eqref{eq:c-freeT} and (\ref{eq04}) it follows that 
\begin{equation}\label{eq05}
z \widetilde{R}_{(\submcfconv{\mu_1}{\nu_1}{\nu_2}{\mu_2}, \nu_1 \boxtimes \nu_2)}(\widetilde{R}^{-1}_{\nu_1 \boxtimes \nu_2}(z)) 
     = \widetilde{R}_{(\mu_1, \nu_1)}(\widetilde{R}^{-1}_{\nu_1}(z)) \widetilde{R}_{(\mu_2, \nu_2)}(\widetilde{R}^{-1}_{\nu_2}(z)).   
\end{equation}
We define new variables $u$, $v$ and $w$ by 
\begin{align}
&\widetilde{R}^{-1}_{\nu_1}(z) = \frac{u}{1 -  \eta_{\nu_1}(u)}, \\
&\widetilde{R}^{-1}_{\nu_2}(z) = \frac{v}{1 -  \eta_{\nu_2}(v)}, \\
&\widetilde{R}^{-1}_{\nu_1 \boxtimes \nu_2}(z) = \frac{w}{1 - \eta_{\nu_1 \boxtimes \nu_2}(w)}. 
\end{align} 
These equalities, combined with (\ref{eq01}) and (\ref{eq02}),  entail 
\begin{equation}
z = \frac{\eta_{\nu_1}(u)}{1 - \eta_{\nu_1}(u)}  = \frac{\eta_{\nu_2}(v)}{1 - \eta_{\nu_2}(v)} = \frac{\eta_{\nu_1 \boxtimes \nu_2}(w)}{1 - \eta_{\nu_1 \boxtimes \nu_2}(w)},  
\end{equation}
and therefore  $\eta_{\nu_1 \boxtimes \nu_2}(w) = \eta_{\nu_1}(u) = \eta_{\nu_2}(v)$. By formula 
(\ref{eq05}) we obtain 
\[
z \frac{\eta_{\submcfconv{\mu_1}{\nu_1}{\nu_2}{\mu_2}}(w)}{1 -  \eta_{\nu_1 \boxtimes \nu_2}(w)} = \frac{\eta_{\mu_1}(u)}{1 - \eta_{\nu_1}(u)}\frac{\eta_{\mu_2}(v)}{1 - \eta_{\nu_2}(v)}. 
\]
Since 
\[
\frac{z}{1-\eta_{\nu_1 \boxtimes \nu_2}(w)} = \frac{z^2}{\eta_{\nu_1 \boxtimes \nu_2}(w)} = \frac{\eta_{\nu_1 \boxtimes \nu_2}(w)}{(1 - \eta_{\nu_1}(u))(1 - \eta_{\nu_2}(v))}, 
\]
the desired formula \eqref{eq003} follows. 
\end{proof}

\begin{cor}
For any distributions $\mu_1,\mu_2,\nu_1$ with $\mu_2(x), \nu_1(x)\ne0$ we have 
\begin{equation}\label{eq004}
\eta_{\submcfconv{\mu_1}{\nu_1}{\mu_2}{\mu_2}}(z) = \eta_{\mu_1} \circ \eta_{\nu_1}^{-1} \circ \eta_{\nu_1 \boxtimes \mu_2}(z). 
\end{equation}
\end{cor}

\begin{rem}
It is worth noting a similarity (at least at a formal level) between Proposition \ref{eq124} and Proposition \ref{prop002}. 
With a formal function 
$f_\mu = \log \circ \eta_\mu \circ \exp$,  Proposition \ref{prop002} reads 
\[
f_{\mu_1 {}_{\nu_1}\!\boxtimes_{\nu_2} \mu_2} = f_{\mu_1} \circ f_{\nu_1} ^{-1} \circ f_{\nu_1 \boxtimes \nu_2} + f_{\mu_2} \circ f_{\nu_2} ^{-1} \circ f_{\nu_1 \boxtimes \nu_2} - f_{\nu_1 \boxtimes \nu_2},  
\]
which is of the same form as Proposition \ref{eq124}. Because of this, additive c-free convolution and multiplicative one show some parallelism. 
Note that similar parallelism between additive and multiplicative convolutions is analyzed in \cite{AA17} in details for free, Boolean and monotone convolutions.  
\end{rem}

%%%%%%%%%%%%%%%%%%%%%%%%%%%%%%%%%%%%%%%%%
\subsection{Some heuristics}  This subsection briefly exposes a certain heuristics on the discovery of $\AP$-free product as an associative product. 

 Considering the simple formula in Corollary \ref{eq125}, 
let us seek for an associative binary operation $\subsetplus$ for pairs of distributions of the form $(\mu_1,\nu_1) \subsetplus (\mu_2,\nu_2) = (\cfconv{\mu_1}{\nu_1}{\mu_2}{\mu_2}, \rho)$, where $\rho = \rho(\mu_1, \mu_2, \nu_1, \nu_2)$ is 
a distribution depending on $\mu_1, \mu_2, \nu_1, \nu_2$.  Assume that $\subsetplus$ is associative.  The associativity implies that 
\[
\cfconv{  (\cfconv{\mu_1}{\nu_1}{\mu_2}{\mu_2})   }{\rho} {\mu_3} {\mu_3}   = \cfconv{\mu_1}{\nu_1}{(\subcfconv{\mu_2}{\nu_2}{\mu_3}{\mu_3})}  { (\cfconv{\mu_2}{\nu_2}{\mu_3}{\mu_3})} . 
\]
By the way, Proposition \ref{eq124} implies that 
\[
\begin{split}
F_{\subcfconv{  (\subcfconv{\mu_1}{\nu_1}{\mu_2}{\mu_2})   }{\rho} {\mu_3} {\mu_3}   }&= F_{\subcfconv{\mu_1}{\nu_1}{\mu_2}{\mu_2}} \circ F_{\rho}^{-1} \circ F_{\rho \boxplus \mu_3} \\ 
&= F_{\mu_1} \circ F_{\nu_1}^{-1} \circ F_{\nu_1 \boxplus \mu_2} \circ F_{\rho}^{-1} \circ F_{\rho \boxplus \mu_3}
\end{split}
\]
and 
\[
\begin{split}
F_{ \subcfconv{\mu_1}{\nu_1}{(\subcfconv{\mu_2}{\nu_2}{\mu_3}{\mu_3})}  { (\subcfconv{\mu_2}{\nu_2}{\mu_3}{\mu_3})} } 
&= F_{\mu_1} \circ F_{\nu_1}^{-1} \circ F_{\nu_1 \boxplus (\subcfconv{\mu_2}{\nu_2}{\mu_3}{\mu_3})}.
\end{split}
\]
Therefore, it holds that $F_{\nu_1 \boxplus \mu_2} \circ F_{\rho}^{-1} \circ F_{\rho \boxplus \mu_3} = F_{\nu_1 \boxplus (\subcfconv{\mu_2}{\nu_2}{\mu_3}{\mu_3})}$, 
or equivalently, 
\begin{equation}\label{eq21}
F_{\rho \boxplus \mu_3}^{-1}\circ F_{\rho} \circ F_{\nu_1 \boxplus \mu_2}^{-1} =  F_{\nu_1 \boxplus (\subcfconv{\mu_2}{\nu_2}{\mu_3}{\mu_3})  }^{-1}.
\end{equation} 
The left hand side is  
\begin{equation}\label{eq22}
\begin{split}
F_{\rho \boxplus \mu_3}^{-1}\circ F_{\rho} \circ F_{\nu_1 \boxplus \mu_2}^{-1}
&= (z + F_{\mu_3}^{-1}\circ F_\rho - F_\rho) \circ F_{\nu_1 \boxplus \mu_2}^{-1} \\
&= F_{\nu_1 \boxplus \mu_2}^{-1} + F_{\mu_3}^{-1}\circ F_\rho \circ F_{\nu_1 \boxplus \mu_2}^{-1} -  F_\rho \circ F_{\nu_1 \boxplus \mu_2}^{-1}
\end{split}
\end{equation}
by using (\ref{rel15}). On the other hand we have 
\[
\begin{split}
F_{\nu_1 \boxplus (\subcfconv{\mu_2}{\nu_2}{\mu_3}{\mu_3})}^{-1} 
&= F_{\nu_1} ^{-1} + F_{(\subcfconv{\mu_2}{\nu_2}{\mu_3}{\mu_3}) }^{-1} -z \\
&= F_{\nu_1} ^{-1} + F_{\nu_2 \boxplus \mu_3}^{-1} \circ F_{\nu_2} \circ F_{\mu_2}^{-1} - z \\ 
&= F_{\nu_1} ^{-1} + F_{\mu_2}^{-1} + F_{\mu_3}^{-1} \circ F_{\nu_2} \circ F_{\mu_2}^{-1} - F_{\nu_2} \circ F_{\mu_2}^{-1} - z \\
&= F_{\nu_1 \boxplus \mu_2} ^{-1} + F_{\mu_3}^{-1} \circ F_{\nu_2} \circ F_{\mu_2}^{-1} - F_{\nu_2} \circ F_{\mu_2}^{-1}.  
\end{split}
\]
Combining (\ref{eq21}) and (\ref{eq22}) yields 
\begin{equation}\label{eq23}
F_{\mu_3}^{-1} \circ F_{\nu_2} \circ F_{\mu_2}^{-1} - F_{\nu_2} \circ F_{\mu_2}^{-1} = F_{\mu_3}^{-1}\circ F_\rho \circ F_{\nu_1 \boxplus \mu_2}^{-1} -  F_\rho \circ F_{\nu_1 \boxplus \mu_2}^{-1}. 
\end{equation}

In view of Corollary \ref{eq125}, formula \eqref{eq23} is satisfied if we define 
\[
\rho = \cfconv{\nu_1}{\nu_1}{\mu_2}{\nu_2}.  
\]  
The above discussion implies that 
\begin{equation}\label{eq:associativity_cfconv1}
\cfconv{    (\cfconv{\mu_1}{\nu_1}{\mu_2}{\mu_2})   }  {  (\subcfconv{\nu_1}{\nu_1}{\mu_2}{\nu_2})  }  {\mu_3}  {\mu_3} =  
\cfconv{\mu_1}{\nu_1}  {(\subcfconv{\mu_2}{\nu_2}{\mu_3}{\mu_3})} {( \cfconv{\mu_2}{\nu_2}{\mu_3}{\mu_3})  }. 
\end{equation}
If we replace $\mu_1$, $\mu_2$, $\mu_3$, $\nu_1$ and $\nu_2$ respectively with $\nu_3$, $\nu_2$, $\nu_1$, $\mu_3$ and $\mu_2$, then we have 
\begin{equation}\label{eq:associativity_cfconv2}
\cfconv{  (\cfconv{\nu_1}{\nu_1}{\mu_2}{\nu_2})  }  {(\subcfconv{\nu_1}{\nu_1}{\mu_2}{\nu_2})   }   {\mu_3} {\nu_3}
=   \cfconv{ \nu_1}{\nu_1}{   (\subcfconv{\mu_2}{\nu_2}{\mu_3}{\mu_3})   } { (\cfconv{\nu_2}{\nu_2}{\mu_3}{\nu_3}) }.  
\end{equation}
Two formulas \eqref{eq:associativity_cfconv1} and \eqref{eq:associativity_cfconv2} exactly mean the associative law of the binary operation  
\[
(\mu_1, \nu_1)\subsetplus (\mu_2, \nu_2) =(\cfconv{\mu_1}{\nu_1}{\mu_2}{\mu_2}, \cfconv{\nu_1}{\nu_1}{\mu_2}{\nu_2}). 
\]

In a similar manner, one sees that the binary operation $\subsetplus$ for triplets of distributions 
\begin{equation}\label{eq:APconvolution}
(\lambda_1, \mu_1, \nu_1)\subsetplus (\lambda_2, \mu_2, \nu_2) =(\cfconv{\lambda_1}{\nu_1}{\mu_2}{\lambda_2}, \cfconv{\mu_1}{\nu_1}{\mu_2}{\mu_2}, \cfconv{\nu_1}{\nu_1}{\mu_2}{\nu_2}) 
\end{equation}
is associative. This is more or less how the $\BGP$-free product and $\AP$-free product were discovered. Once the associativity holds at the level of convolutions, it is quite natural to expect associativity for the product of unital 3-algebraic probability spaces, and indeed it is the case as demonstrated in Theorem \ref{thm:main_associativity}.

\section{The generating function of $\AP$-free cumulants for single variables}\label{sec:generating_function}

%%%%%%%%%%%%%%%%%%%%%%%%%%%%%%%%%%%%%%%%%%%%%%%%%%%%%%%%%%%%%%%%%%
\subsection{The $\AP$-free cumulants and one-parameter convolution groups}
The relation between moments and $\AP$-free cumulants is described in Theorem \ref{thm:m-c-alpha} in the language of ordered set partitions. In some situations, e.g.\ when analyzing the central limit theorem, generating functions are more useful. It is well known that moment generating functions and free / c-free cumulant generating functions are related via basic algebraic operations of  composition, addition,  multiplication and quotient, see Section \ref{sec:convolution}. The Boolean cumulant generating function does not even require composition. On the other hand, relating the monotone  (or c-monotone) cumulant generating function with moment generating functions seems to inevitably involve differential equations \cite{Has3}. It is then natural to expect that differential equations will be useful for $\AP$-freeness too because $\AP$-free cumulants generalize monotone cumulants.    

In this section we relate moment generating functions with $\AP$-free cumulants (and $\BP$-free  and $\GP$-free cumulants) generating functions in the single variable case. 
Let $(\cA,\vp,\psi,\theta)$ be a unital 3-algebraic probability space and $a\in\cA$. Let $(\lambda,\mu,\nu)$ be the distribution of $a$ regarding $(\vp,\psi,\theta)$, i.e.\ $\lambda, \mu, \nu$ be unital linear functionals on $\comp[x]$ such that 
\[
\lambda(x^n) = \vp(a^n), \quad \mu(x^n) = \psi(a^n) \quad \text{and} \quad \nu(x^n) = \theta(a^n), \quad n\in \N.  
\]

For $N\in \{0\}\cup\N$ we denote by $(\lambda_N,\mu_N,\nu_N)$ the distribution of $N.a$ concerning the unital 3-algebraic probability space $(\cU^\Free, \wt\vp^\AP, \wt\psi^\BP, \wt\theta^\GP)$ as defined in Subsection \ref{subsec:alpha_spread}: 
\[
\lambda_N(x^n) = \wt\vp^\AP((N.a)^n), \quad \mu_N(x^n) = \wt\psi^\BP((N.a)^n) \quad \text{and} \quad \nu_N(x^n) = \wt\theta^\GP((N.a)^n), \quad n\in \N.  
\]
 By construction $(\lambda_1,\mu_1,\nu_1)$ coincides with $(\lambda,\mu,\nu)$.  
 Because of the associativity of $\AP$-free convolution and because $N.a$ is the sum of $\AP$-free, identically distributed random variables, the semigroup formula 
 \begin{equation}\label{eq:semigroup}
 (\lambda_M,\mu_M,\nu_M) \subsetplus (\lambda_N,\mu_N,\nu_N) = (\lambda_{M+N},\mu_{M+N},\nu_{M+N})
\end{equation}
holds for all $M,N \in \N \cup\{0\}$. Note that  
\begin{equation} \label{eq:initial}
\lambda_0 = \mu_0 =\nu_0 =\delta_0, 
\end{equation}
where $\delta_0$ is the unital linear functional on $\comp[x]$ defined by $\delta_0(x^n)= \delta_{0,n}, n\in \N\cup\{0\}$.  It can be checked that $(\delta_0,\delta_0,\delta_0)$ is the unit for the convolution $\subsetplus$. 
 
 Because $\lambda_N(x^n),\mu_N(x^n),\nu_N(x^n)$ are polynomials on $N$ without a constant term, we can extend $N$ to real numbers $t$ and then define a distribution $(\lambda_t,\mu_t,\nu_t)$ for $t\in \real$. By the polynomiality, formula \eqref{eq:semigroup} extends to
  \begin{equation}\label{eq:semigroup_real}
 (\lambda_t,\mu_t,\nu_t) \subsetplus (\lambda_s,\mu_s,\nu_s) = (\lambda_{t+s},\mu_{t+s},\nu_{t+s}), \qquad t,s \in \real. 
\end{equation}
 We call  $\{(\lambda_t,\mu_t,\nu_t)\}_{t\in\real}$ the \emph{one-parameter $\AP$-free convolution group associated with $(\lambda,\mu,\nu)$}.

By the definition of cumulants (see Definition \ref{def:cumulants}), we have
 \begin{equation}\label{eq:cumulants_convolution_group}
 \Der \lambda_t(x^n) = K_n^{{\AP}}(a), \quad  \Der  \mu_t(x^n) = K_n^{\BP}(a), \quad \text{and} \quad  \Der \nu_t(x^n) = K_n^{\GP}(a). 
 \end{equation}

\begin{rem} When $(\cA,\vp,\psi,\theta)$ is a 3-$\ast$-algebraic probability space and $a$ is self-adjoint then the distributions $\lambda,\mu,\nu$ are states on $\comp[x]$ if we equip $\comp[x]$ with the involution $\ast$ determined by $x^\ast =x$.  Also, $\lambda_N,\mu_N,\nu_N$ are states on $\comp[x]$ for all $N \in\{0\}\cup\N$. However, for other $t\in\real$, $\lambda_t,\mu_t,\nu_t$ may fail to be states.  
\end{rem}

We derive differential equations for the reciprocal Cauchy transforms $F_{\lambda_t}, F_{\mu_t}, F_{\nu_t}$ (see Section \ref{sec:convolution} for the definition). To begin with, we define the formal (i.e.\ coefficient-wise) partial derivatives of the reciprocal Cauchy transforms  
\[
A_\lambda(z):=\pDer F_{\lambda_t} (z),\quad  B_\mu(z):=\pDer F_{\mu_t} (z) \quad \text{and} \quad  C_\nu(z):=\pDer F_{\nu_t} (z), 
\]  
which are, thanks to \eqref{eq:cumulants_convolution_group}, nothing but the cumulant generating functions of the forms  
\begin{align}
&A_\lambda(z) = - \frac{1}{G_{\lambda_0}(z)^2} \pDer G_{\lambda_t}(z) = - \sum_{n\ge1} \frac{K_n^{{\AP}}(a)}{z^{n-1}}, \\
& B_\mu(z) =  - \sum_{n\ge1} \frac{K_n^{{\BP}}(a)}{z^{n-1}}, \qquad C_\nu(z) =  - \sum_{n\ge1} \frac{K_n^{{\GP}}(a)}{z^{n-1}}. 
\end{align}

\begin{prop}  \label{prop:diff_eq}
The following differential equations hold:  
\begin{align}
&\frac{\partial F_{\lambda_t}}{\partial t} = A_\lambda \circ F_{\mu_t}+  (C_\nu \circ F_{\mu_t})\cdot \left[ \frac{\partial F_{\lambda_t}}{\partial z} -1\right] , \label{eq41} \\ 
&\frac{\partial F_{\lambda_t}}{\partial t} = A_\lambda \circ F_{\nu_t} + (B_\mu \circ F_{\nu_t})\cdot \left[ \frac{\partial F_{\lambda_t}}{\partial z} -1\right] ,  \label{eq42} \\
&\frac{\partial F_{\lambda_t}}{\partial t} = \frac{(B_\mu \circ F_{\nu_t}) \cdot (A_\lambda \circ F_{\mu_t}) - (A_\lambda \circ F_{\nu_t})\cdot (C_\nu \circ F_{\mu_t})}{B_\mu \circ F_{\nu_t} - C_\nu \circ F_{\mu_t}}. \label{eq420} 
\end{align}
Note that \eqref{eq420} holds only when the denominator is nonzero.  
\end{prop}
\begin{proof}
By Proposition \ref{eq124},  
\begin{equation}\label{eq400}
 F_{\lambda_{t+s}} =F_{\subcfconv{\lambda_t}{\nu_t}{\mu_s}{\lambda_s}} = F_{\lambda_t} \circ F_{\nu_t} ^{-1} \circ F_{\nu_t \boxplus \mu_s} + F_{\lambda_s} \circ F_{\mu_s} ^{-1} \circ F_{\nu_t \boxplus \mu_s} - F_{\nu_t \boxplus \mu_s}. 
\end{equation} 
Note that \eqref{eq:initial} asserts that $F_{\lambda_0}(z) = F_{\mu_0}(z) = F_{\nu_0}(z)=z$.  
Differentiating (\ref{eq400}) with respect to $t$ at $t=0$ yields   %set $t=0$, replace $s$ by $t$ and then obtain   
\begin{align}
\frac{\partial F_{\lambda_s}}{\partial s} &= A_\lambda \circ F_{\mu_s} + \pDer F_{\nu_t}^{-1} \circ F_{\nu_t \boxplus \mu_s} + \frac{\partial F_{\lambda_s}}{\partial z} \cdot \frac{\partial (F_{\mu_s}^{-1})}{\partial z} \circ F_{\mu_s} \cdot \pDer F_{\nu_t \boxplus \mu_s}  - \pDer F_{\nu_t \boxplus \mu_s}   \notag \\
&= A_\lambda \circ F_{\mu_s} + \left(\pDer  F_{\nu_t}^{-1}\right) \circ F_{\mu_s} + \frac{\partial F_{\lambda_s}}{\partial z} \cdot \frac{\partial (F_{\mu_s}^{-1})}{\partial z} \circ F_{\mu_s} \cdot \pDer F_{\nu_t \boxplus \mu_s}.   \label{eq:intermediate}
\end{align}
Recall from (\ref{rel15}) that the relation $F_{\nu_t \boxplus \mu_s} ^{-1} = F_{\nu_t} ^{-1} + F_{\mu_s} ^{-1}-z$ holds. Some simple calculations yield  
\begin{equation} \label{eq:free_conv}
\pDer F_{\nu_t \boxplus \mu_s} = (C_\nu \circ F_{\mu_s})\cdot\frac{\partial F_{\mu_s}}{\partial z}.  
\end{equation}
Applying the chain rule for the identities $F_{\mu_s}^{-1} \circ F_{\mu_s} =z$ and $F_{\nu_t} \circ F_{\nu_t}^{-1} =z$ implies 
\begin{align}
&\frac{\partial (F_{\mu_s}^{-1})}{\partial z} \circ F_{\mu_s} = \dfrac{1}{\frac{\partial F_{\mu_s}}{\partial z}}, \\
&\pDer  F_{\nu_t}^{-1} = - C_\nu.  \label{eq:chain_rule}
\end{align} 
Substituting \eqref{eq:free_conv} -- \eqref{eq:chain_rule} into \eqref{eq:intermediate}, we obtain the first desired formula (\ref{eq41}). The second formula (\ref{eq42}) follows by replacing $(\mu_t, \nu_t)$ with $(\nu_t, \mu_t)$; this is allowed because the group property \eqref{eq:semigroup_real} implies
\[
(\lambda_{t+s},\mu_{t+s},\nu_{t+s}) = (\lambda_t,\mu_t,\nu_t) \subsetplus (\lambda_s,\mu_s,\nu_s) = (\lambda_s,\mu_s,\nu_s) \subsetplus (\lambda_t,\mu_t,\nu_t)
\]
and hence, together with the symmetry of c-free convolution,  
\[
\lambda_{t+s} = \cfconv{\lambda_t}{\nu_t}{\mu_s}{\lambda_s} = \cfconv{\lambda_s}{\nu_s}{\mu_t}{\lambda_t}  = \cfconv{\lambda_t}{\mu_t}{\nu_s}{\lambda_s}. 
\]
Finally,  formula (\ref{eq420}) follows  from $(B_\mu \circ F_{\nu_t}) \times (\ref{eq41}) - (C_\nu \circ F_{\mu_t}) \times (\ref{eq42})$. 
\end{proof}

\begin{cor} \label{cor:diff_eq}
The following differential equations hold. 
\begin{align}
&\frac{\partial F_{\mu_t}}{\partial t} = (B_\mu \circ F_{\nu_t})\cdot \frac{\partial F_{\mu_t}}{\partial z},   \label{eq041} \\ 
&\frac{\partial F_{\nu_t}}{\partial t} =  (C_\nu \circ F_{\mu_t})\cdot \frac{\partial F_{\nu_t}}{\partial z},  \label{eq042} \\
&\frac{\partial F_{\mu_t}}{\partial t} = B_{\mu} \circ F_{\mu_t} + (C_\nu \circ F_{\mu_t})\cdot \left[ \frac{\partial F_{\mu_t}}{\partial z} -1\right], \label{eq0042} \\ 
&\frac{\partial F_{\nu_t}}{\partial t} = C_\nu \circ F_{\nu_t}  + (B_\mu \circ F_{\nu_t})\cdot \left[ \frac{\partial F_{\nu_t}}{\partial z}-1\right],  \label{eq0043} \\
&\frac{\partial F_{\mu_t}}{\partial t} = \frac{(B_\mu \circ F_{\nu_t}) ( B_\mu \circ F_{\mu_t} - C_\nu \circ F_{\mu_t})}{B_\mu \circ F_{\nu_t} - C_\nu \circ F_{\mu_t}}, \label{eq0044} \\
&\frac{\partial F_{\nu_t}}{\partial t} = \frac{(C_\nu \circ F_{\mu_t})(B_\mu \circ F_{\nu_t}- C_\nu \circ F_{\nu_t})}{B_\mu \circ F_{\nu_t} - C_\nu \circ F_{\mu_t}}.  \label{eq0045} 
\end{align}
Note that \eqref{eq0044} and \eqref{eq0045} are valid only when the denominator is not zero. 
\end{cor}
\begin{proof} 
Observe that if $\{(\mu_t, \nu_t)\}_{t\in\real}$ is a one-parameter $\BGP$-free convolution group, both $\{(\mu_t, \mu_t, \nu_t)\}_{t \in \real}$ and $\{(\nu_t, \mu_t, \nu_t)\}_{t \in \real}$ become one-parameter $\AP$-free convolution groups. Formulas (\ref{eq041}), (\ref{eq042}), (\ref{eq0042}), (\ref{eq0043}), (\ref{eq0044}) and (\ref{eq0045}) follow by selecting $\lambda_t = \mu_t$ in (\ref{eq42}), $\lambda_t = \nu_t$ in (\ref{eq41}), $\lambda_t = \mu_t$ in (\ref{eq41}), $\lambda_t = \nu_t$ in (\ref{eq42}),  $\lambda_t = \mu_t$ in (\ref{eq420}) and $\lambda_t = \nu_t$ in (\ref{eq420}), respectively. 
\end{proof}

\begin{rem}
The differential equations (\ref{eq41}) and (\ref{eq42}) are identical if $\lambda_t = \mu_t = \nu_t$ for all $t \in \real$.  
In this case, $\{\mu_t\}_{t\in\real}$ is a one-parameter free convolution group and (\ref{eq41}) is reduced to the generalized complex Burgers equation 
\[
\frac{\partial F_{\mu_t}}{\partial t}(z) = B_\mu (F_{\mu_t}(z)) \frac{\partial F_{\mu_t}}{\partial z}(z) 
\]
derived in \cite{V2}. On the other hand, if $\lambda_t = \mu_t$ and $\nu_t = \delta_0$ for all $t\in\real$, then $\{\mu_t\}_{t\in\real}$ is a one-parameter monotone convolution group.   The differential equations (\ref{eq41}) and (\ref{eq42}) are not of the same form (but eventually equivalent) and they get reduced to 
\[
\frac{\partial F_{\mu_t}}{\partial t}(z) = B_\mu ( F_{\mu_t}(z)) \quad \text{and} \quad \frac{\partial F_{\mu_t}}{\partial t}(z) = B_\mu (z) \frac{\partial F_{\mu_t}}{\partial z}(z),  
\]
which appeared in the monotone case \cite{Mur3}. 
\end{rem}

%%%%%%%%%%%%%%%%%%%%%%%%%%%%%%%%%%%%%%%%%%%%%%%%%%%%%%%%%%%
\subsection{The central limit theorem}\label{subsec:CLT} 
As an application of $\AP$-free cumulants, we prove the central limit theorem for $\AP$-free, identically distributed random variables. Here we introduce a notion of convergence of distributions:  a sequence of distributions $\{(\lambda^n, \mu^n,\nu^n)\}_{n\ge1}$ is said to converge to a distribution $(\lambda, \mu, \nu)$ if for every $k \in \N$ we have 
\[
\lim_{n\to\infty}(\lambda^n(x^k), \mu^n(x^k),\nu^n(x^k)) = (\lambda(x^k), \mu(x^k),\nu(x^k)). 
\]

\begin{prop}\label{thm:CLT}
Let $(\cA, \vp, \psi, \theta)$ be a unital 3-algebraic probability space. 
Let $\{a_i \}_{i = 1} ^{\infty}$ be identically distributed $\AP$-free elements in $\cA$ such that $\vp(a_i) = \psi(a_i) = \theta(a_i) = 0$. Set $ \vara := \vp(a_i ^2)$, $\varb:= \psi(a_i^2)$, $\varc:=\theta(a_i^2)$ and 
\[
s_n:= \frac{a_1 + a_2+\cdots + a_n}{\sqrt{n}}. 
\]
Then, as $n\to\infty$,  the distribution of $s_n$ with respect to $(\vp, \psi, \theta)$ converges to the triplet $(\lambda^{(u,v,w)}, \mu^{(v,w)}, \nu^{(v,w)})$ determined by 
\begin{align}
&\lambda^{(u,v,w)}(x^{2p-1} ) =0 \quad \text{and} \quad  \lambda^{(u,v,w)}(x^{2p} ) = \sum_{\pi \in \ONCP_{\rm pair} (2p)} \frac{1}{p!} \vara^{\# \Out(\pi) } \varb^{\# \AA(\pi) } \varc^{\#\BB(\pi)},  \qquad p\in \N, \\
& \mu^{(v,w)} = \lambda^{(v,v,w)} \quad \text{and} \quad \nu^{(v,w)} = \lambda^{(w,v,w)},  
\label{eq:Kesten_moment}
 \end{align}
where $\ONCP_{\rm pair} (n) := \{\pi \in \ONCP(n): \text{every block of $\pi$ has cardinality two}\}$. 
\end{prop}

\begin{proof}%[Proof of Theorem \ref{thm:CLT}] 
We simply denote $(\lambda, \mu, \nu) = (\lambda^{(u,v,w)}, \mu^{(v,w)}, \nu^{(v,w)}).$  Note first that cumulants have the extensivity and homogeneity: 
\begin{align}
&K_p^X(a_1 + \cdots + a_n) = n K_p^X (a_1), \label{eq:extensivity} \\ 
& K_p^X(\lambda a_1) = \lambda^p K_p^X(a_1)
\end{align}
for all $p,n \in \N, \lambda \in \comp$ and $X \in \{\AP, \BP, \GP\}$. Also, the first cumulants are means and the second cumulants are variances: 
\begin{align}
& K_1^{\AP}(a) = \vp(a), && K_2^{\AP}(a) = \vp(a^2) -[\vp(a)]^2, \\
& K_1^{\BP}(a) = \psi(a), && K_2^{\BP}(a) = \psi(a^2) -[\psi(a)]^2,\\
& K_1^{\GP}(a) = \theta(a), && K_2^{\GP}(a) = \theta(a^2) -[\theta(a)]^2  \label{eq:variance}
\end{align}
for all $a \in \cA.$ 

Combining \eqref{eq:extensivity} -- \eqref{eq:variance} together readily yields 
\[
\lim_{n\to\infty}  (K_p^{\AP}(s_n), K_p^{\BP}(s_n), K_p^{\GP}(s_n) )  = 
\begin{cases}
 (\vara, \varb, \varc), & \text{if~} p=2, \\
(0,0,0), & \text{if~} p =1 \text{~or~} p\ge3.  
\end{cases}
\]
Because every moment is a polynomial on cumulants and vice versa, this implies that for every $p\in\N$ the three moments $\vp(s_n^p)$, $\psi(s_n^p)$ and $\theta(s_n^p)$ respectively converge to $\lambda(x^p), \mu(x^p), \nu(x^p)$, where $(\lambda, \mu, \nu)$ is characterized by its cumulants 
\[
 (K_p^{\AP(\lambda, \mu, \nu)}(x), K_p^{\BP(\mu, \nu)}(x), K_p^{\GP(\mu, \nu)}(x) )  = 
\begin{cases}
 (\vara, \varb, \varc), & \text{if~} p=2, \\
(0,0,0), & \text{if~} p =1 \text{~or~} p\ge3.  
\end{cases}
\]
The desired formula for $\lambda(x^{n})$ is exactly the first moment-cumulant formula in Theorem \ref{thm:m-c-alpha}, in which $(\cA,\vp,\psi,\theta)$ is chosen to be $(\comp[x], \lambda,\mu,\nu)$ and $a_1,\dots, a_n$ are all set to be $x$. Formulas for $\mu$ and $\nu$ follow from Proposition \ref{prop:useful}, the principle used several times so far. 
\end{proof}

With the help of generating functions, we further clarify the limit distributions. 

\begin{thm}\label{prop:CLT} Suppose that $v+w\ne0$. The limit distribution $(\lambda, \mu, \nu)=(\lambda^{(u,v,w)}, \mu^{(v,w)}, \nu^{(v,w)})$ in Proposition \ref{thm:CLT} is characterized by the reciprocal Cauchy transforms
\begin{gather}
F_\lambda (z) = \left( 1 - \frac{\vara}{\varb + \varc}\right) z + \frac{\vara z}{\varb + \varc}\sqrt{1 - \frac{2 (\varb + \varc)}{z^2}}, \label{eq00041}\\
F_{\mu}(z) = \left( 1 - \frac{\varb}{\varb + \varc}\right) z + \frac{\varb z }{\varb + \varc}\sqrt{1 - \frac{2 (\varb + \varc)}{z^2}}, \label{eq00042} \\
F_{\nu}(z) = \left( 1 - \frac{\varc}{\varb + \varc}\right) z + \frac{\varc z }{\varb + \varc}\sqrt{1 - \frac{2 (\varb + \varc)}{z^2}},  \label{eq00043}
\end{gather}
where $\sqrt{1 + w}$ is to be interpreted as the formal power series on $w$. 
\end{thm}

\begin{rem} 
The singularity $v+w=0$ turns out to be removable and therefore formulas \eqref{eq00041} -- \eqref{eq00043} can be extended to the case $v+w=0$.  
\end{rem}

\begin{rem} The limit distributions $\lambda,\mu,\nu$ are called the Kesten distributions if $u, v,w \ge0$ which is the case when $(\cA,\vp,\psi,\theta)$ has the positivity structure and $a_i$ are self-adjoint elements. 
The Kesten distribution already appeared in the central limit theorem for c-freeness \cite{BLS} and c-monotone independence \cite{Has3}.  
\end{rem}

\begin{proof}[Proof of Theorem \ref{prop:CLT}] \

\vspace{2mm}
\noindent
{\bf Case 1: $\varb \ne \varc$ and $\varc\ne0$.}  
Let $(\lambda^n,\mu^n,\nu^n)$ be the distribution of $s_n$ and $(\lambda^n_t,\mu^n_t,\nu^n_t), t\in \real$ be the associated one-parameter $\AP$-free convolution group. Also let $(\lambda_t,\mu_t,\nu_t), t\in \real$ be the one-parameter $\AP$-free convolution group associated with $(\lambda,\mu,\nu)$.  Then for every $t\in \real$ and $p\in\N$, 
\[
\lim_{n\to \infty} \lambda_t^n(x^p) = \lambda_t(x^p),  \quad \lim_{n\to \infty} \mu_t^n(x^p) = \mu_t(x^p) \quad \text{and} \quad \lim_{n\to \infty} \nu_t^n(x^p) = \nu_t(x^p) 
\]
because the cumulants of $(\lambda^n,\mu^n,\nu^n)$ converge to those of $(\lambda,\mu,\nu)$ and the moments of $(\lambda^n_t,\mu^n_t,\nu^n_t)$ are polynomials on the cumulants of $(\lambda^n,\mu^n,\nu^n)$ and $t$. This implies that the differential equations in Proposition \ref{prop:diff_eq} and in Corollary \ref{cor:diff_eq} for $\{(\lambda^n_t,\mu^n_t,\nu^n_t)\}_{t\in \real}$ pass to the differential equations for $\{(\lambda_t,\mu_t,\nu_t)\}_{t\in \real}$ in the limit. Therefore, the limit distribution $(\lambda,\mu,\nu)$ satisfies the differential equations in Proposition \ref{prop:diff_eq} and in Corollary \ref{cor:diff_eq}, where
\[
A_\lambda(z) = -\frac{\vara}{z},\quad  B_\mu(z) = -\frac{\varb}{z} \quad \text{and} \quad C_\nu(z) = -\frac{\varc}{z}.
\]
Assume for the moment that $\varb \neq \varc$ and $\varc\ne0$.  
Equations \eqref{eq0044} and (\ref{eq0045}) then make sense because the denominator of the RHS is 
\[
-\frac{\varb}{F_{\nu_t}(z)} + \frac{\varc}{F_{\mu_t}(z)} = -\varb G_{\nu_t}(z)  + \varc G_{\mu_t}(z) = \frac{\varc -\varb}{z} + \left(\text{higher order terms on $\frac1{z}$} \right),  
\]
 which is nonzero as a formal Laurent series. 
After some calculations equations \eqref{eq0044} and (\ref{eq0045}) take the form
\begin{align}
&\frac{\partial F_{\mu_t}}{\partial t}(z) =  \frac{\varb(\varb - \varc)}{-\varb F_{\mu_t}(z) + \varc F_{\nu_t}(z)} \label{eq:mu},   \\ 
&\frac{\partial F_{\nu_t}}{\partial t}(z) =  \frac{\varc(\varb - \varc)}{-\varb F_{\mu_t}(z) + \varc F_{\nu_t}(z)}.  \label{eq:nu}
\end{align}
Therefore, the identity $\varc \frac{\partial F_{\mu_t}}{\partial t}(z) =\varb \frac{\partial F_{\nu_t}}{\partial t}(z)$ follows, which implies that 
\begin{equation} \label{eq:mu_nu} 
F_{\mu_t} (z) = s F_{\nu_t}(z) + (1 - s)z,  \quad \text{where} \quad s = \frac{\varb}{\varc}. 
\end{equation}
Substituting \eqref{eq:mu_nu} into \eqref{eq:mu} and \eqref{eq:nu} we obtain 
\begin{gather}
F_{\mu_t}(z) = \Big( 1 - \frac{s}{ 1 + s}\Big) z + \frac{s z}{1 + s} \sqrt{1 - \frac{2 \varc(1 + s)t}{z^2}}, \\
F_{\nu_t}(z) = \Big( 1 - \frac{1}{1 + s}\Big) z + \frac{z}{1 + s}\sqrt{1 -\frac{2 \varc(1 + s)t}{z^2}}. 
\end{gather}
Because the limit distribution $(\mu, \nu)$ equals $(\mu_1, \nu_1)$, we conclude the desired \eqref{eq00042} and \eqref{eq00043}. 

The remaining  $\lambda$ can be calculated from equation (\ref{eq420}) which now gets reduced to  
\[
\frac{\partial F_{\lambda_t}}{\partial t}(z) =  \frac{\vara(\varb - \varc)}{-\varb F_{\mu_t}(z) + \varc F_{\nu_t}(z)}. 
\]  
By simple calculations we obtain the remaining formula \eqref{eq00041}.  

\vspace{2mm} 
\noindent
{\bf Case 2: the general case.}  According to Proposition \ref{thm:CLT}, the moment $\lambda(x^p)$ is a polynomial on $\vara,\varb,\varc$. 
On the other hand, formula \eqref{eq00041} holds in case $\varb \ne \varc, \varc\ne0$, which provides (by looking at the Cauchy transform $G_\lambda$ and its coefficients) another expression of $\lambda(x^p)$ as a meromorphic function (actually, polynomial) of two variables $u$ and $v+w$. Because these two expressions coincide in case $\varb \ne \varc, \varc\ne0$,  they still coincide for any $u,v,w \in \comp$ (with $v+w\ne0$) by the identity theorem. The other two formulas \eqref{eq00042} and \eqref{eq00043} can be similarly proved. 
\end{proof}

\begin{rem}\label{eq:Wick_formula} 
According to \cite[p.~5]{Sza22+}, there is an alternative formula 
\begin{equation}
\lambda(x^{2p} ) = \sum_{i=0}^p b_{p,i}\left(\frac{\varb+\varc}{2}\right)^i \left( \vara - \frac{\varb+\varc}{2}\right)^{p-i},   
\end{equation}
where $b_{p,i} := \binom{2p}{i} - \binom{2p}{i-1}, i\ge1; b_{p,0}:=1$.  
\end{rem}

%%%%%%%%%%%%%%%%%%%%%%%%%%%%%%%%%%%%
%%%%%%%% Bibliography  %%%%%%%%%%%%%%%%%%%%%
\providecommand{\bysame}{\leavevmode\hbox to3em{\hrulefill}\thinspace}

   \vspace{10mm}
 \noindent  
 \author{Department of Mathematics \\ Hokkaido University \\ Kita 10, Nishi 8, Kita-Ku, Sapporo \\
  Hokkaido, 060-0810, Japan \\ E-mail: thasebe@math.sci.hokudai.ac.jp}

}   %%%% end of "\large" 
\end{document}